\documentclass[11pt]{article}

\usepackage{amscd}      % to be able to use "CD" environment
\usepackage{amssymb}
\usepackage[intlimits]{amsmath}
\usepackage{xypic}      % to be able to use "diagram" environment of
\LaTeXdiagrams          % xypic to typeset commutative diagrams
\usepackage[all,v2]{xy}
\xyoption{2cell}
\UseAllTwocells
\xyoption{frame}
\CompileMatrices

\usepackage{theorem}

\newtheorem{prop}{Proposition}[section]
\newtheorem{lem}[prop]{Lemma}

\newtheorem{cor}[prop]{Corollary}
\newtheorem{them}[prop]{Theorem}

\theorembodyfont{\upshape}

\newtheorem{defn}[prop]{Definition}

\newtheorem{rmk}[prop]{Remark}

\newtheorem{numex}[prop]{Example}

\newenvironment{pf}{\begin{trivlist}\item[]{\sc Proof.}}%
            {\nolinebreak $\Box$ \end{trivlist}}
            {\nolinebreak $\Box$ \end{trivlist}}

\newcommand{\noprint}[1]{}

\newcommand{\teps}{{\tilde{\varepsilon}}}
\newcommand{\tepss}{{\tilde{\varepsilon'}}}

\renewcommand{\tilde}{\widetilde}

\newcommand{\toto}{\rightrightarrows}

\newcommand{\com}{{\scriptscriptstyle\bullet}}
\newcommand{\lcom}{_{\scriptscriptstyle\bullet}}
\newcommand{\lcomm}{_{\scriptscriptstyle\bullet\bullet}}
\newcommand{\upcom}{^{\scriptscriptstyle\bullet}}

\newcommand{\Gg}{{\mathfrak g}}

\newcommand{\zz}{{\mathbb Z}}
\newcommand{\hh}{{\mathbb H}}

\newcommand{\nn}{{\mathbb N}}

\newcommand{\cc}{{\mathbb C}}
\newcommand{\rr}{{\mathbb R}}

\newcommand{\ff}{{\mathbb F}}

\newcommand{\aA}{{\cal A}}
\newcommand{\bB}{{\cal B}}

\newcommand{\cC}{{\cal C}}
\newcommand{\eE}{{\cal E}}
\newcommand{\iI}{{\cal I}}

\newcommand{\hH}{{\cal H}}
\newcommand{\sS}{{\cal S}}

\newcommand{\uU}{{\cal U}}

\newcommand{\rR}{{\cal R}}

\newcommand{\del}{\partial}

\newcommand{\Hom}{\mathop{\rm Hom}\nolimits}

\newcommand{\smalcirc}{\mbox{\tiny{$\circ $}}}

\newcommand{\ldiag}[1]%
       {\makebox[0cm]{${\scriptstyle#1}\downarrow\phantom{\scriptstyle#1}$}}
\newcommand{\ldiagup}[1]%
       {\makebox[0cm]{${\scriptstyle#1}\uparrow\phantom{\scriptstyle#1}$}}
\newcommand{\rdiag}[1]%
       {\makebox[0cm]{$\phantom{\scriptstyle#1}\downarrow{\scriptstyle#1}$}}
\newcommand{\sediagr}[1]%
       {\makebox[0cm]{$\phantom{\scriptstyle#1}\searrow{\scriptstyle#1}$}}
\newcommand{\nediagr}[1]%
       {\makebox[0cm]{$\phantom{\scriptstyle#1}\nearrow{\scriptstyle#1}$}}
\newcommand{\rdiagup}[1]%
       {\makebox[0cm]{$\phantom{\scriptstyle#1}\uparrow{\scriptstyle#1}$}}
\newcommand{\swdiag}[1]%
       {\makebox[0cm]{$\phantom{\scriptstyle#1}\swarrow{\scriptstyle#1}$}}
\newcommand{\sediag}[1]%
       {\makebox[0cm]{${\scriptstyle#1}\searrow\phantom{\scriptstyle#1}$}}
\newcommand{\nediag}[1]%
       {\makebox[0cm]{${\scriptstyle#1}\nearrow\phantom{\scriptstyle#1}$}}

\newcommand{\longiso}{\stackrel{\textstyle\sim}{\longrightarrow}}
\newcommand{\iso}{\stackrel{\sim}{\rightarrow}}

\newcommand{\doublearrowstack}[2]%
 {{{{\scriptstyle#1}\atop{\textstyle\longrightarrow}}\atop{{\textstyle\longright
arrow}\atop{\scriptstyle#2}}}}
\newcommand{\rightleftarrowstack}[2]%
 {{{{\scriptstyle#1}\atop{\textstyle\longrightarrow}}\atop{{\textstyle\longlefta
rrow}\atop{\scriptstyle#2}}}}
\newcommand{\leftrightarrowstack}[2]%
 {{{{\scriptstyle#1}\atop{\textstyle\longleftarrow}}\atop{{\textstyle\longrighta
rrow}\atop{\scriptstyle#2}}}}

\newcommand{\overtoparrow}%
{\makebox[0cm]{\beginpicture
\setcoordinatesystem units <.8cm,.4cm> point at 0 0
\setplotarea x from -3 to 3, y from 0 to 1
\setquadratic
\plot -3 0 0 1 3 0 /
\put{\vector(3,-1){0}}[Bl] at 3 0
\endpicture}}

\newcommand{\underbottomarrow}%
{\makebox[0cm]{\beginpicture
\setcoordinatesystem units <.8cm,.4cm> point at 0 0
\setplotarea x from -3 to 3, y from 0 to 1
\setquadratic
\plot -3 1 0 0 3 1 /
\put{\vector(3,1){0}}[Bl] at 3 1
\endpicture}}

\newcommand{\ses}[5]%
{0\longrightarrow#1\stackrel{#2}{ \longrightarrow}#3\stackrel{#4}{
\longrightarrow}#5\longrightarrow0}

\newcommand{\dt}[6]%
{#1\stackrel{#2}{\longrightarrow}#3 \stackrel{#4}{\longrightarrow}#5
\stackrel{#6}{\longrightarrow} #1[1]}

\newcommand{\cat}[1]%
{(\mbox{\rm #1})}

\newcommand{\gm}{\Gamma}
\newcommand{\tgm}{{\tilde{\Gamma}}}

\newcommand{\lon }{\longrightarrow }

\newcommand{\be }{\begin{eqnarray*}}
\newcommand{\ee }{\end{eqnarray*}}

\newcommand{\DG}{generalized simplicial}
\newcommand{\M}[1]{M_{#1,\bullet}}
\newcommand{\Ma}[1]{M_{\bullet, #1}}
\newcommand{\Mr}[1]{M'_{#1,\bullet}}
\newcommand{\otimesc}{\otimes_{C_0(\gm_0 )}}
\newcommand{\deltaa}{\lambda}

\def\gpd{\,\lower1pt\hbox{$\longrightarrow$}\hskip-.24in\raise2pt
             \hbox{$\longrightarrow$}\,}

\title{The ring structure for equivariant twisted $K$-theory}
\author{ Jean-Louis Tu\\
Universit\'e Paul Verlaine -- Metz\\
LMAM - CNRS UMR 7122\\
 ISGMP, B\^atiment A, Ile du Saulcy\\
 57000 Metz, France\\
 {\sf email: tu@univ-metz.fr}
 \\\\
Ping Xu \thanks{Research partially supported by NSF
       grants DMS03-06665 and DMS-0605725 \& NSA grant H98230-06-1-0047.}\\
        Department of Mathematics\\
         Pennsylvania State University \\
         University Park, PA 16802, USA\\
{\sf email: ping@math.psu.edu }}

\date{}
\begin{document}
\sloppy
\maketitle

\begin{abstract}
We prove, under some mild conditions, that the equivariant twisted $K$-theory
group of a crossed module admits a ring structure if the twisting $2$-cocycle is
$2$-multiplicative. We also give an explicit construction of the transgression
map $T_1: H^{*}(\Gamma\lcom;\aA) \to H^{*-1}((N\rtimes \gm)\lcom;\aA)$ for any
crossed module $N\to \gm$ and prove that any element in the image is
$\infty$-multiplicative. As a consequence, we prove that, under some mild
conditions, for a crossed module $N \to \gm$ and any $e \in
\check{Z}^3(\gm\lcom;\sS^1)$, that the equivariant twisted $K$-theory group
${K}^*_{e,\gm}(N)$ admits a ring structure. As an application, we prove that for
a compact, connected and simply connected Lie group $G$, the equivariant twisted
$K$-theory group $K_{[c], G}^* (G)$, defined as
the $K$-theory group of a  certain groupoid $C^*$-algebra,
 is endowed with a canonical ring structure
$K^{i+d}_{[c],G}(G)\otimes K^{j+d}_{[c],G}(G)\to K^{i+j+d}_{[c], G}(G)$,
where $d=\dim G$ and $[c]\in H^2((G\rtimes G)\lcom;\sS^1)$.
The relation with Freed-Hopkins-Teleman theorem \cite{FHT0} still needs to
be explored.
\end{abstract}

{\small \tableofcontents}

\section{Introduction}

A great deal of interest in twisted equivariant $K$-theory has emerged due to
its close connection to string theory \cite{Witten, Witten1}. In particular, the
recent work of Freed--Hopkins--Teleman \cite{F, FHT02, FHT0, FHT, FHT1}
 concerning the relationship between the twisted equivariant
 $K$-theory of compact Lie groups
and Verlinde algebras has inspired a great deal of activities in this subject.
%%%%%%
It now  becomes increasingly important to develop a general framework
 which allows one to study the ring structure of 
twisted equivariant $K$-theory groups and in particular
to investigate the general criteria
which guarantee the existence of such a ring structure.
%%%%%%

%%%%
This paper serves  this purpose.
%%%%
 More precisely, in this paper we examine the
conditions under which the twisted $K$-theory groups of a crossed module admit a
ring structure. Recall that a crossed module is a groupoid morphism
$$\xymatrix{
N_1\ar[r]^\varphi\ar[d]\ar@<-.6ex>[d] & \gm_1\ar@<-.6ex>[d]\ar[d] \\
N_0\ar[r]^{=} & \gm_0
}$$
where $N_1\toto N_0$ is a bundle of groups, together with an action of $\gm$ on
$N$ by automorphisms satisfying some compatibility conditions (see Definition
\ref{def:cross}). A standard example of a crossed module is as follows. Let
$\gm_1\toto \gm_0$ be a groupoid and $S\gm_1=\{g\in \gm_1\vert\; s(g)=t(g)\}$ be
the space of closed loops in $\gm_1$. Then the canonical inclusion $S\gm \to
\gm$, together with the conjugation action of $\gm$ on $S\gm$, forms a
crossed module. In particular, when $\gm$ is just a Lie group $G$, $S\gm$ is
isomorphic to $G$ with the action being by conjugation. In other words,
$G\stackrel{id}{\to}G$ with the conjugation action is a crossed module. Given a
crossed module $N\to \gm$, since $\gm$ acts on $N$, one forms the transformation
groupoid (also called the crossed product groupoid) $N\rtimes \gm$. In the case
that the crossed module is $S\gm\to \gm$, the transformation groupoid obtained
is called the inertia groupoid and is denoted by $\Lambda\gm$. When $\gm$ is a
Lie group $G$, the inertia groupoid is the standard transformation groupoid
$G\rtimes G\toto G$ with $G$ acting on $G$ by conjugation.
 
In \cite{TXL04}, we developed a general theory of twisted $K$-theory for
differential stacks
%%%%% 
(see also \cite{Ati-Seg, Ati-Seg2} for the case of quotient stacks).
%%%%%
 For a Lie groupoid $X_1\toto X_0$ and $\alpha\in
H^2(X\lcom, \sS^1)$, the twisted $K$-theory groups $K_\alpha^* (X)$ are defined
to be the $K$-theory groups of a certain $C^*$-algebra $C^*_r (X, \alpha )$
associated to the element $\alpha$ (or an $S^1$-gerbe) using groupoid 
$S^1$-central
extensions. However, the construction is not canonical and depends on a choice of
$2$-cocycle $c\in \check{Z}^2(X \lcom;\sS^1)$ representing $\alpha$, though
different choices of $c$ give rise to isomorphic $K$-theory groups. For the
convenience of our investigation, in this paper, we will define twisted
$K$-theory groups using a \v{C}ech $2$-cocycle instead of a cohomology class so
that the twisted $K$-theory groups $K_c^* (X)$ will be canonically defined. For
a Lie groupoid $\gm$ acting on a manifold $N$, and $c\in \check{Z}^2((N\rtimes
\gm)\lcom ;\sS^1)$ a $2$-cocycle of the corresponding transformation groupoid
$N\rtimes \gm$, the twisted equivariant $K$-theory groups are then defined to be
$$K^i_{c, \gm} (N)=K^i_c (N \rtimes \gm).$$
The main question we study in this paper is: For a crossed module $N\to \gm$,
under what condition do the twisted equivariant $K$-theory groups
$K^i_{c,\gm}(N)$ admit a ring structure?

The answer is that $c$ needs to be $2$-multiplicative. Note that since $N\to
\gm$ is a crossed module, $(N\lcom\rtimes \gm)\lcom$ becomes a bi-simplicial
space. Therefore there are two simplicial maps $\del: \check{C}^p((N_q\rtimes
\gm)\lcom;\sS^1) \to \check{C}^{p+1}((N_q\rtimes \gm)\lcom;\sS^1)$ and $\del':
\check{C}^p((N_q\rtimes \gm)\lcom;\sS^1) \to \check{C}^{p}((N_{q+1}\rtimes
\gm)\lcom;\sS^1)$. A $2$-cocycle $c\in \check{Z}^2((N_1\rtimes \gm)\lcom;\sS^1)$
(i.e. $\del c=0$) is said to be $2$-multiplicative if there exist $b\in
\check{C}^1((N_2\rtimes \gm)\lcom;\sS^1)$ and $a\in \check{C}^0((N_3\rtimes
\gm)\lcom;\sS^1)$ such that $\del'c = \del b$, and $\del' b = \del a$. Such a
triple $(c, b, a)$ is called a \emph{multiplicator}. The product structure on
$K^*_{c, \gm} (N)$ depends on the choice of a multiplicator. The main result of
the paper can be summarized as the following

\vskip.2in

\noindent{\bf Theorem A.} {\it
Let $N\stackrel{\varphi}{\to} \gm$ be a crossed module, where $\gm_1\toto \gm_0$
is a proper Lie groupoid such that $s:N_1\to N_0$ is $\gm$-equivariantly
$K$-oriented. Assume that $(c, b, a)$ is a multiplicator, where $c\in
\check{C}^2((N_1\rtimes \gm)\lcom, \sS^1)$, $b\in \check{C}^1((N_2\rtimes
\gm)\lcom, \sS^1)$, and $a\in \check{C}^0((N_3\rtimes \gm)\lcom, \sS^1)$. Then
there is a canonical associative product
$$K^{i+d}_{c,\gm}(N)\otimes K^{j+d}_{c,\gm}(N)\to K^{i+j+d}_{c,\gm}(N),$$
where $d=\dim N_1 -\dim N_0$.}

\vskip.2in

%%%%%
%PING: the following paragraph has been moved to the end of introduction

%Note that the idea of using  multiplicative cocycles
%(called equivariantly primitive in \cite{FHT02}) in constructing
%the product on twisted equivariant $K$-theory has been known
% in the community (see \cite{FHT02, CJMSW} for instance). However,
%it seems that the condition of $2$-multicativity  is new, which is
%very essential for our proof of the associativity of the product constructed.
%%%%% 

The main idea of our approach is to transform this geometric problem into a
problem of $C^*$-algebras, for which there are many sophisticated $K$-theoretic
techniques. As the first step, we give a canonical construction of an equivariant
$S^1$-gerbe (or rather $S^1$-central extension), which should be
of independent interest.

\vskip.2in

\noindent{\bf Theorem B.} {\it
Suppose that $\gm: \gm_1\toto \gm_0$ is a Lie groupoid acting on a manifold $N$
via $J:N\to \gm_0$. Let $\uU$ be a cover of $(N\rtimes\gm)\lcom$. Then any
\v{C}ech $2$-cocycle $c\in \check{Z}^2(\uU, \sS^1 )$ determines a canonical
$S^1$-central extension of the form $\tilde{H}\rtimes \gm\to H\rtimes \gm\toto
M$, where $\tilde{H}\to H\toto M$ is a $\gm$-equivariant  $S^1$-central
extension and $H\toto M$ is Morita equivalent to $N\toto N$, with the
Dixmier-Douady class of
the  extension equal to $[c]\in\check{H}^2((N\rtimes\gm)\lcom;\sS^1)$.}

\vskip.2in

The above theorem allows us to establish a canonical Morita equivalence between
the $C^*$-algebra $C^*_r (N\rtimes \gm, c)$ and the crossed product algebra $A_c
\rtimes_r \gm$, where $A_c$ is a $\gm$-$C^*$-algebra (i.e. a $C^*$-algebra with
a $\gm$-action). This enables us to construct the product structure on
$K^*_{c,\gm}(N)$ with the help of the Gysin map and the external Kasparov
product.

For a $\gm$-equivariantly $K$-oriented submersion $f: M \to N$ between proper
$\gm$-manifolds $M$ and $N$, the Gysin map is a wrong-way functorial map
$$f_!:K^i_{f^*c, \gm}(M)\to K^{i+d}_{c, \gm} (N),$$
where $d= \dim N-\dim M$, which satisfies $g_!\smalcirc f_! = (g\smalcirc f)_!$.
It is standard that any $K$-oriented map $f: M\to N$ yields a Gysin element $f_!
\in KK^d(C_0(M),C_0(N))$ \cite{CS84,HS87}. When $\gm$ is a Lie group, an
equivariant version was proved by Kasparov--Skandalis \cite[\S 4.3]{KS91}: Any
$\gm$-equivariantly $K$-oriented map $f: M \to N$ determines an element $f_! \in
KK^d_\gm(C_0(M),C_0(N))$. A similar argument can be adapted to show that the
same assertion holds when $\gm$ is a Lie groupoid and $KK^*_\gm$ is Le Gall's
groupoid equivariant $KK$-theory \cite{Leg}. As a consequence, our Gysin map can
easily be constructed using such a Gysin element. We note that a different
approach to the Gysin map to (non-equivariant) twisted $K$-theory was recently
studied by Carey--Wang \cite{CW}.

The second ingredient of our construction is the external Kasparov product
\begin{equation}
K^i_{c,\gm}(N)\otimes K^j_{c,\gm}(N) \to K^{i+j}_{p_1^*c+p_2^*c,\gm}(N_2),
\end{equation}
where $\gm$ is a proper Lie groupoid, and $p_1,\ p_2: N_2\to N_1$ are the
natural projections. This essentially follows from the usual Kasparov product
$KK^i_\gm(A,B)\otimes KK^j_\gm(C,D)\to KK^{i+j}_\gm (A\otimes_{C_0 (\gm_0)} C,
B\otimes_{C_0 (\gm_0)} D)$, where $A$, $B$, $C$, $D$ are $\gm$-$C^*$-algebras.
Here again $KK^*_\gm$ stands for the Le Gall's groupoid version of the
equivariant $KK$-theory of Kasparov \cite{Kas88, Leg}.

Theorem A indicates that the ring structure on twisted equivariant $K$-theory
groups relies on ``multiplicators". A natural question now is how multiplicators
arise. In the first half part of the paper, we discuss an important
construction, the so-called transgression maps, which is a powerful 
tool to
produce ``multiplicators". At the level of cohomology, the transgression map for
a crossed module $N\to \gm$ is a map
$$T_1: H^{k}(\gm\lcom; \sS^1 )\to H^{k-1}((N_1\rtimes \gm)\lcom; \sS^1). $$
For instance, when $k=2$, one obtains a map $T_1: H^{3}(\gm\lcom; \sS^1 )\to
H^{2}((N\rtimes \gm)\lcom; \sS^1)$. Any element in the image of $T_1$ is
$2$-multiplicative, so it is reasonable to expect that the corresponding twisted
$K$-theory groups admit a ring structure. To prove this assertion, since our
twisted $K$-theory groups are defined in terms of $2$-cocycles, we must study
the transgression map more carefully at the cochain level. Therefore we put our
construction of the transgression map into a more general perspective which we
believe to be of independent interest.

First, to make our construction more transparent and intrinsic, we introduce the
notion of $\cC$-spaces and their sheaf cohomology, for a category $\cC$. By a
$\cC$-space, we mean a contravariant functor from the category $\cC$ to the
category of topological spaces. One similarly defines $\cC$-manifolds. Here we
are mainly interested in $\cC$-spaces in which $\cC$ is equipped with an
additional \DG{} structure. One standard example of a \DG{} category is the
simplicial category $\Delta$, whose corresponding $\cC$-spaces are simplicial
spaces. Indeed the \DG{} structure on $\cC$ enables us to define sheaf and
\v{C}ech cohomology of a $\cC$-space just as one does for simplicial spaces
\cite{Deligne, Friedlander}. A relevant \DG{} category for our purpose here is
the so-called $\Delta_2$-category, which is an extension of the bi-simplicial
category, i.e., $\Delta\times \Delta$. Indeed $\Delta_2$ has the same objects as
$\Delta\times \Delta$, but contains more morphisms.

Let $M\lcomm$ be a $\Delta_2$-space. Then for any fixed $k\in \nn$, both
$M_{k,\bullet}=(M_{k,l})_{l\in\nn}$ and $M_{\bullet, k}=(M_{l,k})_{l\in\nn}$
are simplicial spaces. Suppose that $\aA_0 \stackrel{d}{\to} \aA_1
\stackrel{d}{\to} \cdots$ is a  complex of abelian sheaves over
$M\lcomm$. Let $C^*(M\lcomm;\aA\lcom)$ (resp. $C^*(\M{0};\aA\lcom)$) be its
associated differential complex on $M\lcomm$ (resp. $\M{0}$). We prove the
following
 
\vskip.2in

\noindent{\bf Theorem C.} {\it
\begin{enumerate}
\item For each $k \in \nn$, there is a map
    $$T_k: C^*(\M{0};\aA\lcom)\to C^{*-k}(M_{k,\bullet};\aA\lcom) $$
    (with $T_0=\mbox{Id}$) such that
    $$T=\sum_{k\ge 0} T_k: C^*(\M{0};\aA\lcom)\to C^*(M\lcomm;\aA\lcom) $$
    is a chain map which therefore induces a morphism
    $$T: H^*(\M{0};\aA\lcom)\to H^*(M\lcomm;\aA\lcom )$$
    on the level of cohomology.
\item In particular, 
    $$T_1: C^*(\M{0};\aA\lcom)\to C^{*-1} (\M{1};\aA\lcom) $$
    is an (anti-)chain map and thus induces a morphism
    $$T_1: H^*(\M{0};\aA\lcom)\to H^{*-1}(M_{1,\bullet};\aA\lcom). $$
\item Similarly, given an abelian sheaf $\aA$ over $M\lcomm$, there is a map
    $$T_k: \check{C}^*(\M{0};\aA)\to \check{C}^{*-k}(\M{k};\aA) $$
    (with $T_0=\mbox{Id}$) such that
    $$T=\sum_{k\ge 0} T_k:
        \check{C}^*(\M{0};\aA)\to \check{C}^{*}(M\lcomm;\aA)$$
    is a chain map which therefore induces a morphism
    $$T: \check{H}^*(\M{0};\aA)\to \check{H}^{*}(M\lcomm;\aA).$$
\item Similarly, for any abelian sheaf $\aA$ over $M\lcomm$,
    $$T_1: \check{C}^*(\M{0};\aA)\to \check{C}^{*-1} (\M{1};\aA)$$
    is an (anti-)chain map and thus induces a morphism
    $$T_1: \check{H}^*(\M{0};\aA)\to \check{H}^{*-1}(M_{1,\bullet};\aA). $$
\end{enumerate}
We call $T$ the \emph{total transgression map} and $T_1$ the \emph{transgression
map}.}

\vskip.2in

For a crossed module $N\stackrel{\phi}{\to}\gm$, one shows that $(N\rtimes \gm
)\lcomm$ is naturally a $\Delta_2$-space. In this case, the transgression maps
can be described more explicitly.

\vskip.2in

\noindent{\bf Theorem D.} {\it
Let $N\stackrel{\phi}{\to}\gm$ be a crossed module and $\aA\lcom$ a 
complex of abelian sheaves over $(N\rtimes \gm)\lcomm$. 
\begin{enumerate}
\item There is a chain map (the total transgression map)
    $$T=\sum_k T_k:
        C^*(\gm\lcom;\aA\lcom)\to C^*((N\rtimes \gm)\lcomm;\aA\lcom).$$
    Moreover
    $$T_k=\sum_{\sigma\in S_{k,l}} \varepsilon(\sigma) \tilde{f}_\sigma^*:
        C^{k+l}(\gm\lcom; \aA\lcom) \lon C^l((N_k\rtimes\gm)\lcom; \aA\lcom), $$
    where $S_{k,l}$ denotes the set of $(k,l)$-shuffles, and the map
    $\tilde{f}_\sigma: N_k\rtimes \gm_l \to \gm_{k+l}$ is given by
    \begin{equation}
    \label{eqn:fsigma2}
    \tilde{f}_\sigma(x_1,\ldots,x_k;g_1,\ldots,g_l)=(u_1,\ldots,u_{k+l}),
    \end{equation}
    where $u_i=g_{\sigma^{-1}(i)}$ if $\sigma^{-1}(i)\ge k+1$, and
    $u_i=\varphi\left(x_{\sigma^{-1}(i)}^{\prod_{\sigma^{-1}(j)>k, j< i}
    g_{\sigma^{-1}(j)}}\right)$ otherwise.
\item There is a transgression map
    $$T_1:
        H^{*}(\gm\lcom;\aA\lcom)\to H^{*-1}((N_1\rtimes \gm)\lcom;\aA\lcom), $$
    which is given, on the cochain level, by
    $$T_1=\sum_{i=0}^{p-1} (-1)^i\tilde{f}_i^* :
        \aA^q(\gm_p)\to \aA^q(N_1\rtimes \gm_{p-1}).$$
    Here the map $\tilde{f}_i: N_1\rtimes \gm_{p-1}\to \gm_p$ is given by
    \begin{equation}\label{eqn:T1dR}
    \tilde{f}_i(x;g_1,\ldots,g_{p-1})=(g_1,\ldots,g_i,
    \varphi(x)^{g_1\cdots g_i} ,g_{i+1},\ldots,g_{p-1}).
    \end{equation}
\end{enumerate}
}

\vskip.2in

Note that the transgression maps have, in various different forms, appeared in
the literature before. For instance, for the crossed module $G\stackrel{id}{\to}
G$ with the conjugation action and $\aA\lcom=\Omega^\com$, the transgression map
$T_1: H^*_G(\bullet )\to H^{*-1}_G (G)$ was studied by Jeffrey \cite{Jeffrey} (see
also \cite{Meinrenken}). The geometric meaning of the transgression $T_1:
\Omega^4_G (\bullet  ) \to \Omega^3_G (G)$ was studied by Brylinski--McLaughlin
\cite{BM}. On the other hand, the suspension map
$ H^4_G(\bullet , \zz)\to H^3(G, \zz)$, which is  the composition of
the transgression $T_1: H^4_G(\bullet, \zz)\to H^3_G(G, \zz)$ with
the canonical map $H^3_G(G, \zz)\to H^3(G, \zz)$,
was shown by Dijkgraaf--Witten  \cite{DW} to induce a geometric 
correspondence between three dimensional Chern-Simons functional
 and Wess-Zumino-Witten models. Such a correspondence
 was further explored recently by Carey et. al. \cite{CJMSW}
using bundle gerbes and was also used in \cite{CW1} 
in their study of fusion of $D$-branes.
 The transgression map for orbifold cohomology was recently studied by
Adem--Ruan--Zhang \cite{ARZ}.

Transgression maps $T_k$ can be used to produce multiplicators. More precisely,
for a crossed module $N\stackrel{\varphi}{\to}\gm$, if $e\in
\check{Z}^3(\gm\lcom;\sS^1)$ and letting $c=T_1 e \in \check{C}^2((N_1\rtimes
\gm)\lcom;\sS^1)$, $b=-T_2 e\in \check{C}^1((N_2\rtimes \gm)\lcom;\sS^1)$, and
$a=-T_3e\in \check{C}^0((N_3\rtimes \gm)\lcom;\sS^1)$, we then prove that $(c,
b, a)$ is a multiplicator. This fact enables us to construct a canonical ring
structure on the $K$-theory groups twisted by elements in
$\check{Z}^3( \gm\lcom;\sS^1)$. More precisely, for any $e\in \check{Z}^3(
\gm\lcom;\sS^1)$, $T_1 e\in \check{Z}^2( (N_1\rtimes\gm)\lcom;\sS^1)$ is
$2$-multiplicative. Define
$$K_{e, \gm}^*(N): =K_{T_1 e, \gm}^*(N). $$
Thus we prove

\vskip.2in

\noindent{\bf Theorem E.} {\it
Let $N\stackrel{\varphi}{\to} \gm$ be a crossed module, where $\gm_1\toto \gm_0$
is a proper Lie groupoid such that $s:N_1\to N_0$ is $\gm$-equivariantly
$K$-oriented.
\begin{enumerate}
\item For any $e\in \check{Z}^3( \gm\lcom;\sS^1)$, the twisted $K$-theory group
    ${K}^*_{e,\gm}(N)$ is endowed with a ring structure
    $$K^{i+d}_{e,\gm}(N)\otimes K^{j+d}_{e,\gm}(N)\to K^{i+j+d}_{e,\gm}(N),$$
    where $d=\dim N_1 -\dim N_0$.
\item Assume that $e$ and $e'\in \check{Z}^3(\gm\lcom;\sS^1)$ satisfy
    $e-e'=\del u$ for some $u\in \check{C}^2(\gm\lcom;\sS^1)$. Then there is a
    ring isomorphism
    $$\Psi_{e',u,e}:K^*_{e,\gm}(N)\to K^*_{e',\gm}(N)$$
    such that
    \begin{itemize}
    \item if $e-e'=\del u$ and $e'-e''=\del u'$, then
        $$\Psi_{e'',u',e'}\smalcirc\Psi_{e',u,e}=\Psi_{e'',u+u',e};$$
    \item for any $v\in \check{C}^1(\gm\lcom;\sS^1)$,
        $$\Psi_{e',u,e}=\Psi_{e',u+\del v,e}.$$
    \end{itemize}
\item
    There is a morphism
    $$H^2(\gm\lcom;\sS^1)\to \mathop{Aut} K^*_{e,\gm}(N).$$
    The ring structure on $K^{*+d}_{e,\gm}(N)$, up to an isomorphism, depends
    only on the cohomology class $[e]\in H^3(\gm\lcom;\sS^1)$. The isomorphism
%!!! I'm not sure I understand what you mean by "induced from ...".
    is unique up to an automorphism of $K^{*+d}_{e,\gm}(N)$ induced from
    $H^2(\gm\lcom;\sS^1)$.
\end{enumerate}
}

As an application, we consider twisted $K$-theory groups of an inertia groupoid.
Let $\gm: \gm_1\toto \gm_0$ be a Lie groupoid and consider the crossed module
$S\gm \to \gm$. As before, $\Lambda \gm: S\gm_1\rtimes \gm_1\toto S\gm_1$
denotes the inertia groupoid of $\gm$. Any element in the image of the
transgression map $T_1: H^3(\gm\lcom; \sS^1)\to H^2(\Lambda \gm\lcom; \sS^1)$ is
$2$-multiplicative. Thus one obtains a ring structure on the corresponding
twisted $K$-theory groups. Since $H^3(\gm\lcom; \sS^1)$ classifies $2$-gerbes,
we conclude that the twisted $K$-theory groups on the inertia
stack twisted by a $2$-gerbe over the stack admits a ring structure.

\vskip.2in

\noindent{\bf Theorem F.} {\it
Let $\gm_1\toto \gm_0$ be a proper Lie groupoid such that $S\gm_1$ is a manifold
and $S\gm_1\to \gm_0$ is $\gm$-equivariantly $K$-oriented (these assumptions
hold, for instance, when $\gm$ is proper and \'etale, or when $\gm$ is a
compact, connected and simply connected Lie group). Let
$d=\dim S\gm_1-\dim\gm_0$.
\begin{enumerate}
\item For any $e\in \check{Z}^3( \gm\lcom;\sS^1)$, the twisted $K$-theory groups
    $K^{*+d}_{e,\gm}(S\gm)$ are endowed with a ring structure
    $$K^{i+d}_{e,\gm}(S\gm)\otimes K^{j+d}_{e,\gm}(S\gm)
        \to K^{i+j+d}_{e,\gm}(S\gm).$$
\item Assume that $e$ and $e'\in \check{Z}^3(\gm\lcom;\sS^1)$ satisfy
    $e-e'=\del u$ for some $u\in \check{C}^2(\gm\lcom;\sS^1)$. Then there is a
    ring isomorphism
    $$\Psi_{e',u,e}:K^*_{e,\gm}(S\gm)\to K^*_{e',\gm}(S\gm)$$
    such that
    \begin{itemize}
    \item if $e-e'=\del u$ and $e'-e''=\del u'$, then
        $$\Psi_{e'',u',e'}\smalcirc\Psi_{e',u,e}=\Psi_{e'',u+u',e};$$
    \item for any $v\in \check{C}^1(\gm\lcom;\sS^1)$,
        $$\Psi_{e',u,e}=\Psi_{e',u+\del v,e}.$$
    \end{itemize}
%!!! What I said before -- this one could use clarification too.
\item There is a morphism
    $$H^2(\gm\lcom;\sS^1)\to \mathop{Aut} K^*_{e,\gm}(S\gm).$$
    The ring structure on $K^{*+d}_{e,\gm}(S\gm)$, up to an isomorphism,
    depends only on the cohomology class $[e]\in H^3(\gm\lcom;\sS^1)$. The
    isomorphism is unique up to an automorphism of $K^{*+d}_{e,\gm}(S\gm )$
    induced from $H^2(\gm\lcom;\sS^1)$.
\end{enumerate}
}

As a special case, when $\gm$ is a compact, connected and simply connected,
simple Lie group $G$, $SG \cong G$, and the $G$-action on $G$ is by conjugation,
then $T_1: H^3 (G\lcom , \sS^1)\to H^2((G\rtimes G)\lcom, \sS^1)$ is an
isomorphism and $H^2 (G\lcom , \sS^1)=0$. Thus, as a consequence, we prove
 the following

\vskip.2in

\noindent{\bf Theorem G.} {\it
Let $G$ be a compact, connected and simply connected, simple Lie group, and
$[c]\in H^2((G\rtimes G)\lcom;\sS^1) \cong \zz$. Then the equivariant twisted
$K$-theory group $K_{[c], G}^* (G)$ is endowed with a canonical ring structure
$$K^{i+d}_{[c],G}(G)\otimes K^{j+d}_{[c],G}(G)\to K^{i+j+d}_{[c], G}(G),$$
where $d=\dim G$, in the sense that there is a canonical isomorphism of the
rings when using any two $2$-cocycles in $\check{Z}^2((G\rtimes G)\lcom;\sS^1)$
which are in the  images of  the transgression $T_1$.}

\vskip.2in

This paper is organized as follows. Section 2 is devoted to preliminaries. In
particular, we introduce \DG-categories and cohomology of \DG-spaces. In Section
3, we give the construction of the transgression maps and discuss their
properties. Section 4 is devoted to the discussion of the ring structures of
twisted equivariant $K$-theory groups.

Note that Freed-Hopkins-Teleman has proved   a remarkable theorem 
that the equivariant twisted $K$-theory group of a compact
connected   Lie  group admits
a ring structure which is isomorphic to
 Verlinde algebra \cite{F, FHT0}.
Indeed the idea of using  multiplicative cocycles
(called equivariantly primitive in \cite{FHT02}) in constructing
the product on twisted equivariant $K$-theory has been known
 in the community.  See \cite{FHT02, FHT1, CJMSW, CW1} for instance.
In \cite{CJMSW, CW1}, there is a notion of
``multiplicative bundle gerbes", which seems  
to be   equivalent  to our ``2-multiplicativity"  except for
the fact that
 ``multiplicative bundle gerbes" in the
sense of  Carey et. al. \cite{CJMSW, CW1} are gerbes over $G$, 
while our gerbes are over $[G/G]$ (i.e. $G$-equivariant gerbes).
In \cite {CW1},   it has been addressed the issue  that
 the multiplicative property is related to the ring structure
on the twisted K-theory groups.
The  idea of considering the ring
structure of $K$-theory groups twisted by classes arising
from the transgression $H^4(BG, \zz)\to H_G^3(G, \zz)$ is
also standard \cite{FHT02, FHT1, CW1}.
 The ring structures on twisted $K$-theory of orbifolds have 
been studied independently by Adem--Ruan--Zhang using a different method
\cite{ARZ}.

However, we would like to stress that the main purpose of 
  our paper is to  emphasize the importance of 
the use of techniques of groupoid $C^*$-algebras and
$KK$-theory in the study of theory of twisted
$K$-theory. This is an advantage of  working
 with  twisted K-theory defined as the K-theory
    group of a certain $C^*$-algebra.
 We also aim to
 provide a different cohomological interpretation of 
    the transgression maps and of 2-multiplicativity in a general framework.

{\bf Acknowledgments.}
We would like to thank several institutions for their hospitality while work 
on this project was done: Penn State University (Tu), and
 Universit\'e Pierre et Marie Curie, Universit\'e de Metz (Xu). 
We also wish to thank Eckhard Meinrenken
and Yong-Bin Ruan for useful discussions  and comments.

\section{Preliminaries}

\subsection{General notations and definitions} 
\label{sec:pre}
Given any category $\cC$ (in particular any groupoid), the collection of objects is denoted by $\cC_{0}$ and the collection
of morphisms is denoted by $\cC_{1}$. We use $\gm$ or $\gm\toto \gm_0$
to denote a groupoid. As usual, $\gm$ is identified with its set of
arrows $\gm_1$.

If $f:x\to y$ is a morphism, then $x$ is called the source of $f$ and
is denoted by $s(f)$, and $y=t(f)$ is  called the target of $f$.
Hence  the composition $fg$ is defined if and only if $s(f)=t(g)$.

Given any $A\subset \cC_{1}$, by  $A^y$, $A_x$ and $A_x^y$
we  denote $A\cap t^{-1}(y)$, $A\cap s^{-1}(x)$ and $A_x\cap A^y$,
respectively.

For all $n\ge 1$, we denote by $\cC_n$ the set of composable $n$-tuples,
i.e. 
$$\cC_n=\{(f_1, \ldots,f_n)|\; s(f_1)=t(f_2),\ldots, s(f_{n-1})=t(f_n)\}.$$

Let $\gm$ be a groupoid and $f:M\to \gm_0$ be a map. We will denote by
$f^*\gm$, or by $\gm[M]$ if there is no ambiguity, {\em the pull-back groupoid}
defined by
$$\gm[M]_0=M, \ \  \gm[M]_1=\{(x,y,g)\in M\times M\times\gm \vert\;
f(x)=t(g),\; f(y)=s(g)\}$$ 
with source and target maps
$t(x,y,g)=x$, $s(x,y,g)=y$, product $(x,y,g)(y,z,h)=(x,z,gh)$ and inverse
$(x,y,g)^{-1}=(y,x,g^{-1})$. In other words, $\gm[M]$ is the fibered product
of the pair groupoid $M\times M$ and $\gm$ over $\gm_0\times \gm_0$.

Let us recall the definition of an action of a groupoid. By definition,
{\em a right action} of  a groupoid $\gm$ on a space $Z$  is given by
\begin{itemize}
\item[(i)] a map $J: Z\to \gm_0$, called the momentum map;
\item[(ii)] a map $Z\times_{\gm_0}\gm:=\{(z,g)\in Z\times \gm\vert\;
J(z)=t(g)\}\to Z$, denoted by $(z,g)\mapsto zg$, satisfying
 $J(zg)=s(g)$,
$z(gh)=(zg)h$ and $z\cdot J(z)=z$
 whenever $J(z)=t(g)$ and $s(g)=t(h)$.
\end{itemize}
Then the {\em transformation groupoid}
(also called {\em the crossed product groupoid})
 $Z\rtimes \gm$ is defined
by $(Z\rtimes \gm)_0=Z$, and
 $(Z\rtimes \gm)_1=Z\times_{\gm_0}\gm$, while
 the   source map,  target map and the product are
$s(z,g)=zg$, $t(z,g)=z$ , $(z,g)(zg,h)=(z,gh)$.

A groupoid $\gm$ is said to be {\em proper} if $(t,s):\gm\to \gm_0\times\gm_0$
is  a proper map. An action of $\gm$ on $Z$ is proper if $Z\rtimes\gm$ is
a proper groupoid.

%\subsection{Semi-direct products}
%\label{sec:semiproduct}
\begin{defn}
Let $N\toto N_0$ and $\gm\toto \gm_0$ be groupoids. We say
 that $\gm$ acts on $N$ by automorphisms if
both $N$ and $N_0$ are right   $\gm$-spaces and
the  actions are compatible  in the following sense
\begin{itemize}
\item the source and target maps $s,t:N \to N_0$ are $\gm$-equivariant,
\item $x^gy^g=(xy)^g$ for all $(x,y,g)\in N\times N\times \gm$ 
whenever either side makes sense. Here $x^g$ denotes
 the action of $g\in \gm$ on $x\in N$.
\end{itemize}
\end{defn}

Given such a pair of groupoids $(N, \gm)$, 
one can form the {\em semi-direct product groupoid}
$N\rtimes \gm$, where the unit space is $N_0$,
the space of morphisms is
$$(N\rtimes \gm)_1=\{(x,g)\in N\times \gm\vert\;
x^g\mbox{ makes sense}\},$$
the target, the source, the multiplication and the
inverse are defined by
$$t(x,g)=t(x),\  s(x,g)=s(x^g),\  (x,g)(y,h)=(xy^{(g^{-1})},gh), 
\ \mbox{and } (x, g)^{-1}=((x^g)^{-1}, g^{-1}).$$

%\subsection{crossed modules}

\begin{defn}
\label{def:cross}
A crossed module is a groupoid morphism
$$\xymatrix{
N\ar[r]^\varphi\ar[d]\ar@<-.6ex>[d] & \gm\ar@<-.6ex>[d]\ar[d] \\
N_0\ar[r]^{=} & \gm_0
}$$
where $N\toto N_0$ is a bundle of groups,
together with an action of $\gm$  on $N$ by automorphisms
 such that
\begin{itemize}
\item[(i)] $\varphi(x^g)=\varphi(x)^g$ for all $x\in N$ and $g\in \gm$
such that $x^g$ makes sense;
\item[(ii)] $x^{\varphi(y)}=x^y$ for all composable pairs $(x,y)\in N_{2}$.
\end{itemize}
Here $\varphi(x)^g:=g^{-1}\varphi(x)g$ and $x^y:=y^{-1}xy$. For short,
a crossed module is denoted by  $N\stackrel{\varphi}{\to}\gm$.
\end{defn}
%In (i) and (ii), note that since $N_0=\gm_0$, the source and target maps
%$s$ and $t:N\to \gm_0$ necessarily coincide, thus
%$\varphi(x)^g:=g^{-1}\varphi(x)g$ and $x^y:=y^{-1}xy$ make sense.

%\begin{tiny}
%A crossed module naturally  gives rise to a double groupoid
%$$\xymatrix{
%N_1 \rtimes \gm_1\ar[r]\ar@<-.6ex>[r]\ar[d]\ar@<-.6ex>[d]
%  & \gm_1\ar@<-.6ex>[d]\ar[d] \\
%N_1\ar[r]\ar@<-.6ex>[r] & N_0=\gm_0
%},$$
%
%Here the first groupoid structure $N_1\rtimes \gm_1\toto N_1$ is given by
%$$s(x,g)=x^g, \  t(x,g)=x \mbox{ and } (x,g)(x^g,h)=(x,gh);$$
%and the  second groupoid structure $N_1\rtimes \gm_1\toto N_1$
% is given by
%$$s'(x,g)=t'(x,g)=g, \ \   (x,g)(y,g)=(xy,g). $$
%
%By abuse of notation, we will also use $N\rtimes \gm$ to denote
%this double groupoid.
%\end{tiny}

A standard example of crossed modules is the  inertia groupoid.
Let $\gm\toto \gm_0$ be a groupoid and
 $S\gm=\{g\in \gm\vert\; s(g)=t(g)\}$ be the space of closed loops  in 
 $\gm$.
Then the canonical inclusion
$S\gm\to \gm$, together with the conjugation
 action of $\gm$  on $S\gm$, 
forms a crossed module, where the crossed-product groupoid $S\gm\rtimes\gm$
is called the inertia groupoid and is denoted by $\Lambda\gm$.

%\begin{tiny}
%Thus
% the inertia groupoid $\Lambda \gm:=S\gm\rtimes \gm$
%is indeed  a double groupoid:
%$$\xymatrix{
%S\gm_1 \rtimes \gm_1\ar[r]\ar@<-.6ex>[r]\ar[d]\ar@<-.6ex>[d] & \gm_1\ar@<-.6ex>[d]\ar[d] \\
%S\gm_1\ar[r]\ar@<-.6ex>[r] & \gm_0
%}.$$
%\end{tiny}

\begin{defn}
Let $N\stackrel{\varphi}{\to}\gm$ and $N'\stackrel{\varphi'}{\to}\gm'$
be crossed modules.
A crossed module morphism
$\tau: (N\stackrel{\varphi}{\to}\gm) \to (N'\stackrel{\varphi'}{\to}\gm')$
is a commutative diagram of groupoid morphisms
$$\xymatrix{
N\ar[r]^\varphi\ar[d]_\tau & \gm\ar[d]^\tau \\
N'\ar[r]^{\varphi'} & \gm'
}$$
satisfying the condition
\begin{equation}
\label{eq:tau}
\tau (x^g)=\tau(x)^{\tau(g)},  \ \ \mbox{ for all compatible } \  x \in N, g\in \gm. 
\end{equation}
\end{defn}

%The following proposition can be easily verified directly.
%
%\begin{prop}
%Let $N\stackrel{\varphi}{\to}\gm$ and $N'\stackrel{\varphi}{\to}\gm'$
%be crossed modules. There is a bijection between
%crossed module morphisms from $(N\stackrel{\varphi}{\to}\gm)$
%to $N'\stackrel{\varphi}{\to}\gm'$ and double
%groupoid morphisms from $N\rtimes \gm$ to $N'\rtimes \gm'$.
%\end{prop}

Given a crossed module $ N\stackrel{\varphi}{\to}\gm$, 
since $\varphi$ maps $N$
to $S\gm$, we have a natural crossed module
morphism from $ N\stackrel{\varphi}{\to}\gm$ to
$S\gm \to \gm$.
%, thus a double groupoid morphism
%from  $N\rtimes \gm$ to $\Lambda \gm$.

\subsection{$\cC$-spaces and sheaf cohomology}
\label{sec:sheaf}
Let $\cC$ be a category.
By a $\cC$-space, we mean  a contravariant functor
from the category $\cC$ to the category
of topological spaces. Similarly,
one defines a $\cC$-manifold.
% (in this sub-section we are mainly interesting
%in the case when $\cC=\Delta$, $\Delta\times\Delta$
%or $\Delta_2$).
 Consider a $\cC$-space $M\lcom$.
% (i.e. a simplicial, or bi-simplicial, or $\Delta_2$-space).
 Let $\cC_{M\lcom}$ 
be the category whose objects are pairs $(i,U)$, with $i\in \cC_{0}$
and $U$ an open subset of $M_i$, such that morphisms  from
$(i,U)$ to  $(j,V)$ consist of those $f\in\Hom_{\cC}(j,i)$ for which
$\tilde{f}(U)\subset V$. By definition, an abelian presheaf on the
$\cC$-space $M\lcom$ is an abelian presheaf on the category
$\cC_{M\lcom}$, i.e.  a contravariant functor
from the category $\cC_{M\lcom}$ to the category of abelian groups.
A presheaf $\aA$ on $M\lcom$ restricts to a presheaf $\aA_i$
on each space $M_i$. We say that $\aA$ is a sheaf if each $\aA_i$ is a
sheaf.

More concretely, a sheaf $\aA$ on $M\lcom$ is given by
a family $(\aA^i)_{i\in \cC_{0}}$
such that $\aA^i$ is a sheaf on $M_i$, together with restriction maps
$\tilde{f}^*:\aA^j(V)\to \aA^i(U)$, for each
$\tilde{f}\in\Hom_{\cC_{M\lcom}}((i,U),(j,V))$, 
satisfying the relation
$(\tilde{f}\circ \tilde{g})^*=\tilde{g}^*\circ \tilde{f}^*$ \cite{Deligne}.
In a similar fashion, one defines the notion of  a sheaf over
 a $\cC$-manifold.
Note that a big sheaf over the site of all smooth manifolds
naturally induces a sheaf on a $\cC$-manifold.
For instance, the sheaf of real-valued smooth functions $\rR$,
the sheaf of $S^1$-valued smooth functions $\sS^1$, and the sheaf of
$q$-forms $\Omega^q$ (for fixed $q$) are examples of such sheaves.

Assume  that $\aA$ is a sheaf on a $\cC$-space $M\lcom$. In order
to define cohomology groups $H^*(M\lcom, \aA)$,
 one needs  an  extra structure on $\cC$.

We  say that a  category $\cC$ is  a \DG\ category  if every object
$k\in \cC_{0}$ is labeled by an integer $\deg(k)\in\nn$
(in other words, there is a functor from the category
$\cC$ to the groupoid $\nn\times \nn\toto \nn$), and  moreover
 it is endowed with a set $A\subset \cC_{1}$ and $\varepsilon:A\to\zz$
satisfying
\begin{itemize}
\item[(i)] $A_k$ is finite for all $k\in \cC_{0}$;
\item[(ii)] for all $f\in A$, $\deg(f)=1$, where
$\deg(f)=\deg(t(f))-\deg(s(f))$;
\item[(iii)] for all $f\in \cC_{1}$,
$$\sum_{{f'\smalcirc f'' = f}\atop{f',f''\in A}}\varepsilon(f')\varepsilon(f'')=0.$$
\end{itemize}
Note that the sum in (iii) is finite due to (i).

Given a \DG\ category $\cC$, a $\cC$-space $M\lcom$
and a sheaf $\aA$ over $M\lcom$, let
$$C^n(M\lcom;\aA)=\oplus_{\deg k =n}\aA(M_k).$$
Then $C^*(M\lcom;\aA)$ is endowed with a degree $1$ differential
$${\del}\omega = \sum_{f\in A_k}\varepsilon(f) \tilde{f}^*\omega, $$
 $\forall \omega\in \aA(M_k)$.  It is simple to check that
 ${\del}^2=0$. Then we obtain a cochain complex 
$(C^*(M\lcom;\aA), \ \del)$. By $Z^*(M\lcom;\aA)$ we denote
its space of cocycles and $\hh^*(M\lcom;\aA)$ the cohomology group.

Now  given a complex of sheaves over $M\lcom$,  bounded below:
$\aA_0\stackrel{d}{\to} \aA_1\stackrel{d}{\to}
 \aA_2\stackrel{d}{\to}\cdots$, 
then
$$C^{p,q}(M\lcom;\aA\lcom):=C^p(M\lcom;\aA_q) =\oplus_{\deg k =p}\aA_q (M_k)$$
is endowed with a  double complex  structure with differentials
$d$ and ${\del}$. We  denote by $\delta$ the
total differential $(-1)^pd+{\del}$,  and by $\hh^*(M\lcom;\aA\lcom)$
its hypercohomology groups.

For a sheaf  $\aA$  on  $M\lcom$, let
 $\aA\lcom$ be a   complex of sheaves over $M\lcom$ which is
 an injective  resolution of $\aA$. Then $H^*(M\lcom;\aA)$
is defined to be $\hh^*(M\lcom;\aA\lcom)$,
 called the sheaf
cohomology group of $M\lcom$ with coefficients in $\aA$.

A particular case is the following:
if $M\lcom$ is a $\cC$-manifold and
 $\aA_q=\Omega^q$:
$\Omega^0 \stackrel{d}{\to} \Omega^1\stackrel{d}{\to} 
 \Omega^2\stackrel{d}{\to}\cdots$ is the de Rham
complex of sheaves, the group $\hh^*(M\lcom;\Omega\upcom)$ 
is called the de Rham cohomology of $M\lcom$ 
and is denoted by $H^*_{dR}(M\lcom)$.
%%\begin{tiny}
%In our cases of interest (i.e. $\cC=\Delta$ or $\Delta^2$ or $\Delta_2$),
%\end{tiny}
It is isomorphic to  $H^*(M\lcom;\rr)$.

\subsection{Simplicial spaces}
\label{sec:simplicial}
Recall that the simplicial category, denoted by $\Delta$,
has as  objects the set of non-negative integers, and $\Hom_\Delta(k,k')$
is the set of non-decreasing maps from  $[k]$ to $[k']$, where
$[k] =\{0,\ldots,k\}$. A $\Delta$-space is thus called a
 simplicial (topological) space, and a $\Delta$-manifold
is a simplicial manifold.
% is a contravariant functor from the

In a down-to-earth term,
a simplicial space is given by a sequence $M\lcom=(M_n)_{n\in\nn}$ of
spaces,  and for each $f\in\Hom_\Delta (k,n)$, we are given a  map
(called face  or degeneracy map depending
which of $k$ and $n$ is larger)  $\tilde{f}: M_n\to M_k$ such that
 $\tilde{f}\smalcirc\tilde{g}=\widetilde{g\smalcirc f}$.

Similarly,  denote by $\bar{\Delta}$ the category obtained from
$\Delta$ by identifying $f:[k]\to [n]$ with $f':[k]\to [n]$ whenever
both $f$ and $f'$ are constant. We will call $\bar{\Delta}$
the reduced simplicial category.

A groupoid naturally gives rise to a  simplicial space.
To see this, consider the pair groupoid
$[n]\times [n]\toto [n]$. For a groupoid $\gm\toto \gm_0$,
let $\gm_n=\Hom ([n]\times [n], \gm )$ be the space of homomorphisms
from the pair groupoid $[n]\times [n]\toto [n]$ to $\gm$.
 Any $f\in\Hom_\Delta(k,n)$ gives rise to a
groupoid homomorphism from $[k]\times [k]\toto [k]$ to $[n]\times [n]\toto [n]$, again denoted by $f$. It thus, in turn, induces a 
map $\tilde{f}: \gm_n(= \Hom ([n]\times [n], \gm ))\to
 \gm_k (=\Hom ([k]\times [k], \gm ))$, which is the ``face/degeneracy" map.
Note that  $\gm_n$ can be identified with the space of
composable $n$-tuples: 
$\gm_n=\{(g_1,\dots,g_n)\vert\; g_1\cdots g_n\mbox{ makes sense}\}$
since the groupoid $[n]\times [n]\toto [n]$ is
generated by elements $(i-1,i)$ ($1\le i\le n$).
Hence  any groupoid morphism from $[n]\times [n]\toto [n]$ to $\gm\toto \gm_0$ is uniquely determined
by the image of each element $(i-1,i)$,
 which is  denoted by $g_i$ ($1\le i\le n$).

Moreover, the simplicial space structure descends to a reduced
simplicial structure when the source and target maps coincide,
i.e. when $\gm \toto  \gm_0$ is a bundle of groups.

Recall that the simplicial category
$\Delta$ is equipped with  a natural \DG\ 
category structure. The degree
map is obviously the identity map $\Delta_{0}\to \nn$.
For all $k\in\nn$, let $\varepsilon^k_i:[k]\to [k+1]$ 
be the unique increasing
map which omits $i$ ($i=0,\ldots,k+1$):
$$\varepsilon^k_i (0)=0, \ldots, \varepsilon^k_i (i-1)=i-1, \varepsilon^k_i
(i)=i+1, \ldots, \varepsilon^k_i (k)=k+1.$$
 We will omit the superscript $k$
if there is no ambiguity. Let $\varepsilon(\varepsilon_i)=(-1)^i$.
Then the pair $(A,\varepsilon)$, where $A_k=\{\varepsilon^k_i\vert\;
i\in [k+1]\}$ is a \DG\ structure on $\Delta$.
For sheaf  cohomology  of simplicial manifolds, we refer the
reader to \cite{tu05, Dupont} for details. 

%Similarly, a double groupoid gives rise to a bi-simplicial space, i.e.,
%a contravariant functor from $\Delta\times\Delta$ to the category of
%topological spaces. Indeed, for all $(k,l)\in\nn^2$, let $G(k,l)$
%be the double groupoid 
%$$\xymatrix{
%([k]\times [k])\times ([l]\times [l])\ar[r]\ar@<-.6ex>[r]\ar[d]\ar@<-.6ex>[d]
% & [k]\times ([l]\times [l])\ar@<-.6ex>[d]\ar[d] \\
%([k]\times [k])\times [l]\ar[r]\ar@<-.6ex>[r] & [k]\times [l] 
%}.$$
%
%Given any double groupoid $\gm$, let $\gm_{k,l}=\Hom(G(k,l),\gm)$
%be the space of double groupoid morphisms from $G(k,l)$
%to $\gm$. More explicitly, $\gm_{k,l}$ is the space of arrays
%$(g_{i,j})_{0\le i\le k\atop 0\le j\le l}$ such that $g_{i-1,j}$ is composable
%with $g_{i,j}$ for the first groupoid law, $g_{i,j-1}$ is composable
%with $g_{i,j}$ for the second groupoid law, $g_{0,j}$ is a unit for the first
%groupoid law and $g_{i,0}$ is a unit for the second groupoid law.
%
%Then any morphism $f$ in $\Hom_{\Delta\times\Delta}
%((k,l), (k',l'))$ induces a double groupoid morphism 
%$f:G(k,l)\to G(k',l')$, which in turn  induces
%a pull-back map $\tilde{f}:\gm_{k',l'}\to \gm_{k,l}$.
%In particular,
%any crossed module $N\to \gm $ gives rise to a bi-simplicial space,
%which is in
%fact the one induced by the double groupoid $N\rtimes \gm$.

Suppose now that $\cC$ and $\cC'$ are two \DG\ categories.
 Then the product  $\cC''=\cC'\times \cC$
is naturally a  \DG\  category,
where   $\deg(k,l)=
\deg(k)+\deg(l)$, $A''=A'\times \{1\} \cup \{1\}\times A$
and $\varepsilon(f,1)=\varepsilon(f)$ for all $(f,1)\in {A_k}'\times {\cC}_l$,
$\varepsilon(1,g)=(-1)^{\deg k}\varepsilon(g)$ for all
$(1,g)\in {\cC_k}'\times {A}_l$.

In particular,  $\Delta\times \Delta$ is a \DG\ category. More precisely, if
$\del'=\sum_{i=0}^{k+1} (-1)^i {\tilde{\varepsilon'}}_i^*$ is the
differential with respect to the first simplicial structure (as above
$\varepsilon'_i:[k]\to [k+1]$ is the increasing map that omits $i$)
and $\del=\sum_{i=0}^{l+1} (-1)^i{\tilde{\varepsilon}_i}^*$ is the
differential with respect to the second simplicial structure,
then $\underline{\del}=\del'+(-1)^k\del$ is the differential for
the bi-simplicial structure.

\subsection{$\ff$-spaces}
Let us now consider an extension of $\bar{\Delta}$. Denote
by $\ff$ the category with the same objects as $\Delta$, 
but $\Hom_\ff(k,n)$
consists of $\Hom(\ff_k,\ff_n)$, where $\ff_n$ is the free group on
$n$ generators. More concretely, any element of $\Hom(\ff_k,\ff_n)$
is given by a $k$-tuple $w=(w_1,\ldots,w_k)$ of words in $x_1^{\pm 1},
\ldots, x_n^{\pm 1}$. Denoting by $\hat{w}$ the corresponding
element of $\Hom_\ff(k,n)$, we have the relation
$\widehat{w'\smalcirc w}
=\hat{w'}\smalcirc\hat{w}$, where $(w'\smalcirc w)_i=w'\smalcirc w_i$ 
consists of the word obtained from $w_i$ by substituting each occurrence
 of $x_j$ by $w'_j$.

To see that $\ff$ extends the category $\bar{\Delta}$,
for any  $f\in\Hom_{\bar{\Delta}}(k,n)$,
let $\bar{f}\in\Hom_\ff(k,n)$ be the morphism 
$$\bar{f}(x_i) =x_{f(i-1)+1}\cdots x_{f(i)}, \ \ \ i=1, \ldots, k$$
 (with the convention $\bar{f}(x_i)=1$ if $f(i-1)=f(i)$).
%$$\bar{f}(x_1, \ldots,  x_i, \ldots, x_k)=
 One immediately checks that $\overline{f\smalcirc g}
=\bar{f}\smalcirc\bar{g}$.
 Another way to explain the inclusion
$\bar{\Delta}\subset \ff$ is as follows.
 Any $f\in\Hom_{\bar{\Delta}} (k,n)$ gives rise to a
groupoid homomorphism from $[k]\times [k]\toto [k]$ to $[n]\times [n]\toto [n]$,
 again denoted by $f$. Let $\iota: [n]\times [n]\to \ff_n$ be the
unique groupoid morphism such that $(i-1,i)$ maps to $x_i$. Then
$\bar{f}$ is the unique group homomorphism such that the diagram
$$\xymatrix{
[k]\times [k] \ar[r]^\iota\ar[d]_f & \ff_k\ar[d]^{\bar{f}} \\
[n]\times [n]\ar[r]^\iota & \ff_n
}$$
commutes.

As above,  a $\ff$-(topological) space is 
 a contravariant functor from
$\ff$ to the category of topological spaces. If $G$ is a topological group,
then we obtain
 an associated $\ff$-space by setting $G_n=\Hom(\ff_n,G)\;
(\cong G^n)$.
In particular, since $\ff$ extends the category $\bar{\Delta}$,
$G\lcom=(G_n)_{n\in\nn}$ is a reduced simplicial space and therefore
 a simplicial space. The simplicial structure can be 
seen,  as in Section \ref{sec:simplicial},  by considering
$G$ as a groupoid.

%A simpler way to view that simplicial space is to observe that
%$G_n\cong \Hom([n]\times [n],G)$, and, as above, that
%any nondecreasing map $f:[k]\to [n]$
%gives rise to a groupoid morphism $[k]\times [k]\to [n]\times [n]$,
%hence to a pull-back map $\tilde{f}:G_{n}\to G_k$.

%More generally, any groupoid $\gm\toto \gm_0$
%gives rise to a simplicial space,
%and one checks that
%$\gm_n=\{(g_1,\dots,g_n)\vert\; g_1\cdots g_n\mbox{ makes sense}\}$,
%since any groupoid morphism from $[n]\times [n]$ to $\gm$
%is uniquely determined
%by the image of each element $(i-1,i)$ ($1\le i\le n$).
%Moreover, the simplicial space structure descends to a reduced
%simplicial structure when the source and target maps coincide,
%i.e. when $\gm\to \gm_0$ is a bundle of groups.

%Similarly, a double groupoid gives rise to a bi-simplicial space, i.e.
%a contravariant functor from $\Delta\times\Delta$ to the category of
%topological spaces. Indeed, for all $(k,l)\in\nn^2$, let $G(k,l)$
%be the double groupoid $([k]\times [k])\times ([l]\times [l])$.
%Given any double groupoid $\gm$, let $\gm_{k,l}=\Hom(G(k,l),\gm)$,
%then any double groupoid morphism $f:G(k,l)\to G(k',l')$ induces
%a pull-back map $\tilde{f}:\gm_{k',l'}\to \gm_{k,l}$.

\subsection{$\ff\Delta$-spaces}
We  now introduce  a category $\ff\Delta$.
Objects are pairs $(k,l)\in\nn^2$. To describe morphisms, let us 
introduce some notations: let
 $X_{k,l}$ be the groupoid $\ff_k\times ([l]\times [l])
\toto [l]$, the product of the free group $\ff_k$ with  the pair groupoid
$[l]\times [l]\toto [l]$. Then we define $\Hom((k,l),(k',l'))$
as the set of groupoid morphisms $f:X_{k,l}\to X_{k',l'}$ such that
the restriction of $f$ to the unit space, again denoted by $f:[l]\to [l']$,
is a nondecreasing function. In particular, for $k=0$ we recover the
simplicial category $\Delta$ and for $l=0$ we recover the category $\ff$.
We also note that  the sub-category of $\ff\Delta$ consisting  of morphisms
$f:X_{k,l}\to X_{k',l'}$ of the form $f=(f_1,f_2)$, where
$f_1:\ff_k\to\ff_{k'}$ is a group morphism and $f_2:[l]\times [l]\to
[l']\times [l']$ is a groupoid morphism whose restriction to the unit spaces
$[l]\to [l']$ is nondecreasing, is exactly isomorphic to the
 product category $\ff\times \Delta$.

To understand the   category $\ff\Delta$ in a more concrete way,
consider the following arrows of the groupoid $X_{k,l}$:
\begin{equation}
\label{eq:F}
\tilde{a}=(a, 0, 0), \ \gamma_i=(1,0,i),
\end{equation}
 where
$a\in \ff_k $ and $i=0, \ldots, l$.
They generate $X_{k,l}$ since
any arrow in $X_{k,l}$ can be written in a unique way as
 
\begin{equation}
\label{eq:aij}
(a, i, j)=\gamma_i^{-1}\tilde{a}\gamma_j,
\end{equation}
 where $a\in\ff_k$.
Consider any morphism in $\Hom_{\ff\Delta}((k,l),(k',l'))$,
 whose restriction to the unit space
is  denoted by $f:[l]\to [l']$. Assume that under this
morphism, we have $\tilde{a}\mapsto (\psi (a), f(0), f(0))\in X_{k',l'}$
and $\gamma_i\mapsto (u_i, f(0), f(i))\in X_{k',l'}$, 
where $\psi\in\Hom(\ff_k,\ff_{k'})$,
$f\in \Hom_{\Delta}(l,l')$, and $u=(u_0,\ldots,u_l)\in(\ff_{k'})^{l+1}$.
Thus 
$$(a, i, j)= \gamma_i^{-1}\tilde{a}\gamma_j \mapsto
(u_i^{-1}\psi (a)u_j, f(i), f(j)). $$
 Note that
%Then, any element of $\Hom_{\ff\Delta}((k,l),(k',l'))$ is of the form
%$\gamma_i^{-1}a\gamma_j\mapsto
%\gamma_{f(i)}^{-1}u_i^{-1}\psi(a)u_j\gamma_{f(j)}$,
the triple $(\psi,u,f)$ is uniquely determined modulo the equivalence
relation: 
$(\psi,u,f)\sim (\psi',u',f)$  if $\psi'(a)=\psi(a)^v$ and
$u'_i=v^{-1}u_i$ for some $v\in \ff_k$.

We summarize the above discussion in the following

\begin{prop}
\label{pro:FD}
$\Hom_{\ff\Delta}((k,l),(k',l'))$ can be identified with
triples $(\psi,u,f)$, where
$\psi\in\Hom(\ff_k,\ff_{k'})$, $f\in \Hom_{\Delta}(l,l')$,
and $u=(u_0,\ldots,u_l)\in(\ff_{k'})^{l+1}$,
 modulo the equivalence relation $(\psi,u,f)\sim (\psi',u',f)$
 if and only if $\psi'(a)=\psi(a)^v$ and
$u'_i=v^{-1}u_i$ for some $v\in \ff_k$.

The composition law of morphisms is then $(\psi',u',f')\smalcirc
(\psi,u,f)=(\psi'',u'',f'')$,
 where $\psi''=\psi'\smalcirc\psi$, $f''=f'\smalcirc f$ and 
$u''_i=\psi'(u_i)u'_{f(i)}$.
\end{prop}

\subsection{$\Delta_2$-spaces}
Next we define a category $\Delta_2$ as follows: objects are pairs
of integers $(k,l)\in\nn^2$. $\Hom_{\Delta_2}((k,l),(k',l'))$
consists of triples $(a,b,c)$ such that $a\in \{\emptyset\}\cup
\Hom_\Delta(k,k')$, $b\in \{\emptyset\}\cup\Hom_\Delta(l,k')$,
$c\in\Hom_\Delta(l,l')$, and either $a=\emptyset$ or $b=\emptyset$.

We define the composition as follows.
$$(a',\emptyset,c')\smalcirc (a,\emptyset,c)=(a'\circ a,\emptyset,c'\circ c),
(a',\emptyset,c')\smalcirc (\emptyset,b,c)=(\emptyset,a'\circ b, c'\circ c),
(\emptyset,b',c')\smalcirc (a,b,c)=(\emptyset,b'\circ c,c'\circ c).$$
The associativity can be checked easily and is left to the reader.

It is clear that the bi-simplicial category $\Delta\times\Delta$
embeds into $\Delta_2$ by $(a,c)\in\Hom_\Delta(k,k')\times\Hom_\Delta
(l,l')\cong \Hom_{\Delta\times\Delta}((k,k'),(l,l'))\mapsto
(a,\emptyset,c)\in \Hom_{\Delta_2}((k,k'),(l,l'))$.
\par\medskip

Let us now define a category $\bar{\Delta}_2$, which has the same 
objects as $\Delta_2$, and whose  morphisms are
 obtained from  morphisms of $\Delta_2$ by identifying $(\emptyset,b,c)$,
$(\emptyset,b',c)$, $(a,\emptyset,c)$, $(a',\emptyset,c)$
 whenever $a$, $a'$, $b$ and
$b'$ are constant functions. The resulting element is
denoted by $0_{(k,l),(k', l'),c}$, or  simply by $0_c$ if there is no ambiguity.
 One checks directly that this definition
makes sense and that $0_{c'}\smalcirc (a,b,c)=0_{c'\smalcirc  c}$,
$(a',b',c')\smalcirc 0_c = 0_{c'\smalcirc  c}$.

The category $\bar{\Delta}_2$ embeds into $\ff\Delta$ as a subcategory.
 Indeed,
to a triple $(a,b,c)\in \Hom_{\Delta_2}((k,l),(k',l'))$
one associates a triple $F(a,b,c)=
(\psi_a,u_b,c)\in\Hom_{\ff\Delta}((k,l),(k',l'))$
as follows. Denote by $x_i$ the generators of $\ff_k$, and let
$y_i=x_1\cdots x_i$, with the convention $y_0=1$.
Let
\begin{equation}\label{eqn:psiu}
\psi_a(y_i)=y_{a(0)}^{-1}y_{a(i)}\mbox{ and }
(u_b)_i=y_{b(i)},
\end{equation}
where,  by convention,  $a(i)=0$ if $a=\emptyset$.
Then $(a,b,c)\mapsto (\psi_a,u_b,c)$
is injective, and a simple calculation shows that
$F((a',b',c')\smalcirc (a,b,c))=F(a',b',c')\smalcirc F(a,b,c)$.

Note also that $\bar{\Delta}\times \Delta\subset \bar{\Delta}_2$
by $(a,c)\mapsto (a,\emptyset,c)$.

The above discussion can be  summarized by
the following  diagram, where all maps are embeddings except for
the two horizontal arrows on the left:

$$\xymatrix{
\Delta\times \Delta\ar[r]\ar[d] & \bar{\Delta}\times \Delta \ar[r]\ar[d] &
\ff\times\Delta\ar[d]\\
\Delta_2\ar[r]&\bar{\Delta}_2 \ar[r] & \ff\Delta
}$$

\subsection{The $\ff\Delta$- and bi-simplicial spaces 
associated to a crossed module}

We show that a crossed module $N\stackrel{\varphi}{\to}\gm$ naturally
 gives rise to a $\ff\Delta$-space. Recall that every groupoid $\gm$
gives rise to a crossed module $S\gm\to \gm$, with $\gm$
acting on $S\gm$ by conjugations (see Section \ref{sec:pre}).
Let $(N\rtimes \gm)_{k,l}$ be the
space  of morphisms of crossed modules from $(SX_{k,l}\to X_{k,l})$
to $(N\to \gm)$. Since a groupoid morphism  $f:X_{k,l}\to X_{k',l'}$
induces a crossed module morphism from
$(SX_{k,l}\to X_{k,l})$ to $(SX_{k',l'}\to X_{k',l'})$,
we obtain a map $\tilde{f}:(N\rtimes \gm)_{k',l'}\to (N\rtimes \gm)_{k,l}$.
Hence, $(N\rtimes \gm)\lcomm = ((N\rtimes \gm)_{k,l})_{(k,l)\in\nn^2}$
is endowed with a structure of $\ff\Delta$-space.

To see this $\ff\Delta$-space $(N\rtimes \gm)\lcomm$ in
down-to-earth terms,
%Recall that a crossed module morphism 
%$\tau: (N\stackrel{\varphi}{\to}\gm) \to (N'\stackrel{\varphi}{\to}\gm'$
%is a commutative diagram of groupoid morphisms
%$$\xymatrix{
%N\ar[r]^\varphi\ar[d]_\tau & \gm\ar[d]^\tau \\
%N'\ar[r]^{\phi'} & \gm'
%}$$
%satisfying the condition
%\begin{equation}
%\tau (x^g)=\tau(x)^{\tau(g)},  \ \ \mbox{ for all compatible } \  x \in N, g\in \gm 
%\end{equation}
%
%Hence one can check that  an element of $(N\rtimes \gm)_{k,l}$ is determined
%by a pair $(h,g)$, where $h\in\Hom(\ff_k,N)$ and $g=
%(g_0,\ldots,g_l)$, with $g_i\in \gm$, and
%$t(\varphi(h(a)))=t(g_0)=\cdots=t(g_l)$ for all $a\in \ff_k$.
%To see the precise correspondence,
 let $\tau \in (N\rtimes \gm)_{k,l}$
be an arbitrary  element, which corresponds to a 
commutative diagram of groupoid morphisms
$$\xymatrix{
SX_{k, l}\ar[r]\ar[d]_\tau & X_{k, l}\ar[d]^\tau \\
N\ar[r]^{\varphi} & \gm
}$$
satisfying Eq. \eqref{eq:tau}. 
Consider $\ff_k$ as a subgroupoid of $SX_{k,l}\toto [l]$ by
identifying $\ff_k$ with $\ff_k\times \{(0, 0)\}$, and  let
$h: \ff_k\to N$ be the restriction of $\tau: SX_{k,l}\to N$ to $\ff_k$.
And let $g_i =\tau (\gamma_i), \ i=0, \ldots, l$, where
$\gamma_i$ is as in Eq. \eqref{eq:F}.
Using Eq. \eqref{eq:aij}, for elements in $SX_{k,l}$, the map
$\tau$ is  given by
$$ \tau (a, i, i)=\tau ((a, 0, 0)^{\gamma_i})=h(a)^{g_i}, \  \ \forall 
a\in \ff_k, i=0, \ldots, l, $$
while for elements in  $X_{k,l}$, the map $\tau$ is then given by
$$\tau (a, i, j)=\tau ( \gamma_i^{-1}\tilde{a}\gamma_j)=g_i^{-1}\tau(\tilde{a})
g_j=g_i^{-1}\varphi(h(a))g_j, \ \ \forall a\in   \ff_k, i, j=0, \ldots, l.$$

%The morphism from $SX_{k,l}=\ff_k\times [l]$ to $N$ is
%$(a,i)\mapsto h(a)^{g_i}$ and the morphism from
%$X_{k,l}$ to $\gm$ is $(a,i,j)\mapsto g_i^{-1}\varphi(h(a))g_j$.

Note that the pair $(h,g)$ is not unique
and it is  uniquely determined modulo the equivalence
relation $(h,g)\sim (h',g')$, if
 $h'(a)=h(a)^r$ and
$g'_i=r^{-1}g_i$ for some $r\in \gm$.
Also note that $\Hom (\ff_k, N)$ can be naturally
identified with $N_k$ by 
identifying $h\in \Hom (\ff_k, N)$ with $(h(w_1), \ldots, h(w_k ))\in
N_k$, where $w_1, \ldots, w_k$ are generators of $\ff_k$.
Hence it follows that $(N\rtimes \gm)_{k,l}$ is isomorphic, as a space,
 to the quotient of 
\begin{equation}
\label{eq:ngm}
\{(x_1,\ldots,x_k;g_0,\ldots,g_l)\in N^k\times \gm^{l+1}|\;
t(x_i)=t(g_j)\;\forall i,j\}
\end{equation}
modulo the equivalence relation
$$(x_1,\ldots,x_k;g_0,\ldots,g_l)\sim(r^{-1}x_1r,\ldots,r^{-1}x_kr;r^{-1}
g_0,\ldots,r^{-1}g_l).$$

A simple calculation shows that
the $\ff\Delta$-structure on $(N\rtimes \gm)\lcomm$ is given by
$(\psi,u,f)^{\tilde{\;}}(h',g')=(h,g)$, where
\begin{equation}\label{eqn:hpsi}
h(a)=h'(\psi(a))\mbox{ and }g_i=h'(u_i)g'_{f(i)}.
\end{equation}
Here $(\psi,u,f)$ is a triple
defining a morphism in  $\Hom_{\ff\Delta}((k,l),(k',l'))$ 
as in Proposition \ref{pro:FD}.
Since any  $\ff\Delta$-space  is automatically
a bi-simplicial space, $(N\rtimes \gm)_{k,l}$  is
naturally a bi-simplicial space.
On the other hand, for any
 fixed $k$, the groupoid $\gm$ acts on the space
$N_k$. Hence we obtain  a simplicial space 
\begin{equation}\label{sim.ma}
\xymatrix{
\ldots N_k \rtimes \gm_{2} 
\ar[r]\ar@<1ex>[r]\ar@<-1ex>[r] & N_k \rtimes \gm_{1}\ar@<-.5ex>[r]\ar@<.5ex>[r]
&N_k\rtimes \gm_{0}\,,}
\end{equation}
where
$$N_k\rtimes\gm_l=\{(x_1,\ldots,x_k;g_1,\ldots,g_l)\in N^k\times \gm^l\vert\;
t(x_1)=\cdots=t(x_k)=t(g_1), g_1\cdots g_l\mbox{ makes sense}\}.$$
Moreover, for any  fixed $l$, we get a simplicial structure on
$N\lcom\rtimes \gm_l$ since $N$ is a groupoid. In fact,
$N\lcom\rtimes \gm\lcom$ is a bi-simplicial space.
%To describe the bi-simplicial structure on
%$(N\rtimes \gm)_{k,l}$, consider
%$$(N\rtimes \gm)_{k,l}^{\Delta\times \Delta}=
%\Hom (G(k, l), N\rtimes \gm).$$
%Then
%$$(N\rtimes \gm)_{k,l}^{\Delta\times \Delta}=
%\{(x_1,\ldots,x_k;g_1, \ldots ,g_l)|s(x_1)=\ldots=s(x_k)=t(g_1),
%g_1\cdots g_l \mbox{ makes sense }\}$$
%
 Introduce a map
\begin{eqnarray}
\Phi: N_k\rtimes \gm_l&\lon & (N\rtimes \gm)_{k,l} \nonumber\\
(x_1,\ldots,x_k;g_1,\ldots,g_l)&\mapsto&
[(x_1,\ldots,x_k;t(x_1),g_1,g_1g_2,\ldots,g_1\cdots g_l)].
\end{eqnarray}
It is clear that $\Phi$ establishes an isomorphism
between  $N_k\rtimes \gm_l$ and $(N\rtimes \gm)_{k,l}$
as topological spaces.
Thus it follows that  $N\lcom\rtimes \gm\lcom$
inherits a  $\Delta_2$-structure, which
is the  pull-back of the  $\Delta_2$-structure
on $(N\rtimes \gm)\lcomm$ via $\Phi$.
The proposition below describes  some of the morphisms of
this $\Delta_2$-structure, which are useful in the following
sections.

\begin{prop}\label{prop:ngm}
For any  $f\in \Hom_{\Delta_2}((k,l),(k',l'))$,
let $\tilde{f}: N_{k'}\rtimes \gm_{l'}\to N_k\rtimes \gm_l$,
$(x'_1,\ldots,x'_k;g'_1,\ldots,g'_l)\mapsto (x_1,\ldots,x_k;g_1,\ldots,g_l)$, 
be its corresponding map. Then
\begin{enumerate}
\item if $f=(\emptyset,b,c)$, then $g_i=
(x'_{b(i-1)+1}\cdots x'_{b(i)})^{g_1'\cdots g'_{c(i-1)}}
g'_{c(i-1)+1}\cdots g'_{c(i)}$ and $x_i=t(g'_{c(0)+1})$; and
\item if $f=(a,\emptyset,c)$, then $x_i=x'_{a(i-1)+1}\cdots x'_{a(i)}$ and
$g_i=g'_{c(i-1)+1}\cdots g'_{c(i)}$.
\end{enumerate}
\end{prop}
\begin{pf}
This follows from  an elementary verification, which is left to  the
reader.
\end{pf}

The following is an immediate consequence of Proposition \ref{prop:ngm} (2).

\begin{cor}
The map $\Phi$ is an isomorphism of bi-simplicial spaces, where
the bi-simplicial structure on the left-hand side
is the standard one described above and the bi-simplicial structure on the
right-hand side is induced  from the $\Delta_2$-structure.
\end{cor}

%The map
%$$\Phi: N_k\rtimes \gm_l\lon (N\rtimes \gm)_{k,l}$$
%defined by $(x_1,\ldots,x_k;g_1,\ldots,g_l)\mapsto
%(x_1,\ldots,x_k;t(x_1),g_1,g_1g_2,\ldots,g_1\cdots g_l)$
%is an isomorphism of bi-simplicial spaces, where
%the bi-simplicial structure on the
%right-hand side is induced  from the $\ff\Delta$-structure.
%\end{prop}
%\begin{pf}
%This is an elementary verification. Let us sketch it.
%Using Eqs. (\ref{eqn:hpsi}) and (\ref{eqn:psiu}), we find that the
%$\Delta_2$-structure on $N\lcom\rtimes \gm\lcom$ pulled-back by $\Phi$
%from the one on $(N\rtimes\gm)\lcomm$ is given by
%$\tilde{f}(x'_1,\ldots,x'_{k'};g'_1,\ldots,g'_{l'})=(x_1,\ldots,x_k;
%g_1,\ldots,g_l)$, where:

%\begin{itemize}
%\item[(i)]
%if $f=(\emptyset,b,c)$, then $g_i=
%(x'_{b(i-1)+1}\cdots x'_{b(i)})^{g_1'\cdots g'_{c(i-1)}}
%g'_{c(i-1)+1}\cdots g'_{c(i)}$ and $x_i=t(g'_{c(0)+1})$;
%\item[(ii)]
%if $f=(a,\emptyset,c)$, then $x_i=x'_{a(i-1)+1}\cdots x'_{a(i)}$ and
%$g_i=g'_{c(i-1)+1}\cdots g'_{c(i)}$.
%\end{itemize}
%
%Then (ii) exactly proves the proposition.
%\end{pf}

\subsection{\v{C}ech Cohomology}
Let us now define \v{C}ech cohomology. We refer
the reader to  \cite{tu05} for the particular
case of $\cC=\Delta$; proofs of assertions for  the categories
$\Delta^2$ and $\Delta_2$ are essentially the same and will be omitted.
Given a category $\cC$ equipped with a \DG\ structure $A$,
assume that  we are given a sub-category $\cC'$ such that
\begin{itemize}
\item[(i)] $\cC'$ contains $A$;
\item[(ii)] ${\cC'}^k$ is finite for all $k\in \cC_0$.
\end{itemize}
Note that in this case $A^k$ is necessarily finite for all $k\in \cC_0$; conversely,
if $A^k$ is finite for all $k\in \cC_0$, then the sub-category generated by
$A$ satisfies (i) and (ii) above.
For instance, in the case of $\cC=\Delta$, one can take $\cC'$ to be
the pre-simplicial category $\Delta'$, i.e.
$\Hom_{\Delta'}(k,k')$ consists of (strictly) increasing maps $[k]\to [k']$.
For $\cC=\Delta_2$, $\cC'$ will be the set of degree $\ge 0$ 
morphisms (recall that $\deg(f)=\deg(t(f))-\deg(s(f))$).
The reason why we define $\cC'$ this way is that we need morphisms
$f_\sigma$ (see Eq. (\ref{eqn:fsigma})) to belong to $\cC'$.

An \emph{open cover} of a $\cC$-space $M\lcom$ is a collection $(\uU^k)$,
indexed by $k\in \cC_{0}$ such that $\uU^k=(U^k_i)_{i\in I_k}$
is an open cover of the topological space $M_k$. A $\cC'$-cover
is an open cover, together with a $\cC'$-structure on $I\lcom$ such that
for all $f\in {\cC'}_{1}$ and all $i\in I_{t(f)}$,
$\tilde{f}(U_i)\subset U_{\tilde{f}(i)}$.
Given any open cover, there is a canonical $\cC'$-cover which is finer.
Indeed, let $\Lambda_k$ be the set of families
$\lambda=(\lambda_f)_{f\in {\cC'}_k}$ such that $\lambda_f\in I_{s(f)}$.
Let $V^k_\lambda = \cap_{f\in {\cC'}_k} \tilde{f}^{-1}
(U^{s(f)}_{\lambda_f})$. This is an open subset of $M_k$ since it is
the intersection of finitely many open subsets. Moreover, $\Lambda\lcom$
is endowed with a $\cC'$-structure, by
$(\tilde{h}\lambda)_g=\lambda_{h\smalcirc g}$, $\forall h\in \cC'_1$,
and it is straightforward
to check that the cover $(\sigma\uU^k)$ defined by $\sigma\uU^k:
=(V^k_\lambda)_{\lambda\in\Lambda_k}$, is a $\cC'$-cover, called
the $\cC'$-refinement of  $(\uU^k)$.

\par\medskip
Now given a $\cC'$-cover $(\uU^k)$, let $M'_k=\amalg_{i\in I_k} U^k_i$.
Then $M'\lcom$ is endowed with a $\cC'$-structure.
Moreover, any sheaf $\aA$ on $M\lcom$ induces a sheaf on $M'\lcom$,
again denoted by $\aA$, by $\aA(U)=\prod_{i\in I_k}\aA(U\cap U^k_i)$
($U$ is any open subset of $M'_k$).

Since $\cC'$ is a \DG\ category, one can define $C^*(M'\lcom;\aA)$,
$Z^*(M'\lcom;\aA)$ and $\hh^*(M'\lcom;\aA)$ as in Section \ref{sec:sheaf}.
These groups will be denoted by $C^*(\uU;\aA)$,
$Z^*(\uU;\aA)$ and $H^*(\uU;\aA)$ respectively.
Then the \v{C}ech cohomology groups
$\check{H}^*(M\lcom;\aA)$ are, by definition, the inductive limit of
the groups $H^*(\uU;\aA)$, when $\uU$ runs over all $\cC'$-covers
of $M\lcom$ (the inductive limit being taken in the generalized sense
of limits of functors, since a $\cC'$-cover may be refined to another
by several different ways).

Note that the \v{C}ech cohomology groups do not depend on the choice of $\cC'$
satisfying (i) and (ii) above. Indeed, let $\cC''$ be the category
generated by $A$. Since any $\cC'$-cover is a $\cC''$-cover, and since
any $\cC''$-cover admits a $\cC'$-cover which is finer, it follows
easily that the \v{C}ech cohomology groups defined using $\cC'$
coincide with those defined using $\cC''$.

When $\cC=\Delta$, $\Delta^2$ or $\Delta_2$, the \v{C}ech cohomology
groups can also be seen as the cohomology groups of a ``canonical
\v{C}ech complex'' $\check{C}^*(M\lcom;\aA)$, which is, by definition,
the inductive limit of the complexes $C^*(\sigma\uU;\aA)$,
where $\uU$ runs over  all covers of the form $\uU^k=(U^k_x)_{x\in M_k}$,
and $U^k_x$ is an open neighborhood of $x$; the cover
$\uU'$ is said to be finer than $\uU$ if $(U')^k_x\subset U^k_x$
for all $k$ and $x\in M_k$ (see \cite{tu05} for details).
In the sequel, \v{C}ech cochains (resp. \v{C}ech cocycles)
should be understood in the above sense.

%such that for all $q\ge 0$
%$\aA_q$ is an injective sheaf over $M\lcom$ (or, more generally,
%such that $H^n(M_p;\aA_q)=0$ for all $n\ge 1$ and $p\ge 0$).
%Let $C^{p,q}=(\aA_q(M_p))$ ($p,q\ge 0$). Then $C^{\bullet,\bullet}$
%is endowed of a structure of double complex as follows. Let
%$\del=\sum_{i=0}^p (-1)^i\teps_i^*$, where
%for all $n\ge 0$ and $0\le i \le n$,
%we denote by $\varepsilon_i^n:[n-1]\to [n]$
%the unique increasing map which omits $i$ (the superscript $n$
%will be omitted if there is no ambiguity).
%
%Similarly, if $M\lcomm$ is a bi-simplicial space, let
%$C^{k,l,q}=\aA_q(M_{k,l})$.
%Denote by $\del'=\sum_{i=0}^k(-1)^i\tilde{\varepsilon'}_i$ the
%differential with respect to the first simplicial structure, and
%by $\del=\sum_{i=0}^l (-1)^i\tilde{\varepsilon}_i$ the differential
%with respect to the second simplicial structure. Then
%$H^*(M\lcom;\aA)$ is the cohomology of the complex $(C^{k,l,q})$
%with the differential $\delta=(-1)^{k+l}d+\del'+(-1)^{k}\del$.
%
%Finally, if $M\lcom$ is a $\Delta_2$-space, then we define its
%cohomology as the cohomology of the underlying bi-simplicial
%space.

\section{The transgression maps}

The purpose of this section is to show that there
is a natural  transgression map on the level of cochains
for  the cohomology of a $\Delta_2$-space.
As a consequence, we prove that
for a crossed module   $N\stackrel{\varphi}{\to}\gm$
there exist transgression maps
$$T: \hh^*(\gm\lcom;\aA\lcom)\to \hh^*((N\rtimes \gm)\lcomm;\aA\lcom),$$
and 
$$T_1: \hh^*(\gm\lcom;\aA\lcom)\to \hh^{*-1}((N\rtimes\gm)\lcom;\aA\lcom),$$
and similarly for  \v{C}ech cohomology.
  Throughout this section, $M\lcomm$ denotes a $\Delta_2$-space.

\subsection{Construction of the transgression maps}

For any  fixed $k\in \nn$, consider the restriction of the category
 $\Delta_2$ to the  objects of
the form $(k,l)$ ($l\in\nn$), and to morphisms of the form
$\mbox{Id}\times (f\times f):
\ff_k\times [l]\times [l]\to \ff_k\times [l']\times [l']$,
where $f:[l]\to [l']$ is non-decreasing.
This category is isomorphic to $\Delta$. Hence we obtain
 a simplicial space
$M_{k,\bullet}=(M_{k,l})_{l\in\nn}$.
%In particular, if $N\to\gm$ is a crossed module, denote by $M\lcomm$
%the associated $\Delta_2$-space. Let $N_k
%=\{(x_1,\ldots,x_k)\in N^k\vert\; s(x_1)=\cdots=s(x_k) (=t(x_1)=\cdots
%=t(x_k))\}$. Then $N_k$ is endowed with an action of $\gm$,
%and the simplicial space associated to the crossed-product groupoid
%$N_k\rtimes \gm$ is precisely $M_{k,\bullet}$. As a special case,
%for $k=0$ we get $M_{0,\bullet}=\gm\lcom$ and for $k=1$ we get
%$M_{1,\bullet}=(N_1\rtimes\gm)\lcom$.
%When $N=S\gm$, $M_{1,\bullet}=(S\gm\rtimes \gm)\lcom=(\Lambda\gm)\lcom$.
%Denote, in particular,  by $\gm\lcom$ the simplicial space $M_{0,\bullet}$.
Similarly, $M_{\bullet, k}=(M_{l, k})_{l\in\nn}$ is also a 
simplicial space.

Let
%\begin{equation}
%\label{eqn:complex-sheaves}
$\aA_0\stackrel{d}{\to}\aA_1\stackrel{d}{\to}\cdots$
%\end{equation}
be a  complex of abelian
sheaves over $M\lcomm$, and   $C^*(M\lcomm;\aA\lcom)$
(resp. $C^*(\M{0};\aA\lcom)$)
 its  associated cochain
complex 
on $M\lcomm$ (resp. $\M{0}$) as in Section \ref{sec:sheaf}.
 The main goal of this section
is to  construct  chain maps 
$T: C^*(\M{0};\aA\lcom)\to C^*(M\lcomm;\aA\lcom)$
and $T_1:  C^*(\M{0};\aA\lcom)\to C^{*-1} (\M{1}; \aA\lcom)$.
Thus
 we  construct   natural transgression
maps $T: \hh^*(\M{0};\aA\lcom)\to \hh^* (M\lcomm;\aA\lcom)$,
 and $T_1:  \hh^*(\M{0};\aA\lcom)\to \hh^{*-1} (\M{1}; \aA\lcom)$.
As a consequence, for  an abelian  sheaf $\aA$ over $M\lcomm$, we
obtain transgression
maps  $T: H^*(\M{0};\aA)\to H^*(M\lcomm;\aA)$ and
$T_1:  H^*(\M{0};\aA)\to H^{*-1} (\M{1}; \aA)$.
 Similarly, there are chain maps
$T: \check{C}^*(\M{0};\aA)\to \check{C}^{*}(M\lcomm;\aA)$
and $T_1:  \check{C}^*(\M{0};\aA)\to \check{C}^{*-1} (\M{1}; \aA)$,
which induce transgression maps for \v{C}ech cohomology
$T: \check{H}^*(\M{0};\aA)\to \check{H}^{*}(M\lcomm;\aA)$ and
$T_1:  \check{H}^*(\M{0};\aA)\to \check{H}^{*-1} (\M{1}; \aA)$
as well.

We  first give the  construction for sheaf cohomology.
First of all, we  need to introduce some notations.

Recall that a $(k,l)$-shuffle is a permutation $\sigma$ of
$\{1,\ldots,k+l\}$ such that $\sigma(1)<\cdots<\sigma(k)$ and
$\sigma(k+1)<\cdots<\sigma(k+l)$. One can represent a $(k,l)$-shuffle
by a sequence of balls, $k$ of which being black and $l$ of which being white.
More precisely, $B_\sigma=\sigma(\{1,\ldots,k\})$ and
$W_\sigma=\sigma(\{k+1,\ldots,k+l\})$. The signature of $\sigma$
can be easily computed by the formula
$\varepsilon(\sigma)=(-1)^{\sum_{i\le k} \sigma(i)-i}$.

Denote by $S_{k,l}$ the set of $(k,l)$-shuffles.
For any  $\sigma\in S_{k,l}$, we define
$f_\sigma \in\Hom_{\Delta_2}((0,k+l),(k,l))$ by
\begin{equation}\label{eqn:fsigma}
f_\sigma=(\emptyset,b^\sigma,c^\sigma),
\end{equation}
%where $0$ stands for the zero map $[0]\to [k]$,
where $b^\sigma$ is the map $[k+l]\to [k]$  given
by $b^\sigma(i)=\#(B_\sigma\cap\{1,\ldots,i\})$,
and $c^\sigma$ is the map $[k+l]\to [l]$
given by $c^\sigma (i)=\#(W_\sigma\cap\{1,\ldots,i\})$, $i=0, \ldots, k+l$.
Thus   $f_\sigma$   induces $\tilde{f}_\sigma: M_{k,l}\to M_{0, k+l}$.
%\begin{tiny}
%\begin{numex}
%In the case of the $\Delta_2$-space associated to a crossed module,
%if $\sigma$ is the $(2,2)$-shuffle $(1,3,2,4)$, then
%$\tilde{f}_\sigma(x_1,x_2,g_1,g_2)=(\varphi (x_1) ,g_1, \varphi (x_2)^{g_1},
%g_2)$.
%\end{numex}
%\end{tiny}
Therefore, we obtain  a map $\tilde{f}^*_\sigma:
\aA^q(M_{0, k+l})\to \aA^q(M_{k,l})$. Taking
the direct sum over all $l$ and $q$, we obtain a map
$$\tilde{f}_\sigma^*: C^*(\M{0};\aA\lcom)\to C^{*-k}(M_{k,\bullet};\aA\lcom).$$

Set 
$$T_k=\sum_{\sigma\in S_{k,l}} \varepsilon(\sigma)\tilde{f}^*_\sigma:
 C^*(\M{0};\aA\lcom)\to C^{*-k}(M_{k,\bullet};\aA\lcom)
$$
(with $T_0=\mbox{Id}$), and
$$T=\sum_{k\ge 0} T_k:  C^*(\M{0};\aA\lcom)\to C^*(M\lcomm;\aA\lcom) $$
 using the decomposition 
$C^*(M\lcomm;\aA\lcom)=\oplus_{k\ge 0} C^{*-k}(M_{k,\bullet};\aA\lcom)$.

For any fixed $k\geq 0$, by $\partial$ and $\partial'$ we denote the
differentials
$$\partial: C^{p, q}(\M{k}, \aA\lcom) (=\aA_q (M_{k, p}))
\lon C^{p+1, q}(\M{k}, \aA\lcom) (=\aA_q (M_{k, p+1})), $$
and 
$$\partial': C^{p, q}(\Ma{k}, \aA\lcom) (=\aA_q (M_{p, k}))
\lon C^{p+1, q}(\Ma{k}, \aA\lcom) (=\aA_q (M_{p+1, k})), $$
respectively, and by  $\delta_k=(-1)^pd+\del$, 
we denote the  total differential of the double complex
$C^{p,q}(M_{k,\bullet};\aA\lcom)$.
Note that  $C^*(M\lcomm;\aA\lcom)= \sum_{p+q+k=*} \aA_q (M_{k, p})$
and the total differential is $(-1)^{k+p}d+ (-1)^k\partial+\partial'$.

\begin{lem}
\label{lem:3.2}
Assume that  $M\lcomm$ is a $\Delta_2$-space, and
 $\aA\lcom$ is  a complex of abelian
sheaves over $M\lcomm$. Then
\begin{eqnarray}
\label{eqn:delT}
\del'T_k &=& T_{k+1}\del+(-1)^k\del T_{k+1} \ \ \mbox{ and}\\
\label{eqn:deltaT}
\del'T_k &=& T_{k+1}\delta_0+(-1)^k\delta_{k+1} T_{k+1},
\end{eqnarray}
where  both sides are maps from $C^* (\M{0}; \aA\lcom )$
to $C^{*-k} (\M{k+1}; \aA\lcom )$.
% FIND WHAT IS THE EXACT ASSUMPTION HERE FOR CECH COHOMOLOGY, 
%\begin{equation}
%\del'T_k =T_{k+1}\del+(-1)^k\del T_{k+1}
%\end{equation}
\end{lem}
\begin{pf}
Let us first show that Eq. \eqref{eqn:deltaT}
 follows from Eq. \eqref{eqn:delT}.
For any  $\omega\in C^{k+l,q}(\M{0}, \aA\lcom )$, we have
\be
T_{k+1}\delta_0 \omega+(-1)^{k}\delta_{k+1} T_{k+1}\omega
&=&T_{k+1}((-1)^{k+l} d +\del)\omega + (-1)^k ((-1)^{l-1} d+\del)T_{k+1}\omega\\
&=&T_{k+1}\del\omega+(-1)^k\del T_{k+1}\omega.
\ee
Here we have
 used the fact that $T_{k+1}$ commutes with $d$ since 
$T_{k+1}$ is
a summation of pull-back maps.

Now we prove Eq. \eqref{eqn:delT}.  For any 
$\omega\in C^{k+l,q}(\M{0}, \aA\lcom )$,
$$T_{k+1}\del \omega =\sum_{\sigma\in S_{k+1,l}}\sum_{j=0}^{k+l+1}
(-1)^j \varepsilon(\sigma) f_\sigma^*\teps_j^*\omega.$$
In the sum above, we distinguish three cases:

1) ($j=0$, $\sigma^{-1}(1)\le k+1$)\footnote{thus, $\sigma^{-1}(1)=1$}
or ($j=p$, $\sigma^{-1}(k+l+1)\le k+1$)\footnote{thus,
$\sigma^{-1}(k+l+1)=k+1$}
or ($1\le j\le k+l$, $\sigma^{-1}(j)\le k+1$, $\sigma^{-1}(j+1)\le k+1$);

2) ($j=0$, $\sigma^{-1}(1)\ge k+2$)\footnote{thus, $\sigma^{-1}(1)=k+2$},
($j=k+l+1$, $\sigma^{-1}(k+l+1)\ge k+2$)\footnote{thus,
$\sigma^{-1}(k+l+1)=k+l+1$}
or ($1\le j \le k+l$, $\sigma^{-1}(j)\ge k+2$, $\sigma^{-1}(j+1)\ge k+2$);

3) $1\le j \le k+l$ and
\begin{itemize}
\item[(a)] either $\sigma^{-1}(j)\le k+1$ and $\sigma^{-1}(j+1)\ge k+2$
\item[(b)] or $\sigma^{-1}(j)\ge k+2$ and $\sigma^{-1}(j+1)\le k+1$.
\end{itemize}

We  show that the terms in 1) are equal to $\del'T_k\omega$,
the terms in 2) are equal to $(-1)^{k+1}\del T_{k+1}$ and the terms in
3a) cancel out with those in 3b).
\par\medskip

Let us examine the terms in 1). We have
\begin{eqnarray*}
\del'T_k\omega &=&\sum_{m=0}^{k+1}(-1)^m{{\tepss}}^{*}_m T_k\omega\\
&=&\sum_{m=0}^{k+1}\sum_{\tau\in S_{k,l}} (-1)^m \varepsilon(\tau)
{{\tepss}}^{*}_m f_\tau^*\omega.
\end{eqnarray*}

Given $(j,\sigma)$ as in 1), we define $(m,\tau)$ as follows:
$m=\sigma^{-1}(j)$, with the convention $\sigma^{-1}(0)=0$,
and $\tau$ is uniquely determined by the equation
$$\varepsilon_j\smalcirc\tau = \sigma\smalcirc \varepsilon_m.$$
In other words, if the shuffle $\sigma$ is represented by a sequence
of $p=k+l+1$ balls, $k+1$ of which being black and $l$ of which being
white,
then $\tau$ is obtained from $\sigma$ by removing the $j$-th
one (which is black).

We need to check the following equalities:
\begin{itemize}
\item[(i)] $(-1)^m\varepsilon(\tau)=(-1)^j\varepsilon(\sigma)$, and
\item[(ii)] $\tilde{f}_\sigma^*\teps_j^*\omega=
{{\tepss}}^{*}_m\tilde{f}_\tau^*\omega$.
\end{itemize}

To show (i), let $p=k+l+1$ and 
$\sigma_{j,p}$  the circular permutation
$(j,j+1,\ldots,p)$. Then
$\varepsilon_p\smalcirc \tau = \sigma_{j,p}^{-1}\smalcirc\varepsilon_j\smalcirc\tau
=\sigma_{j,p}^{-1}\smalcirc\sigma\smalcirc\varepsilon_m
=(\sigma_{j,p}^{-1}\smalcirc\sigma\smalcirc\sigma_{m,p})\smalcirc\varepsilon_p$.
Thus $\varepsilon(\tau)=\varepsilon(\sigma_{j,p}^{-1}
\smalcirc\sigma\smalcirc \varepsilon_{m,p})=(-1)^{j-m}\varepsilon(\sigma)$.

To show (ii), it suffices to prove that
$f_\sigma\smalcirc \varepsilon_j={\varepsilon'_m}\smalcirc f_\tau$.
Now
\begin{eqnarray*}
f_\sigma\smalcirc\varepsilon_j&=&(0,b_\sigma,c_\sigma)\smalcirc
(\mbox{Id},0,\varepsilon_j)
 = (0,b_\sigma\smalcirc\varepsilon_j,c_\sigma\smalcirc\varepsilon_j),  \  \ \mbox{ and}\\
\varepsilon'_m\smalcirc f_\tau &=& (\varepsilon_m,0,\mbox{Id})
\smalcirc (0,b_\tau,c_\tau) = (0,\varepsilon_m\smalcirc b_\tau,c_\tau).
\end{eqnarray*}
Hence it remains to check that $b_\sigma\smalcirc\varepsilon_j
=\varepsilon_m\smalcirc b_\tau$ and $c_\sigma\smalcirc\varepsilon_j=
c_\tau$, i.e. that

\begin{eqnarray*}
\#(B_\sigma\cap\{1,\ldots,\varepsilon_j(i)\})
  &=&\varepsilon_m(\#(B_\tau\cap\{1,\ldots,i\})), \ \ \mbox{and}  \\
\#(W_\sigma\cap\{1,\ldots,\varepsilon_j(i)\})
  &=&\#(W_\tau\cap\{1,\ldots,i\}),
\end{eqnarray*}
which is immediate from the description of $B_\tau$ and $W_\tau$
in terms of $B_\sigma$ and $W_\sigma$.
\par\bigskip

The terms in 2) are treated in a similar fashion: define $(m,\tau)$
by $m=\sigma^{-1}(j)$ with the convention $\sigma^{-1}(0)=k+1$,
and $\tau$ by the equation
$\varepsilon_j\smalcirc\tau = \sigma\smalcirc\varepsilon_m$.
\par\medskip
For the terms in 3), it is easy to check that the term corresponding to
$(j,\sigma)$ cancels out with $(j,\tau_{j,j+1}\smalcirc\sigma)$,
where $\tau_{j,j+1}$ is the transposition which exchanges $j$ and $j+1$.
\end{pf}

One can introduce transgression maps for \v{C}ech cohomology
in a similar fashion. Namely,
if $\uU:=(\uU^{k, l})$ is a $\Delta_2'$-cover of $M\lcomm$,
then for any fixed $k$,  $(\uU^{k, l})$ is a pre-simplicial
cover of $\M{k}$, and  $(\uU^{l,k})$ is a pre-simplicial
cover of $M_{\bullet, k}$. They are denoted, respectively,
by $\uU^{k, \bullet}$ and $\uU^{\bullet, k}$.

Let $M'_{k,l}=\amalg_{i\in I_{k, l}} U^{k, l}_i$.
Then $M'\lcomm$ is endowed with a $\Delta'_2$-structure.
Hence for any fixed $k$, $M'_{k, \bullet}$ and
$M'_{\bullet, k}$ are pre-simplicial spaces.
For any  $\sigma\in S_{k,l}$, since   $f_\sigma
\in\Hom_{\Delta_2}((0,k+l),(k,l))$,
one has  a map   $\tilde{f}_\sigma: M'_{k,l}\to M'_{0, k+l}$. 
Thus $\tilde{f}_\sigma^*: \aA (M'_{0, k+l})\to \aA (M'_{k,l} )$.
Set 
$$T_k=\sum_{\sigma\in S_{k,l}} \varepsilon(\sigma)\tilde{f}^*_\sigma:
\check{C}^*(\uU^{0, \bullet}, \aA) (= {C}^*(\Mr{0};\aA))\to
\check{C}^*(\uU^{k, \bullet}, \aA)(= C^{*-k}(M'_{k,\bullet};\aA))
$$
(with $T_0=\mbox{Id}$), and 
$$T=\sum_{k\ge 0} T_k:  \check{C}^*(\uU^{0, \bullet}, \aA)(= {C}^*(\Mr{0};\aA))
\to  \check{C}^*(\uU, \aA) (=C^*(M'\lcomm;\aA)), $$
 using the decomposition
$C^*(M'\lcomm;\aA)=\oplus_{k\ge 0} C^{*-k}(M'_{k,\bullet};\aA)$.
For any fixed $k\geq 0$, by $\partial$ and $\partial'$,  
 we denote,  respectively,  the differentials
$$\partial: \check{C}^p(\uU^{k, \bullet}, \aA) (=\aA (M'_{k, p}))
\lon \check{C}^{p+1}(\uU^{k, \bullet}, \aA) (=\aA (M'_{k, p+1})), $$
and
$$\partial': \check{C}^p(\uU^{ \bullet , k}, \aA)  (=\aA (M'_{p, k}))
\lon \check{C}^{p+1}(\uU^{\bullet, k}, \aA)  (=\aA (M'_{p+1, k})). $$
Note that  $ \check{C}^*(\uU, \aA)=
 C^*(M'\lcomm;\aA)= \sum_{p+k=*} \aA (M'_{k, p})$,
and the total differential is $ \partial+\partial'$.

\begin{lem}
\label{lem:3.2-1}
Assume that  $M\lcomm$ is a $\Delta_2$-space.
If $\aA$ is  an  abelian sheaf over $M\lcomm$ and $(\uU^{k, l})$ is a
 $\Delta_2'$-cover of $M\lcomm$, then
\begin{equation}
\label{eq:10}
\del'T_k =T_{k+1}\del+(-1)^k\del T_{k+1},
\end{equation}
where  both sides are maps from $\check{C}^*(\uU^{0, \bullet}, \aA)$
to $\check{C}^{*-k}(\uU^{k+1, \bullet}, \aA)$.
\end{lem}

As an immediate consequence of Lemma \ref{lem:3.2}-\ref{lem:3.2-1}, we  obtain the following

\begin{them}
\label{thm:transgression}
Assume that  $M\lcomm$ is a $\Delta_2$-space.
\begin{enumerate}
 \item If $\aA\lcom$ is  a  complex of abelian
sheaves over $M\lcomm$, then
\begin{equation}
T: C^*(\M{0};\aA\lcom)\to C^{*} (M\lcomm;\aA\lcom)
\end{equation}
 is a chain map. Therefore it induces a morphism
on the level of cohomology
\begin{equation}
T: \hh^*(\M{0};\aA\lcom)\to \hh^*(M\lcomm;\aA\lcom ).
\end{equation}
 Similarly, given an abelian sheaf $\aA$ over $M\lcomm$,  then
\begin{equation}
T: \check{C}^*(\M{0};\aA)\to \check{C}^{*}(M\lcomm;\aA)
\end{equation}
is a chain map, and therefore it induces a morphism
\begin{equation}
T: \check{H}^*(\M{0};\aA)\to \check{H}^{*}(M\lcomm;\aA). 
\end{equation}
\end{enumerate}
\end{them}
We call $T$  the {\em total transgression map}.
\medskip

Considering the special  case $k=1$ in Lemma  \ref{lem:3.2} and
Lemma \ref{lem:3.2-1}, we immediately obtain  the following

\begin{them}
\label{thm:T1}
Under the same hypothesis as in Theorem \ref{thm:transgression},
\begin{enumerate}
 \item if $\aA\lcom$ is  a  complex of abelian
sheaves over $M\lcomm$, then the map
\begin{equation}
T_1: C^*(\M{0};\aA\lcom)\to C^{*-1} (\M{1};\aA\lcom)
\end{equation}
satisfies the relations:
\begin{eqnarray}\label{eqn:T1}
T_1\del+\del T_1=0, \\ \nonumber
T_1\delta+\delta T_1 =0;
\end{eqnarray}
and thus induces a  morphism:
$$T_1: \hh^*(\M{0};\aA\lcom)\to \hh^{*-1}(M_{1,\bullet};\aA\lcom). $$
\item Similarly,  for any abelian sheaf $\aA$ over  $M\lcomm$,
\begin{equation}
T_1: \check{C}^*(\M{0};\aA)\to \check{C}^{*-1} (\M{1};\aA)
\end{equation}
satisfies 
\begin{equation}
\label{eq:T1par}
T_1\del+\del T_1=0,
\end{equation}
and therefore we have a  morphism:
$$T_1: \check{H}^*(\M{0};\aA)\to \check{H}^{*-1}(M_{1,\bullet};\aA). $$
\end{enumerate}
\end{them}
We call $T_1$ the  {\em transgression map}.

\subsection{Transgression maps for crossed modules}
\label{sec:3.2}

Let  $N\stackrel{\phi}{\to}\gm$ be  a crossed module.
Denote by $M\lcomm$
its associated $\Delta_2$-space $(N\rtimes \gm )\lcomm$. 
Let $N_k
=\{(x_1,\ldots,x_k)\in N^k\vert\; s(x_1)=\cdots=s(x_k) (=t(x_1)=\cdots
=t(x_k))\}$. Then $N_k$ is endowed with a $\gm$-action,
and the simplicial space associated to the crossed-product groupoid
$N_k\rtimes \gm$ is precisely $M_{k,\bullet}$. As a special case,
for $k=0$ we get $M_{0,\bullet}=\gm\lcom$ and for $k=1$ we get
$M_{1,\bullet}=(N\rtimes\gm)\lcom$.
When $N=S\gm$, $M_{1,\bullet}=(S\gm\rtimes \gm)\lcom=(\Lambda\gm)\lcom$.
According to Theorem \ref{thm:transgression},  for
any  complex of  sheaves $\aA\lcom$ over $(N\rtimes \gm )\lcomm$,
 we have the  total transgression  map
$$T: \hh^*(\gm\lcom;\aA\lcom)\to \hh^*((N\rtimes \gm)\lcomm;\aA\lcom).$$
In particular, we have
$$T: H_{dR}^*(\gm\lcom)\to H^*_{dR}((N\rtimes \gm)\lcomm).$$
 Recall that $T=\sum_k T_k$, where
$$T_k=\sum_{\sigma\in S_{k,l}} \varepsilon(\sigma) \tilde{f}_\sigma^*:
C^{k+l}(\gm\lcom; \aA\lcom)
\lon C^l((N_k\rtimes\gm)\lcom; \aA\lcom). $$
Using Proposition~\ref{prop:ngm}, we obtain the following explicit
formula for $\tilde{f}_\sigma$:

\begin{prop}
Let  $N\stackrel{\phi}{\to}\gm$ be  a crossed module.
Then $\forall \sigma \in S_{k, l}$, the map
$\tilde{f}_\sigma: N_k\rtimes \gm_l \to \gm_{k+l}$
is given by
\begin{equation}
\label{eqn:fsigma2a}
\tilde{f}_\sigma(x_1,\ldots,x_k;g_1,\ldots,g_l)=(u_1,\ldots,u_{k+l}),
\end{equation}
where $u_i=g_{\sigma^{-1}(i)}$ if $\sigma^{-1}(i)\ge k+1$, and
$u_i=\varphi\left(x_{\sigma^{-1}(i)}^{\prod_{\sigma^{-1}(j)>k, j< i}
g_{\sigma^{-1}(j)}}\right)$ otherwise.
\end{prop}

\begin{numex}
If $\sigma$ is the $(2,2)$-shuffle $(1,3,2,4)$, then
$$\tilde{f}_\sigma(x_1,x_2,g_1,g_2)=(\varphi (x_1) ,g_1, \varphi (x_2)^{g_1},
g_2). $$
\end{numex}

We also have the transgression map:
$$T_1: \hh^{*}(\gm\lcom;\aA\lcom)\to  \hh^{*-1}((N\rtimes \gm)\lcom;\aA\lcom). $$
In particular, for any groupoid
$\gm\toto \gm_0$, we have  the transgression map
$$T_1: \hh^*(\gm\lcom;\aA\lcom)\to \hh^{*-1}((\Lambda\gm)\lcom;\aA\lcom).$$

\begin{prop}
Let  $N\stackrel{\phi}{\to}\gm$ be  a crossed module.
Then the transgression
$T_1: \hh^{*}(\gm\lcom;\aA\lcom)\to \hh^{*-1}((N\rtimes \gm)\lcom;\aA\lcom)$
is given, on the cochain level,  by 
$$T_1=\sum_{i=0}^{p-1} (-1)^i\tilde{f}_i^* : 
\aA^q(\gm_p)\to \aA^q(N\rtimes \gm_{p-1}),$$
where the map  $\tilde{f}_i: N\rtimes \gm_{p-1}\to \gm_p$ is
given by
\begin{equation}\label{eqn:T1dRa}
\tilde{f}_i(x;g_1,\ldots,g_{p-1})=(g_1,\ldots,g_i,
\varphi(x)^{g_1\cdots g_i} ,g_{i+1},\ldots,g_{p-1}).
\end{equation}
\end{prop}

\begin{rmk}
In the case of de Rham cohomology,  the transgression
$T_1:\Omega^q(\gm_p)\to \Omega^q(N\rtimes \gm_{p-1})$ is
defined by $T_1=\sum_{i=0}^{p-1} (-1)^i\tilde{f}_i^*$,
where $\tilde{f}_i$ is given  by Eq. \eqref{eqn:T1dR}.
More generally, the map
$T_k:\Omega^q(\gm_p)\to \Omega^q(N_k\rtimes \gm_{p-k})$
takes the form $T_k=\sum_{\sigma\in S_{k,p-k}}\varepsilon(\sigma)
\tilde{f}_\sigma^*$, where  $\tilde{f}_\sigma$ is
given   by Eq.  (\ref{eqn:fsigma2}).
\end{rmk}

\begin{rmk}
Consider the case that $\gm$ is a Lie  group $G$.
Then $\Lambda G$ is the transformation groupoid $G\rtimes G\toto G$,
where $G$ acts on $G$ by conjugation.
On the cochain level, the total transgression map  
 for the de Rham cohomology becomes
\begin{equation}
T: \Omega^n_{dR}( G\lcom )\lon \Omega^n_{dR}((G\rtimes G)\lcomm),
\end{equation}
and  the  transgression map is
\begin{equation}
T_1: \Omega^n_{dR} ( G\lcom )\lon \Omega^{n-1}_{dR}((G\rtimes G)\lcom).
\end{equation}
On the other hand, there is a Bott-Shulman map \cite{Bott, BSS}
$S(\Gg^*)^G\to  \Omega^n_{dR} ( G\lcom )$. We have the
following

\begin{quote}
{\bf  Conjecture} 
The equivariant Bott-Shulman map $\Phi$ constructed by Jeffrey
\cite{Jeffrey, Meinrenken}  are compositions of the
usual Bott-Shulman map with the total transgression map
$T$,  under which $T_1$ becomes $\Phi_1$ in \cite{Jeffrey, Meinrenken}.
\end{quote} 

%Indeed, our $T_k$ can be thought of the simplicial  counterpart
%of $\Phi_k$ in Section 5.3  of \cite{Jeffrey} (see also Section B8
%of \cite{Meinrenken}). However our $T_k$ is constructed
%on the level of cochains, while $\Phi_k$ is only
%constructed on the level of cocycles (which are elements of $S(\Gg^*)^G$).
%In fact, one sees that Eq. \eqref{eqn:deltaT}
% reduces to Eqs. (5.7-5.9) in  \cite{Jeffrey} when applying to a cocycle
%of $\delta_0$.
More precisely, the composition of the  Bott-Shulman map
$S^p(\Gg^*)^G\to  \Omega^{2p}_{dR} ( G\lcom )$ with  $T_k:
\Omega^{2p}_{dR} ( G\lcom )\to \Omega^{2p-k}_{dR} ((G\rtimes G^k)\lcom )$
yields a  map 
$S^p(\Gg^*)^G\to\Omega^{2p-k}_{dR} ((G\rtimes G^k)\lcom )$,
which should be the simplicial counterpart of the maps $\Phi_G^{(k)}$
in Section B8 of \cite{Meinrenken}. 
\end{rmk}

To end this subsection, we  give  
a geometric description of the transgression map
$T_1: H^2(\gm\lcom;\sS^1)\to H^1(\Lambda\gm\lcom;\sS^1)$,
which was introduced in \cite{TX}. 

Given a topological groupoid $\gm$ and an element
 $\alpha\in H^2(\gm\lcom;\sS^1)$,
represent  $\alpha$
  by  an $S^1$-central extension
\begin{equation}\label{eqn:centralexth2}
S^1\to \tilde{\gm'}\to \gm'.
\end{equation}
Let $L=\tilde{\gm'}\times_{S^1}\cc$ be the complex line bundle associated
to the $S^1$-principal bundle $\tilde{\gm'}\to \gm'$,  and
 $L'\to S\gm'$ its restriction  to the subspace $S\gm'$. Then $L'\to S\gm'$
is a $\gm'$-equivariant bundle,  where the $\gm'$-action 
is given by  $\gamma\cdot \xi = \tilde{\gamma}\xi\tilde{\gamma}^{-1}$,
for all compatible $\gamma \in \gm'$ and $\xi \in L'$. Here
 $\tilde{\gamma}\in\tilde{\gm'}$ is any lift of $\gamma$. Thus
$L'\to S\gm'$ is an $S^1$-bundle over $\Lambda\gm'\toto S\gm'$,
and hence it determines an element in
 $H^1(\Lambda\gm';\sS^1)$.
Since $\gm'$ is Morita equivalent to $\gm$, $\Lambda\gm'$ is Morita equivalent
to $\Lambda\gm$. Thus we obtain an element in
 $H^1(\Lambda\gm;\sS^1)$, which is denoted by $\tau(\alpha)$.
It is a simple exercise to check that the  element $\tau(\alpha)$
does not depend on
the choice of the $S^1$-central extension (\ref{eqn:centralexth2}).
Hence we obtain a  map
$$\tau: H^2(\gm\lcom;\sS^1)\lon H^1(\Lambda\gm\lcom;\sS^1). $$

\begin{prop}
The map $T_1$ coincides with the opposite of the map $\tau$.
\end{prop}

\begin{pf}
Represent an element $\alpha\in H^2(\gm\lcom;\sS^1)$
by an $S^1$-central extension as in Eq. (\ref{eqn:centralexth2}) above.
Without loss of generality, we may assume
 that $\gm'=\gm$. Let $(U_i)$ be an open cover of $\gm$
such that  $\tilde{\gm}\to \gm$ admits a section over each  $U_i$,
denoted by $\varphi_i: U_i\to \tilde{\gm}$. Then the \v{C}ech
2-cocycle $c_{ijk} \in \check{C}^2 (\uU, \sS^1)$ associated to
the central extension is defined by the equation
\begin{equation}\label{eqn:cijk}
\varphi_i(g)\varphi_j(h)=\varphi_k(gh)c_{ijk}(g,h)
\end{equation}
($g\in U_i$, $h\in U_j$, $gh\in U_{k}$).
Therefore, for all $x\in S\gm$ and $g\in \gm$ such that $t(g)=s(x)$,
$x\in U_i$, $g\in U_j$, $xg\in U_k$, $x^g \in U_l$, we have
\begin{eqnarray}
\label{eqn:cijk1}
\varphi_i(x)\varphi_j(g)&=&\varphi_k(xg)c_{ijk}(x,g)\\
\label{eqn:cijk2}
\varphi_j(g)\varphi_l(x^g)&=&\varphi_k(xg)c_{jlk}(g,x^g).
\end{eqnarray}
Comparing Eq. (\ref{eqn:cijk1}) and Eq. (\ref{eqn:cijk2}), we get
\begin{equation}\label{eqn:cpsi}
\varphi_j(g)\varphi_l(x^g)\varphi_j(g)^{-1}=\varphi_i(x)
\frac{c_{jlk}(g,x^g)}{c_{ijk}(x,g)}.
\end{equation}
The quantity $\frac{c_{jlk}(g,x^g)}{c_{ijk}(x,g)}$ thus does not depend
on $j$ and $k$. Let us denote it by $\psi_{il}(x,g)$. Then Eq.
 (\ref{eqn:cpsi})
reads
$$\gamma\cdot\varphi_l(s(\gamma))=\varphi_i(t(\gamma))\psi_{il}(\gamma)$$
for all $\gamma\in(\Lambda\gm)_{U_l}^{U_i}$, where $(\gamma,p)\mapsto
\gamma\cdot p$ denotes the action of $\gm$ on the principal bundle
$S\tilde{\gm}$. This shows that $\tau[c]=[\psi]$.

On the other hand, $(T_1c)_{il}(x,g)=\frac{c_{ijk}(x,g)}{c_{jkl}(g,x^g)}
=\psi_{il}(x,g)^{-1}$. Hence $T_1= -\tau$. 
\end{pf}

\subsection{Multiplicative cochains}
Following the notations of Section \ref{sec:3.2}, we assume
that $N\stackrel{\varphi}{\to}\gm$ is a crossed module. 

\begin{defn}
Let $\aA$ be an abelian sheaf over the bi-simplicial space
$N\lcom\rtimes\gm\lcom$. A
cochain $c\in \check{C}^n((N_k\rtimes \gm)\lcom ;\aA)$ is said
to be $r$-multiplicative ($r\ge 1$) if $\del' c=\del b$ for some
$(r-1)$-multiplicative $b\in \check{C}^{n-1}((N_{k+1}\rtimes \gm)\lcom; \aA)$.
For $r=0$, every cochain is said to be 0-multiplicative.

In other words, $c$ is $r$-multiplicative if
$\del'c$ is a coboundary for the
total complex $(\check{C}^p((N_{k+q}\rtimes \gm)\lcom;\aA),\del,\del')_{p\ge 0, 1\le q\le r}$.
\end{defn}

 From now on, we assume that $\aA$ is a sheaf over the $\Delta_2$-space
$(N\rtimes\gm)\lcomm$, and therefore it is also a sheaf over
the bi-simplicial space $N\lcom \rtimes\gm\lcom$.
The following
 lemma  shows that the maps $T_k$ enable us to produce many
multiplicative cochains.

\begin{lem}
Assume  that $\omega\in \check{C}^n(\gm\lcom;\aA)$ satisfies $\del\omega=0$.
Then,  for all $k\ge 1$, $T_k\omega\in 
\check{C}^{n-k}((N_k\rtimes \gm)\lcom, \aA)$ is $\infty$-multiplicative.
\end{lem}

\begin{pf}
We show,  by induction on $r$,
 that $T_k\omega$ is $r$-multiplicative for all
$k$ and all $\omega$ as in the lemma. For $r=0$,  this is obvious.
Let $r\ge 1$, and assume that  the assertion is true up to $(r-1)$. 
According to Lemma \ref{lem:3.2-1}, $\del'T_k\omega = T_{k+1}\del\omega
+(-1)^k \del T_{k+1}\omega = \pm \del T_{k+1}\omega$.
Since $T_{k+1}\omega$
is $(r-1)$-multiplicative by the induction assumption, it follows that
$T_k\omega$ is $r$-multiplicative.
\end{pf}

\begin{rmk}
\label{rmk:3.10}
It is easy to check that $c$ is $r$-multiplicative if and only if
$c+\partial c'$ is $r$-multiplicative. Hence  one can talk about
$r$-multiplicative cohomology classes. Therefore, any element in the
image of $T_1: \check{H}^n(\gm\lcom;\aA)\to \check{H}^{n-1}((N\rtimes\gm)\lcom;\aA)$ is $\infty$-multiplicative.
In particular, for any groupoid $\gm$, any element in the image 
of the transgression map  of $T_1: \check{H}^n(\gm\lcom;\aA)\to
 \check{H}^{n-1}(\Lambda \gm\lcom;\aA)$ is $\infty$-multiplicative.
\end{rmk}

For  $\aA=\sS^1$ and $n=3$, we  are led to the following
\begin{cor}
\label{cor:3.11}
Let $N\stackrel{\varphi}{\to}\gm$ be a crossed module.
For any $e\in \check{Z}^3(\gm\lcom;\sS^1)$, let    $c=T_1 e \in
\check{C}^2((N\rtimes \gm)\lcom;\sS^1)$,
 $b=-T_2 e\in \check{C}^1((N_2\rtimes \gm)\lcom;\sS^1)$ and $a=-T_3e\in
\check{C}^0((N_3\rtimes \gm)\lcom;\sS^1)$. Then we have
\begin{equation}
\label{eq:multiplicator}
\del c=0, \ \ \del'c = \del b, \ \ \mbox{ and }  \del' b = \del a.
\end{equation}
\end{cor}
Such a triple $(c, b, a)$ is called a multiplicator.

\begin{defn}
\label{def:multiplicator}
For a crossed module $N\stackrel{\varphi}{\to}\gm$, a {\em multiplicator}
is a  triple $(c, b, a)$, where $c\in
\check{C}^2((N\rtimes \gm)\lcom;\sS^1)$,
 $b\in \check{C}^1((N_2\rtimes \gm)\lcom;\sS^1)$ and $a\in
\check{C}^0((N_3\rtimes \gm)\lcom;\sS^1)$ such that 
Eq. \eqref{eq:multiplicator} holds.
\end{defn}

\begin{rmk}
Multiplicators are closely related to   ``multiplicative bundle gerbes"
studied in   \cite{CJMSW, CW1}. They  are  also closely related to
 ``simplicial gerbes" of Brylinski \cite{brylinski1}, which
was later further developed by Stevenson under the name ``simplicial
bundle gerbes" \cite{Stevenson}. 
\end{rmk}

\subsection{Compatibility of the maps $T_1$ in \v{C}ech and  de Rham
cohomology}

Let   $N\stackrel{\varphi}{\to} \gm$ be  a crossed module,  and
 $N\toto N_0$
and $\gm \toto\gm_0$ are proper Lie groupoids.
The purpose of this subsection is to compare the two
transgression maps $T_1: \check{H}^{i+1}(\gm\lcom;\sS^1)\to 
\check{H}^{i}((N\rtimes\gm)\lcom;\sS^1)$
and $T_1: H^{i+1}_{dR}(\gm\lcom)\to H^i_{dR}((N\rtimes\gm)\lcom)$.
First we state two general results.
\begin{lem}\label{lem:Tdel}
Let $0\to\iI\to\aA\to\bB\to 0$ be an exact sequence of abelian sheaves
over a  $\Delta_2$-space $M\lcomm$. Then the
diagram
$$\xymatrix{
{H}^i(\M{0};\bB)\ar[r]^{-T_1}\ar[d]^\deltaa
  & {H}^{i-1}(\M{1};\bB)\ar[d]^\deltaa\\
{H}^{i+1}(\M{0};\iI)\ar[r]^{T_1}
  & {H}^i(\M{1};\iI)
}$$
commutes, where $\deltaa$ denotes the boundary maps.
 A similar assertion holds for  \v{C}ech cohomology as well.
\end{lem}
\begin{pf}
Follows immediately from the definition of the boundary maps and 
the relation (\ref{eqn:T1}).
\end{pf}

%\begin{tiny}
%\begin{lem}\label{lem:Tdel}
%Let $0\to\iI\to\aA\to\bB\to 0$ be an exact sequence of abelian sheaves
%over the $\Delta_2$-space $N\lcom\rtimes\gm\lcom$. Then the
%diagram
%$$\xymatrix{
%{H}^i(\gm\lcom;\bB)\ar[r]^{-T_1}\ar[d]^\delta
%  & {H}^{i-1}((N_1\rtimes\gm)\lcom;\bB)\ar[d]^\delta\\
%{H}^{i+1}(\gm\lcom;\iI)\ar[r]^{T_1}
%  & {H}^i((N_1\rtimes\gm)\lcom;\iI)
%}$$
%commutes ($\delta$ denotes the boundary maps). A similar assertion holds
%for  \v{C}ech cohomology.
%\end{lem}
%\end{tiny}

%\begin{tiny}
%\begin{lem}\label{lem:cechsheaf}
%Given a sheaf $\aA$ over a  $\Delta_2$-space $M\lcomm$,
%the  following  diagram
%$$\xymatrix{
%H^i(\gm\lcom;\aA)\ar[r]^{T_{1}}\ar[d]_\cong &
%H^{i-1}((N_1\rtimes\gm)\lcom;\aA)\ar[d]^\tau \\
%\check{H}^i(\gm\lcom;\aA)\ar[r]^{T_{1}} &  
%\check{H}^{i-1}((N_1\rtimes\gm)\lcom;\aA)
%}$$
%commutes.
%\end{lem}
%\end{tiny}

\begin{lem}\label{lem:cechsheaf}
Given an abelian  sheaf $\aA$ over a  $\Delta_2$-space $M\lcomm$,
the  following  diagram
$$\xymatrix{
H^i(\M{0};\aA)\ar[r]^{T_{1}}\ar[d]_\cong & H^{i-1}(\M{1};\aA)\ar[d]^\cong  \\
\check{H}^i(\M{0};\aA)\ar[r]^{T_{1}} &  \check{H}^{i-1}(\M{1};\aA)
}$$
commutes.
\end{lem}
\begin{pf}
For any sheaf $\aA$ on a space $X$, denote by $\hat{\aA}$ the sheaf
such that $\hat{\aA}(U)=\prod_{x\in U}\aA_x$, i.e. a section of
$\hat{\aA}$ over $U$ consists of a (generally discontinuous) map
$x\mapsto f(x)\in \aA_x$.
As in \cite[Lemma~7.2]{tu05}, let $I(\aA)^{k,l} =
\teps_{0*} \hat{\aA}^{k,l+1}$
be the image sheaf \cite{Godement} of $\hat{\aA}^{k,l+1}$ 
 under the smooth map
$\teps_0:M_{k,l+1}\to M_{k,l}$, where $\teps_0$ is the map 
corresponding to $\varepsilon_0:[l]\to [l+1]$, $i\mapsto i+1$, $\forall i$.

Define  a $\Delta_2$-structure
on $I(\aA)\upcom$ as follows. For any  $f=(a,b,c)\in\Hom_{\Delta_2}
((k,l),(k',l'))$, let $f'=(a',b',c')  \in \Hom_{\Delta_2}((k,l+1),(k',l'+1))$
be the unique element satisfying the conditions
  \begin{equation}
\label{eq:ff}
\varepsilon_0\smalcirc f= f'\smalcirc\varepsilon_0,\;
b'(0)=c'(0)=0.
\end{equation}
More precisely, $a'=a$, $b'(i+1)=b(i)$, $c'(0)=0$, and
$c'(i+1)=c(i)+1\;, \forall i\ge 0$. 
For any $U\subset M_{k,l}$, $V\subset M_{k',l'}$ 
such that $\tilde{f}(V)\subset U$, from Eq. \eqref{eq:ff} it follows
 that $\tilde{f'}(\teps_0^{-1}(V))\subset \teps_0^{-1}(U)$.
Set 
$$\tilde{f}*=\tilde{f'}*: I(\aA)^{k', l'}(V) (\cong
\hat{\aA}^{k', l' +1} (\teps_0^{-1}(V)))\lon I(\aA)^{k', l'}(U) (\cong
\hat{\aA}^{k', l' +1} (\teps_0^{-1}(U))) .$$

It is simple to see that this indeed defines a $\Delta_2$-structure
on $I(\aA)\upcom$.
%\comment{Is there way to make the above argument working?
%Otherwise, just delete it and keep the old argument below}
%
%FINDD
%
%As in \cite[Lemma~7.2]{tu05}, we define a sheaf $I(\aA)^{k,l}$
%such that for every open set $U\subset M_{k,l}$,
%$I(\aA)^{k,l}(U)$ consists of sections $\varphi$ (continuous or not)
%over ${\teps_0}^{-1}(U)$
%of the sheaf $\aA^{k,l+1}$,
%such that $\varphi(x)\in
%\aA_x^{k,l+1}$ for all $x\in{\teps_0}^{-1}(U)$.
%(We recall that $\varepsilon_0:[l]\to [l+1]$ is the map $i\mapsto i+1$,
%$\forall i$, and  $\teps_0:M_{k,l+1}\to M_{k,l}$ its induced map.)
%
%Define the $\Delta_2$-structure
%on $I(\aA)\upcom$ as follows. For any  $f=(a,b,c)\in\Hom_{\Delta_2}
%((k,l),(k',l'))$,
%$U\subset M_{k,l}$, $V\subset M_{k',l'}$ such that $\tilde{f}(V)\subset U$
%and for any  $\varphi\in I(\aA)^{k,l}(U)$, we need to produce a section
%$\tilde f(\varphi)$ in $I(\aA)^{k',l'}(V)$.
%Let  $f'=(a',b',c')$ be  the unique element in
%$\Hom_{\Delta_2}((k,l+1),(k',l'+1))$ such that $\varepsilon_0\smalcirc f=
%f'\smalcirc\varepsilon_0$ and $c'(0)=0$.
%More precisely, $a'=a$, $b'=b$, $c'(0)=0$, and \comment{why $b'=b$?}
%$c'(i+1)=c(i)+1\;, \forall i\ge 0$.
%Then the section $\tilde{f}(\varphi)$ is defined to be
%$$x\in{\teps_0}^{-1}(V)\mapsto \varphi ({\tilde{f}}'(x))\in
%\aA^{k,l+1}_{{\tilde{f}}'(x)}\mapsto ({\tilde{f}}')^*\varphi({\tilde{f}}'(x))
%\in \aA_x^{k,l+1}.$$

Consider the diagram
$$\xymatrix{
H^i(M_{0,\com};\aA)\ar[r]^{T_1}&H^{i-1}(M_{1,\com};\aA)\\
H^{i-1}(M_{0,\com};U(\aA))\ar[r]^{T_1} \ar[u]^\deltaa&
 H^{i-2}(M_{1,\com};U(\aA))\ar[u]^\deltaa
}$$
where $U(\aA)=\aA/I(\aA)$ and the boundary maps $\deltaa$ come
from the exact sequence of sheaves
$0\to I(\aA)\to\aA\to U(\aA)\to 0$.

The diagram is anti-commutative, and the left vertical map is an isomorphism
for $i\ge 2$ (see \cite{tu05}). Similar assertions  hold
for \v{C}ech
cohomology as well. Therefore, it suffices  to prove the assertion
for $i=1$. The general conclusion follows from
induction. Let us sketch an  easy and direct argument below.

We have that $H^1(M_{0,\com};\aA)\cong \check{H}^1(M_{0,\com};\aA)$,
which  consists of continuous
sections $\gamma\mapsto \varphi(\gamma)\in \aA^{0,1}_\gamma$
satisfying the equation
 $\tilde{\varepsilon}_1^*\varphi=\tilde{\varepsilon}_0^*\varphi
+\tilde{\varepsilon}_2^*\varphi$. 
Also, $H^0(M_{1,\com};\aA)\cong \check{H}^0(M_{1,\com};\aA)$, which
is the space of $\gm$-invariant sections $\aA^{1,0}(M_{1,0})^{inv}
=\ker(\tilde{\varepsilon_1}^*-\tilde{\varepsilon}_0^*)$.

The maps $T_1: H^1(M_{0,\com};\aA)\to H^0(M_{1,\com};\aA)$ 
and 
$T_1: \check{H}^1(M_{0,\com};\aA)\to \check{H}^0(M_{1,\com};\aA)$
are both equal to $\tilde{\rho}^*$, where $\rho\in\Hom_{\Delta_2}((0,1),(1,0))$
is the morphism $(\emptyset,\mbox{Id},0)$
(note  that  in the case of a crossed module $N\stackrel{\varphi}{\to}\gm$,
$\tilde{\rho}$ coincides with the map $\varphi:N\to\gm$).

Let us explain why $T_1=\tilde{\rho}^*$ for instance in sheaf cohomology.
Represent an element of $H^1(M_{0,\com};\aA)$ by a map $\gamma\mapsto \varphi(\gamma)$
as above. Let $(\aA_\com)$ be an injective resolution of $\aA$. In the complex
$\oplus_{p,q}\aA_q(M_{0,p})$, the cohomology class $[\varphi]$ is represented
by the cochain $\varphi\in \aA(M_{0,1})$, where $\aA(M_{0,1})$ is viewed
as a subspace of $\aA_0(M_{0,1})\oplus \aA_1(M_{0,0})$ using the injection
$\aA\to \aA_0$. By definition of $T_1$, we have $T_1(\varphi)=\tilde{f}^*_\sigma
\varphi\in \aA(M_{1,0})$ where $\sigma\in S_{1,0}$ is the unique
$(1,0)$-shuffle. By definition of $f_\sigma$ (see (\ref{eqn:fsigma})),
we have $f_\sigma=\rho$.
\end{pf}

\begin{prop}
\label{pro:3.8}
Let   $N\stackrel{\varphi}{\to} \gm$ be  a crossed module,  and
 $N\toto N_0$
and $\gm \toto\gm_0$  proper Lie groupoids.
Denote by $\deltaa$ the  boundary maps in the 
induced long exact sequence of  cohomology corresponding to  the short
exact sequence of sheaves
$$0\to \zz\to\rR\to\sS^1\to 0.$$

The following diagram
$$\xymatrix{
\check{H}^i(\gm\lcom;\sS^1)\ar[r]^{-T_1}\ar[d]^\deltaa
  & \check{H}^{i-1}((N\rtimes\gm)\lcom;\sS^1)\ar[d]^\deltaa\\
\check{H}^{i+1}(\gm\lcom;\zz)\ar[r]^{T_1}\ar[d]^{\iota}
  & \check{H}^i((N\rtimes\gm)\lcom;\zz)\ar[d]^\iota\\
\check{H}^{i+1}(\gm\lcom;\rr)\ar[r]^{T_1}\ar[d]^\psi
  & \check{H}^i((N\rtimes\gm)\lcom;\rr)\ar[d]^\psi\\
H^{i+1}_{dR}(\gm\lcom)\ar[r]^{T_1}
  & H^i_{dR}((N\rtimes\gm)\lcom)
}$$
commutes. 
Here the maps $\iota$ are $\check{H}^*(-;\zz)\to
 \check{H}^*(-;\zz)\otimes\rr$, 
and  $\psi$ are the  natural isomorphisms.
Moreover, if $i\ge 2$, the maps $\deltaa$ are isomorphisms.
\end{prop}
\begin{pf}
The fact that the first square commutes follows from Lemma~\ref{lem:Tdel}.
The boundary maps $\deltaa$ are isomorphisms because
$\check{H}^i(\gm\lcom;\rR)\cong {H}^i(\gm\lcom;\rR)=0$,
for all $i\ge 1$,  since  $\gm$ is proper \cite{Crainic}.
Commutativity of the second square is obvious. Commutativity of the
third square follows from  Lemma~\ref{lem:cechsheaf},
since  the  de Rham cohomology is isomorphic to
 the sheaf cohomology of $\rr$.
\end{pf}

Consider the special case that  $N=\gm=G$ is a compact,  connected and
simply connected simple Lie group, with $G$ acting on $N=G$ by conjugation.
It is standard that 
$$\check{H}^3(G\lcom;\sS^1)\cong {H}^3(G\lcom;\sS^1)=\zz $$
and
$$ \check{H}^2(( G\rtimes G)\lcom;\sS^1)\cong {H}^2(( G\rtimes G)\lcom, \zz)
=\zz.$$
As a consequence of Proposition \ref{pro:3.8}, we see that 
the transgression map 
$T_1:\check{H}^3(G\lcom;\sS^1)\to \check{H}^2(( G\rtimes G)\lcom;\sS^1)$
 is indeed an isomorphism.

\begin{cor}\label{cor:T1iso}
Let $G$ be a  compact,  connected and
simply connected simple Lie group. Then
the transgression map
$T_1: \check{H}^3(G\lcom;\sS^1)\to \check{H}^2((G\rtimes G)\lcom;\sS^1)$
is an isomorphism. 
\end{cor}
\begin{pf}
According to Proposition \ref{pro:3.8}, 
it suffices to show that
$T_1: H^4_{dR}(G\lcom)\to H^3_{dR}((G\rtimes G)\lcom)$
maps the canonical generator to the canonical generator.
% This is a direct verification, using the explicit formula for $T_1$:
It is well known that the canonical generator of $H^4_{dR}(G\lcom)$
is given by $\lambda-\Omega\in \Omega^2(G\times G)\oplus \Omega^3(G)$
\cite{BSS, Weinstein},
where $\lambda=(\theta_1, \bar{\theta}_2)$, $\Omega=\frac{1}{6}
(\theta,  [\theta,\theta])$ and $\theta$ (resp. $\bar{\theta}$)
denotes the left (resp. right) Maurer-Cartan form on $G$.
Here $(\cdot, \cdot)$ stands for the Killing form on $\Gg$.
According to Eq.  (\ref{eqn:T1dR}), the image  of $\lambda-\Omega$ 
by $T_1$ is $\omega-\Omega \in \Omega^2(G\times G)\oplus \Omega^3(G)$,
where $\omega= \lambda-f^* \lambda$ and  $f(x,g)=(g,x^g)$. A simple calculation shows that
\begin{equation}
\label{eq:quasi}
\omega|_{(x, g)} =- [(g^* \bar{\theta},  Ad_{x} g^* \bar{\theta})
+(g^* \bar{\theta} , x^{*}(\theta +\bar{\theta} ))],
\end{equation}
where $(x, g)$ denotes the coordinate in $G\times G$, and $x^* \theta$ and
$g^* \theta$ are, respectively,
  the $\Gg$-valued one-forms
on $G\times G$ obtained by pulling back
$\theta$ via the first and second projections, and similarly
for $x^* \bar{\theta} $.
According to \cite{BXZ}, $\omega-\Omega$ is indeed  the canonical generator
 of $H^3_{dR}((G\rtimes G)\lcom)$. \footnote{Note that the transformation
groupoid in \cite{BXZ} uses the left conjugation  action.
The map $\phi: G\rtimes G\to G\ltimes G$ 
 establishing the groupoid isomorphisms from
our convention to that in  \cite{BXZ} is given
by $\phi (x, g)=(g, g^{-1}xg)$.} 
\end{pf}

%\subsection{\v{C}ech cohomology}
%Let us recall a few basic facts on \v{C}ech cohomology of simplicial
%(topological) spaces. Let $M\lcom=(M_n)_{n\in\nn}$ be a simplicial
%space, and let $\aA\upcom$ be a sheaf on $M\lcom$.
%
%A cover $\uU$ of $M\lcom$ is a sequence $\uU^n$ of covers of $M_n$.
%If $\uU^n=(U^n_i)_{i\in I_n}$, one defines a differential complex
%as follows: denote by $P_n^k$ the set of nondecreasing maps $f:
%[k]\to [n]$. Let $\hat{I}_n$ be the set of families $(i_f)_{f\in P_n}$,
%such that $i_f\in I_k$ if $f\in P_n^k$. Let
%$$U_i=\cap_{0\le k \le n}\cap_{f\in P_n^k} \tilde{f}^{-1}(U^k_{i_f}).$$
%Let $C^n(\uU;\aA)$ be the set of families $(c_i)_{i\in \hat{I}_n}$
%such that $c_i\in \aA^n(U_i)$ (i.e. $c_i$ is a section of the sheaf
%$\aA^n$ over the open set $U^n_i\subset M_n$).
%
%The cover $((U_i)_{i\in \hat{I}_n})_{n\in\nn}$ of $M\lcom$ is pre-simplicial:
%for all nondecreasing $g:[m]\to [n]$, one defines $\tilde{g}:\hat{I}_n
%\to\hat{I}_m$ by $(\tilde g (i))_f=i_{g\circ f}$; we have
%$\tilde h\tilde g = \widetilde{gh}$, and $\tilde{g}(U_i)
%\subset U_{\tilde{g}(i)}$.
%
%Then the formula
%$$(\del c)_f=\sum_{i=0}^{n+1} (-1)^i\teps_i^* c_{\teps_i(f)}$$
%($c\in C^n(\uU;\aA\upcom)$) defines a differential, and the
%cohomology groups are denoted by $H^n(\uU;\aA\upcom)$.
%The inductive limit of the groups $H^n(\uU;\aA\upcom)$ when $\uU$ runs
%over finer and finer covers is denoted by $\check{H}^n(M\lcom;\aA\upcom)$
%and is called the \v{C}ech cohomology of $M\lcom$ with coefficients
%in $\aA\upcom$. It is isomorphic to $H^n(M\lcom;\aA\upcom)$ if
%$M_n$ is paracompact for all $n$.

\begin{rmk}
In \cite{BM1},  Brylinski and McLaughlin
 studied the  transgression map  from $H^4 (BG)$
 to $H^3 (BLG)$,  and the relation with Segal-Witten reciprocity law.
It would be interesting to investigate the connection
with our transgression map here.
\end{rmk}

\section{Ring structure on the equivariant  twisted $K$-theory}

In this section, we assume that 

\begin{itemize}
\item[(i)] $N\stackrel{\varphi}{\to} \gm$ is
a crossed module, where $N\toto N_0$ and $\gm\toto \gm_0$
 are proper  Lie groupoids such that the source (or  target) map $N\to N_0$
is $\gm$-equivariantly  $K$-oriented (see definition below); and
\item[(ii)] $(c, b, a)$ is a multiplicator,
where $c\in \check{C}^2((N\rtimes \gm)\lcom, \sS^1)$,
$b\in \check{C}^1((N_2\rtimes \gm)\lcom, \sS^1)$ and
$a\in \check{C}^0((N_3\rtimes \gm)\lcom, \sS^1)$.
That is,
$$\del c=0, \ \ \del'c = \del b, \ \ \mbox{ and }  \del' b = \del a.$$
\end{itemize}

The main 
purpose of this section is to construct a ring structure
on the twisted  equivariant
$K$-theory groups $K^*_{c,\gm}(N)$ for   any 
given multiplicator $(c, b, a)$,  and investigate
 how the ring
structure depends on the choice of these multiplicators.

\subsection{2-cocycles and $S^1$-central extensions}
\label{subsec:2cocycle}

Let us briefly recall how an
 $\sS^1$-valued 2-cocycle determines
an $S^1$-central extension. Given a groupoid $\gm\toto\gm_0$, a cover
$\uU=(\uU^k)_{k\in\nn}$ of $\gm\lcom$ with $\uU^k=(U^k_j)_{j\in I_k}$,
and a 2-cocycle $c\in\check{Z}^2(\sigma\uU;\sS^1)$,
$c$ can be viewed as a family $(c_i)$, with
$i=(i_0,i_1,i_2,i_{01},i_{02},i_{12},i_{012})\in I_0^3\times I_1^3\times I_2$.
Here $c_i$ is an $S^1$-valued function on some open set of $\gm_2$.
The $2$-cocycle $c$ determines an $S^1$-central extension
\begin{equation}
\label{eq:c}
S^1\to \tilde{\gm}^c\stackrel{\pi}{\to}
\gm^c\toto \amalg_{i \in I_0} U^0_i,
\end{equation}
where $\gm^c\toto \amalg_{i \in I_0} U^0_i$ is the pull-back
 groupoid
$\gm[\uU^0]\toto \amalg_{i \in I_0} U^0_i$
given by $\gm[\uU^0]=\{(i,j,\gamma)|\;i,j\in I_0\mbox{ and }\gamma\in
{\gm}^{U^0_i}_{U^0_j}\}$. The groupoid $\tilde{\gm}^c
\toto \amalg_{i \in I_0} U^0_i$ is

\begin{equation}\label{eqn:tgmc}
\tilde{\gm}^c=
\{(k,i,j,\gamma,\lambda)\vert\; i,j\in I_0, k\in I_1,
\gamma\in\gm_{U^0_j}^{U^0_i}\cap U^1_k,\lambda\in S^1\}/\sim,
\end{equation}
where $\sim$ is an appropriate equivalence relation,
$\pi(k,i,j,\gamma,\lambda)=(i,j,\gamma)$, and the product on
$\tilde{\gm}^c$ is defined by
\begin{equation}
\label{eqn:prod-ext}
(i_{01},i_0,i_1,\gamma_1,\lambda)(i_{12},i_1,i_2,\gamma_2,\mu)
=(i_{02},i_0,i_2,\gamma_1\gamma_2,\lambda\mu c_i(\gamma_1,\gamma_2)).
\end{equation}
We refer the reader to \cite{tu05} for details.

Recall that $\check{C}^n(\gm\lcom;\sS^1)
=\lim_{\uU}\check{C}^n(\sigma\uU;\sS^1)$,
where $\uU$ runs over all covers of the form
$(U^k_\gamma)_{\gamma\in\gm_k}$, $\gamma\in U^k_\gamma$,
and $\uU'$ is finer than $\uU$ if ${U'}^k_\gamma \subset U^k_\gamma$
for all $k$ and all $\gamma\in\gm_k$.

\begin{lem}\label{lem:cocycle-iso-gpoid}
For any  open cover $\uU$ of $\gm\lcom$,
if $c^1$, $c^2\in \check{Z}^2(\sigma\uU;\sS^1)$, and  $b\in \check{C}^1(\sigma\uU;\sS^1)$
satisfy $c^1-c^2=\del b$, then there is
an isomorphism of $S^1$-central extensions
$$\Phi^\uU_{c^2,b,c^1}:  (\tilde{\gm}^{c_1}\stackrel{\pi}{\to}
\gm^{c_1}) \lon (\tilde{\gm}^{c_2}\stackrel{\pi}{\to}
\gm^{c_2})$$
satisfying the following properties:
\begin{itemize}
\item[(a)] $\Phi^\uU_{c^3,b',c^2}\smalcirc \Phi^\uU_{c^2,b,c^1}=
\Phi^\uU_{c^3,b +b',c^1}$, where $c^1-c^2=\del b$ and $c^2-c^3=\del b'$;
\item[(b)] if  $b=\del a$ for some
$a\in \check{C}^0(\sigma\uU;\sS^1)$, then
 $\Phi^\uU_{c,b,c}: \tilde{\gm}^{c}\to  \tilde{\gm}^{c}$
is the multiplier by $a(t(g))a(s(g))^{-1}$, $\forall g \in \tilde{\gm}^{c}$;
and
\item[(c)] if $\uU'$ is a refinement of $\uU$ and ${c'}^1$, ${c'}^2$ and
$b'$ are the images of $c^1$, $c^2$ and $b$ under  the canonical map
$\check{C}^*(\sigma\uU;\sS^1)\to \check{C}^*(\sigma\uU';\sS^1)$, then
the following diagram
$$\xymatrix{
{\tilde{\gm}}^{{c'}^1}\ar[r]\ar[d]_{\Phi^{\uU'}_{{c'}^2,b',{c'}^1}} & 
{\tilde{\gm}}^{c^1}\ar[d]^{\Phi^\uU_{c^2,b,c^1}} \\
{\tilde{\gm}}^{{c'}^2}\ar[r] &  {\tilde{\gm}}^{c^2}
}$$
commutes, where the horizontal maps are Morita morphisms of
$S^1$-central extensions.
\end{itemize}
\end{lem}
\begin{pf}
In terms of  notation (\ref{eqn:tgmc}), 
we  denote by $\varphi_k^1$ and $\varphi_k^2$ the maps
 $$\varphi_k^1 (i,j,\gamma)
=[(k,i,j,\gamma,1)]\in \tilde{\gm}^{c^1}, \ \ 
\varphi_k^2 (i,j,\gamma)
=[(k,i,j,\gamma,1)]\in \tilde{\gm}^{c^2}.$$
 Then $\varphi_k^1$
and $\varphi_k^2$ are local sections
 of the $S^1$-principal bundles
 $\tilde{\gm}^{c^1}\stackrel{\pi}{\to}{\gm^{c^1}}$
and $\tilde{\gm}^{c^2}\stackrel{\pi}{\to}{\gm^{c^2}}$, respectively.
The relation $c^2-c^1=\del b$ means that
\begin{eqnarray}
\label{eqn:c1c2b}
c^2_i(\gamma_1,\gamma_2) c^1_i(\gamma_1,\gamma_2)^{-1}
&=&b_{i_0i_1i_{01}}(\gamma_1)b_{i_1i_2i_{12}}(\gamma_2)
b_{i_0i_2i_{02}}(\gamma_1\gamma_2)^{-1},\\
&&\forall 
i=(i_0,i_1,i_2,i_{01},i_{02},i_{12},i_{012})\in I_0^3\times I_1^3\times I_2 .
\end{eqnarray}

By applying Eq.  (\ref{eqn:prod-ext}) to the $2$-cocycles
 $c^1$ and $c^2$, we are led to
\begin{equation}
\label{eqn:prod-ext2}
\varphi^j_{i_{01}}(i_0,i_1,\gamma_1)
\varphi^j_{i_{12}}(i_1,i_2,\gamma_2)
=\varphi^j_{i_{02}}(i_0,i_2,\gamma_1 \gamma_2)c^j_i(\gamma_1,\gamma_2), \ \  j=1, \ 2.
\end{equation}

Define a map $\Phi^\uU_{c^2,b,c^1}:
\tilde{\gm}^{c^1}\to \tilde{\gm}^{c^2}$ by
\begin{equation}
\label{eqn:defphi1}
[(k,i,j,\gamma,\lambda)]  \lon [(k,i,j,\gamma,\lambda b_{ijk}(\gamma)^{-1})].
\end{equation}
That is,
\begin{equation}
\label{eqn:defphi}
\lambda\varphi^1_k(i,j,\gamma)\mapsto \lambda\varphi^2_k (i,j,\gamma)
b_{ijk}(\gamma)^{-1}
\end{equation}
for all $\lambda\in S^1$ and
$\gamma\in U^1_k\cap\gm^{U^0_i}_{U^0_j}$.

Let us check that $\Phi^\uU_{c^2,b,c^1}$ is well-defined, i.e. 
it is independent of the choice of $k$. 
To simplify notations, after replacing $\gm$ by $\gm[\uU_0]$
and the cocycles $c^j$ by ${c'}^j_{i_{01}i_{12}i_{02}}((i_0,i_1,\gamma_1),
(i_1,i_2,\gamma_2))=c^j_{i_0i_1i_2i_{01}i_{12}i_{02}}(\gamma_1,\gamma_2)$,
we can assume
 that the cover $\uU^0$ consists of the entire open set $\gm_0$.
Thus Eq.  (\ref{eqn:defphi}) reads

\begin{equation}
\label{eqn:defphi2}
\lambda\varphi^1_k(\gamma)\mapsto \lambda\varphi_k^2(\gamma)
b_k(\gamma)^{-1} .
\end{equation}

Using Eq. (\ref{eqn:prod-ext2}), we get $\varphi^j_m(s(\gamma))
\varphi^j_m(s(\gamma))=\varphi^j_m(s(\gamma))
c^j_{mmm}(s(\gamma),s(\gamma))$. Thus one can identify
$\varphi^j_m(s(\gamma))$ with  the complex number
$c^j_{mmm}(s(\gamma),s(\gamma))$.
By  Eq. (\ref{eqn:prod-ext2}) again, we obtain  $\varphi^j_l(\gamma)
\varphi^j_m(s(\gamma))=\varphi^j_k(\gamma)c^j_{lmk}(\gamma,s(\gamma))$.
For $j=1$, after replacing $\varphi^j_m(s(\gamma))$ by 
$c^j_{mmm}(s(\gamma),s(\gamma))$, we have 
\begin{equation}
\label{eqn:phi1}
\varphi^1_l(\gamma)=\lambda_1\varphi^1_k(\gamma)
\end{equation}
with $\lambda_j= c^j_{lmk}(\gamma,s(\gamma))/c^j_{mmm}(s(\gamma),s(\gamma))$
($j=1,2$).
>From Eq. (\ref{eqn:c1c2b}), we obtain that
$\lambda_2=\lambda_1b_l(\gamma)b_k(\gamma)^{-1}$. Therefore

\begin{equation}
\label{eqn:phi2}
\varphi^2_l(\gamma)b_l(\gamma)^{-1}=\lambda_1\varphi^2_k(\gamma)
b_k(\gamma)^{-1}.
\end{equation}

Comparing Eqs. (\ref{eqn:defphi2})-(\ref{eqn:phi2}),
we see that $\Phi^\uU_{c^2,b,c^1}$ is indeed well-defined.

The fact that $\Phi^\uU_{c^2,b,c^1}$ is a groupoid morphism now immediately
follows from Eqs. (\ref{eqn:c1c2b})-(\ref{eqn:prod-ext2}).
Relation (a) is also clear from Eq. (\ref{eqn:defphi}).
\par\medskip

Let us prove (b). Again, for simplicity,
 we may  assume that $\uU^0$ is a cover consisting of
just one  open set $\gm_0$. Hence,   by assumption,
 $b_i(g)=a(t(g))a(s(g))^{-1}$
for all $i\in I_1$ and $g\in U^1_i$. Hence the automorphism of
the central extension $S^1\to {\tilde{\gm}}^c\to \gm^c$ given by
Eq. \eqref{eqn:defphi1} is 
$$\Phi^\uU_{c,b,c}(\gamma)=\gamma
a(s(\gamma))a(t(\gamma))^{-1}\ ,\quad  \forall \gamma
\in {\tilde{\gm}}^c $$
\par\medskip
Let us show (c). By definition, a refinement is a map
$\psi:I'_n\to I_n$ ($n\in\nn$) such that ${U'}^n_i\subset
U^n_{\psi(i)}$.
Consider the following diagram
$$\xymatrix{
{\tilde{\gm}}^{{c'}^1}\ar[r]\ar[d]_{\Phi^{\uU'}_{{c'}^2,b',{c'}^1}} & 
{\tilde{\gm}}^{c^1}\ar[d]^{\Phi^\uU_{c^2,b,c^1}} \\
{\tilde{\gm}}^{{c'}^2}\ar[r] &  {\tilde{\gm}}^{c^2}
}$$
Here the horizontal maps are defined by
$[k,i,j,\gamma,\lambda]\mapsto [\psi(k),\psi(i),\psi(j),\gamma,\lambda]$,
which are clearly Morita morphisms of $S^1$-central extensions.
It is immediate to check that the above diagram indeed commutes.
\end{pf}

\subsection{The $C^*$-algebra associated to  a 2-cocycle}
\label{sec:4.2}
Let us recall the construction of the $C^*$-algebra associated to an 
$S^1$-central extension
$$S^1\to \tilde{\gm}\to \gm \toto \gm_0,$$
where $\gm\toto \gm_0$ is a Lie groupoid, or, more generally, a
locally compact groupoid endowed with a Haar
system $(\lambda^x)_{x\in \gm_0}$.
Let $L=\tgm\times_{S^1}\cc$  be the complex line 
bundle associated to the $S^1$-principal bundle $\tgm\to \gm$.
Then $L\to \gm\toto M$
is equipped with   an associative bilinear product
\begin{eqnarray*}
L_g\otimes L_h&\to& L_{gh} \ \ \ \forall (g,h)\in \gm_2\\
(\xi,\eta)&\mapsto&\xi\cdot\eta
\end{eqnarray*}
and an antilinear  involution
\begin{eqnarray*}
L_g&\to& L_{g^{-1}}\\
\xi&\mapsto&\xi^*
\end{eqnarray*}
satisfying the following properties:
\begin{itemize}
\item the restriction of the line bundle to the unit space $M$
is isomorphic to the trivial line bundle $M\times\cc\to M$;
\item $\forall \xi,\eta\in L_g$, $\langle \xi,\eta\rangle =
\xi^*\cdot\eta\in L_{s(g)}\cong \cc$ defines a scalar product;
\item $(\xi\cdot\eta)^* = \eta^*\cdot \xi^*$.
\end{itemize}

 The space $C_c(\gm,L)$ of
continuous, compactly supported sections of the line bundle $L\to\gm$
is endowed with a  convolution product:
$$(\xi *\eta )(g)=\int_{h\in\gm^{t(g)}} \xi (h)\cdot \eta (h^{-1}g)
d\lambda^{t(g)}(h),$$
and  the adjoint
$$\xi^*(g)=(\xi(g^{-1}))^*,$$
where $\xi (h)\cdot \eta (h^{-1}g)$ is understood
as  the product $L_h\otimes L_{h^{-1}g}\to L_g$.

For all $x\in \gm_0$, let $\hH_x$ be the Hilbert space obtained
by completing $C_c(\gm,L)$ with respect to the scalar
product
$$\langle \xi,\eta\rangle = (\xi^* * \eta)(x)
=\int_{g\in\gm^x}\langle \xi(g^{-1}),\eta(g^{-1})\rangle d\lambda^x(g).$$
Let $(\pi_x(f))\xi
=f*\xi$, $\forall f\in C_c (\gm, L)$, $\xi\in \hH_x$.
Then $f\mapsto \pi_x(f)$ is a $*$-representation of $C_c (\gm,L)$
in $\hH_x$. The $C^*$-algebra $C^*_r(\gm,L)$ is, by definition, the
completion of $C_c(\gm,L)$ with respect to the norm:
$\sup_{x\in \gm_0}\|\pi_x(f)\|$.

%\subsection{The $C^*$-algebra associated to a 2-cocycle}
%\label{subsec:alg-2-cocycle}

Let $\gm\toto \gm_0$ be a locally compact groupoid with Haar system.
Let $c\in \check{Z}^2(\sigma\uU;\sS^1)$ be a \v{C}ech 2-cocycle on an
open cover $\uU$ of $\gm$. Then $c$ determines
 an $S^1$-central
extension  as given by Eq. \eqref{eq:c},
where $\gm^c$
is Morita equivariant to $\gm$.
Denote by $L^c=\tgm^c\times_{S^1}\cc$ the associated
complex line bundle.
By $C^*_r (\gm, c)$, we denote the $C^*$-algebra $C^*_r(\gm^c,L^c)$ 
and call the $C^*$-algebra associated to the 2-cocycle
$c$.

\subsection{$S^1$-equivariant gerbes}

The main result of this subsection is  the following

\begin{them}
\label{pro:equ}
Assume  that $\gm\toto \gm_0$ is a Lie  groupoid acting on a manifold $N$
via $J:N\to \gm_0$.
Let $\uU$ be a cover of $(N\rtimes\gm)\lcom$.
Any \v{C}ech $2$-cocycle $c\in \check{Z}^2(\uU, \sS^1 )$
canonically determines an $S^1$-central extension of
the form $\tilde{H}\rtimes \gm\to H\rtimes \gm\toto M$, where
$\tilde{H}\to H\toto M$ is a $\gm$-equivariant  $S^1$-central extension
and $H\toto M$ is Morita equivalent to $N\toto N$, such that the class
of this central extension is equal to
$[c]\in\check{H}^2((N\rtimes\gm)\lcom;\sS^1)$.
\end{them}

We first prove the following special case.

\begin{prop}
\label{pro:eq}
Let $\gm\toto \gm_0$ be   a  Lie 
groupoid, and  $\uU$ a cover of $\gm\lcom$.
Any \v{C}ech $2$-cocycle $c\in \check{Z}^2(\uU, \sS^1 )$
canonically determines an $S^1$-central extension of
the form $\tilde{H}\rtimes \gm\to H\rtimes \gm\toto M$, where
$\tilde{H}\to H\toto M$ is a $\gm$-equivariant  $S^1$-central extension
and $H\toto M$ is Morita equivalent to $\gm_0\toto\gm_0$, such 
that the class of this central extension is equal to
$[c]\in\check{H}^2(\gm\lcom;\sS^1)$.
\end{prop}

\begin{lem}
\label{lem:pull}
Assume that $M$ is a  right $\gm$-space
via the momentum map  $J: M \to \gm_0$, which is  a surjective
submersion. Then the pull-back groupoid $\gm [M]\toto M$ via $J$ 
is naturally isomorphic to the crossed-product
 $H\rtimes \gm \toto M$, where $H\toto M$ is 
the equivalence groupoid $M\times_{\gm_0}M\toto M$,
and $\gm$ acts on $H\toto M$ by natural automorphisms.
\end{lem}
\begin{pf}
By definition, $\gm [M]=\{(x, y, g)|J(x)=t(g), J(y)=s(g)\}$. 
One checks directly that $\chi:  H\rtimes \gm \to \gm [M]$
given by $(x, y, g)\mapsto (x, yg^{-1}, g)$ is
indeed a groupoid isomorphism.
\end{pf}
{\bf Proof of Proposition \ref{pro:eq}}  
\begin{pf}
From Section~\ref{subsec:2cocycle}, we know that $c$  canonically
  determines
\begin{itemize}
\item[(a)] a surjective submersion $\pi: U\to \gm_0$; and
\item[(b)] an $S^1$-central extension
\begin{equation}
\label{eqn:central-ext}
S^1\to \tilde{\gm}\stackrel{p}{\to}\gm[U] \toto  U,
\end{equation}
\end{itemize}
where $U=\amalg_{i\in I_0} U^0_i$.  Let
 $M=U\times_{\gm_0}\gm$.
It is clear that $\gm$ acts on  $M$ by $(x,g)\cdot\gamma=(x,g\gamma)$
with the momentum map $\sigma :M\to \gm_0$ given by
$\sigma(x,g)=s(g)$, $\forall (x,g)\in U\times_{\gm_0}\gm$.
Thus according to Lemma  \ref{lem:pull}, the pull
back groupoid $\sigma^*\gm$ is isomorphic
to the crossed-product $H\rtimes \gm\toto M$,
where $H=M\times_{\sigma, \gm_0, \sigma}M=\{(x',y')\in M\times M\vert
\sigma(x')=\sigma(y')\}$ and $\gm$ acts on $H$ by
natural automorphisms.  Now consider the map
$\tau:M\to \gm_0$ given by $\tau(x,g)=t(g)$, $\forall 
(x,g)\in U\times_{\gm_0}\gm$. It is clear that $\tau$ is also a surjective submersion.
Introduce a map
\begin{eqnarray}
\rho: \sigma^*\gm &\lon &\tau^* \gm\\
\big((x, g), (y, h), r)\big)&\lon &\big((x, g), (y, h), grh^{-1} \big)  .
\label{eq:rho}
\end{eqnarray}
It is simple to check that $\rho$ is a groupoid isomorphism.
It thus follows that   the groupoids $H\rtimes \gm \toto M$ 
and $\tau^* \gm\toto M$ are isomorphic. The isomorphism
$\psi: H\rtimes \gm\to \tau^* \gm$ 
is given by  $\psi((x,g),(y,h),\gamma)=((x,g),(y,h\gamma),gh^{-1})$,
which is the composition of the isomorphism $\chi$
in the proof of Lemma \ref{lem:pull} with
 the isomorphism $\rho$ as in Eq. \eqref{eq:rho}.
On the other hand, since 
$\tau=\pi\smalcirc p_1$,  where $p_1:M\to U$
 is the natural projection,
then $\tau^*\gm  \stackrel{p_1}{\to}  \pi^*\gm(=\gm[U])$
is a Morita morphism of groupoids. Here, by abuse of notations,
we use the same symbol $p_1$ to denote the groupoid
morphism on the pull-back groupoids induced by
$p_1: M\to U$. Therefore, it follows that 
$$\xymatrix{
H\rtimes \gm\ar[r]^\varphi\ar[d]\ar@<-.6ex>[d] & \gm [U]\ar@<-.6ex>[d]\ar[d] \\
M\ar[r]^{p_1} &  U
}$$
is a Morita morphism of groupoids, where  $\varphi= p_1\smalcirc \psi$
can be expressed  explicitly by $\varphi((x,g),(y,h),\gamma)=(x,y,gh^{-1})$.
Let 
$$\tilde{\gm}'= \varphi^* \tilde{\gm}=
(H\rtimes \gm) \times_{\gm [U]} \tilde{\gm}.$$
Then it is clear that
\begin{equation}\label{eqn:ext-hgm}
S^1\to \tilde{\gm}'\to H\rtimes\gm \toto M
\end{equation}
is an $S^1$-central extension Morita equivalent to 
(\ref{eqn:central-ext}), and therefore it defines the
same cohomology class  $[c]$. 

Identify $H$ with the subgroupoid  $H'$ of $H\rtimes \gm$
via the embedding $i:h\mapsto (h, J_1 (h))$, where
$J_1: H\to \gm_0$ is the momentum map of the
$\gm$-action. Note that $\gm$ naturally acts
on  $H'$ by
$(h, J_1 (h))^r=(h^r, s(r))$, which coincides
with the natural $\gm$-action on $H$ under  the
 identification $H'\iso H$. Let $j=\varphi\circ i:H\to\gm [U]$.
Now let  
$$\tilde{H}= j^*  \tilde{\gm} = H'\times_{\gm [U]} \tilde{\gm}$$
be the pull back $S^1$-bundle, via $j$,
 of  $\tilde{\gm}\stackrel{p}{\to}\gm[U]$
as given  in Eq. \eqref{eqn:central-ext}.

The $\gm$-action on $H'$ induces a $\gm$-action on $\tilde{H}$.
Indeed, $\tilde{H}$ can be  naturally identified with
$\{(h,\tilde{\gamma})\in H\times\tilde{\gm}\vert\;
j(h)=p(\tilde{\gamma})\}$, under which the $\gm$-action is
given  by $(h, \tilde{\gamma})^\gamma=(h^\gamma,\tilde\gamma)$.
This action is well-defined since $j$ is $\gm$-invariant.
It is clear from the construction that
\begin{equation}
\label{eqn:central-ext2}
S^1\to \tilde{H}\stackrel{q}{\to} H \toto M
\end{equation}
is a $\gm$-equivariant  $S^1$-central extension.
 Therefore one obtains
an $S^1$-central extension
\begin{equation}
\label{eq:hgm}
S^1\to \tilde{H}\rtimes\gm \stackrel{q}{\to}
H\rtimes\gm \toto M,
\end{equation}
which is clearly isomorphic to  (\ref{eqn:ext-hgm}). Indeed
the  map $f: \tilde{H}\rtimes\gm \to \tilde{\gm}$, 
 defined by $f(h,\tilde{\gamma},\gamma) =\tilde{\gamma}$,
is an $S^1$-equivariant groupoid morphism and
 the diagram
$$\xymatrix{
{\tilde{H}\rtimes\gm} \ar[r]^f \ar[d]^q
  & {\tilde{\gm}}\ar[d]^p\\
H\rtimes\gm\ar[r]^\varphi & \gm[U]
}$$
commutes.
\end{pf}

Theorem \ref{pro:equ} thus follows from
Proposition \ref{pro:eq} combining with the following

\begin{lem}
Assume  that $ \gm \toto \gm_0$ is a Lie  groupoid acting on a manifold $N$
via $J:N\to \gm_0$. If a groupoid $H\toto M$ admits
an action of $N\rtimes \gm$ by automorphisms,
then it must admit a $\gm$-action by automorphisms.
Moreover, the crossed product groupoids
$H\rtimes (N\rtimes \gm)$ and $H\rtimes \gm$
are canonically isomorphic.
\end{lem}
\begin{pf}
This follows from a direct verification, which we
will leave to the reader.
\end{pf}

\begin{cor}
\label{cor:Ac}
Under the same hypothesis as in Theorem
\ref{pro:equ}, $c$ canonically determines a $\gm$-$C^*$-algebra
$A_c$ and a Morita equivalence $\nu_c$ from the $C^*$-algebra
$C^*_r(N\rtimes\gm,c)$ of the  $S^1$-central extension
(\ref{eqn:central-ext}) to the crossed product $A_c\rtimes_r\gm$.

Moreover, under the same assumptions as in
Lemma~\ref{lem:cocycle-iso-gpoid}, there is a $\gm$-equivariant
Morita equivalence $\alpha_{c^2,b,c^1}$ from $A_{c^1}$ to $A_{c^2}$
such that the diagram
$$\xymatrix{
C^*_r(\gm,c^1)\ar[r]^{\nu_{c^1}}\ar[d]^{\Phi^\uU_{c^2,b,c^1}}&
A_{c^1}\rtimes_r\gm\ar[d]^{\alpha_{c^2,b,c^1}\rtimes\gm}\\
C^*_r(\gm,c^2)\ar[r]^{\nu_{c^2}}&
A_{c^2}\rtimes_r\gm
}$$
commutes. Moreover, the above diagram is compatible with refinements
of covers as in Lemma~\ref{lem:cocycle-iso-gpoid}.
\end{cor}

\begin{pf}
After replacing $\gm$ by $N\rtimes \gm$, we may assume that
 $N=\gm_0$.
Recall from the proof of Proposition~\ref{pro:eq} that the cover
$\uU$ on which $c$ is defined determines an \'etale map
$\pi_c:U^c\to \gm_0$, a groupoid $H^c$, an action of $\gm$ on
$H^c$ by automorphisms,  and  the map $j_c$ which is the composition
$H^c\stackrel{i_c}{\to}H^c\rtimes \gm
\stackrel{\varphi_c}{\to} \gm[U^c]$.

The cocycle $c$ determines an $S^1$-central extension
$S^1\to \tilde{\gm}^c\stackrel{p_c}{\to}\gm[U^c]$,
 whose $C^*$-algebra
is $C^*_r(\gm,c)$. The pull back of this $S^1$-bundle
$S^1\to \tilde{H^c}\to H^c$ by the map
$j_c$ is a $\gm$-equivariant extension as shown
in Eq. \eqref{eqn:central-ext2}.
 We take $A_c$ to be its $C^*$-algebra,
which is thus endowed with an action of $\gm$.

Since there is an $S^1$-equivariant isomorphism $f_c: 
\tilde{H}^c\rtimes \gm \to \varphi_c^*\tilde{\gm}^c$ (see the proof
of Proposition~\ref{pro:eq}), we obtain a canonical
 Morita equivalence $\nu_c$.
Since $j_c$  depends only  on $\uU$, the isomorphism
$\Phi^\uU_{c^2,b,c^1}:(S^1\to\tilde{\gm}^{c^1}\to \gm[U^{c^1}])
\to (S^1\to\tilde{\gm}^{c^2}\to \gm[U^{c^2}])$ determines a $\gm$-equivariant
isomorphism $\alpha_{c^2,b,c^1}: (S^1\to \tilde{H}^{c^1}\to H^{c^1})
\to (S^1\to\tilde{H}^{c^2}\to H^{c^2})$, and therefore  an isomorphism
$\alpha_{c^2,b,c^1}\rtimes \gm : (S^1\to \tilde{H}^{c^1}\rtimes \gm
\to H^{c^1}\rtimes\gm )\to (S^1\to\tilde{H}^{c^2}\rtimes\gm
\to H^{c^2}\rtimes\gm)$. It is easy to check that this isomorphism
coincides with  $\varphi_c^*\Phi^\uU_{c^2,b,c^1}$, via the
identification $f_c:\tilde{H^c}\rtimes\gm\to
\varphi_c^*\tilde{\gm}^c$. Hence the commutative diagram follows.
 The last
assertion is easy to check.
\end{pf}

\begin{rmk}
Let $\gm$ be a Lie groupoid. Recall \cite{KMRW98,TXL04} that the Brauer group
$Br(\gm)$ of $\gm$ is defined as the group of Morita 
equivalence classes of $\gm$-equivariant locally trivial bundles of $C^*$-algebras over $\gm_0$
whose fibers are the algebra of  compact operators on a  Hilbert space.
It is known that $Br(\gm)$ is isomorphic to $H^2(\gm;\sS^1)$ \cite{tu05}.
The  corollary above enables us
to  realize such a bundle geometrically
for a  given  \v{C}ech 2-cocycle $c\in \check{Z}^2(\gm;\sS^1)$.
\end{rmk}

We will need the following functorial property of the  correspondence
$c\mapsto A_c$:

\begin{prop}
\label{prop:fAc}
Let $\gm$ be a  Lie groupoid acting on manifolds $N'$ and $N$,
and $f:N'\to N$  a $\gm$-equivariant map. Assume  that $c\in
\check{Z}^2(\uU,S^1)$, where $\uU$ is a cover of $(N\rtimes\gm)\lcom$.
Then $A_{f^*c}$ is canonically isomorphic to
$f^*A_c \cong C_0 (N') \otimes_{C_0 (N)} A_c$.
Moreover, this isomorphism is compatible with refinements of covers,
and with $\Phi_b$ and $\alpha_b$ as in Corollary~\ref{cor:Ac}.
\end{prop}

(Note  that $A_c$  can be considered as a field of $C^*$-algebras over
$N$, and hence  $f^*A_c$ is a field of $C^*$-algebras over $N'$.)

\begin{pf}
We use the same notations as in the proof of
Proposition \ref{pro:eq}. The groupoid $H$ is a  fibration over $N$,
with each fiber $H_x=\{((i,\gamma),(j,\gamma'))\in
I_0\times\gm^x\times I_0\times \gm^x\vert\; s(\gamma)\in U^0_i,
s(\gamma')=U^0_j\}$.
Similarly, the groupoid $H'$ is a fibration over $N'$, and one immediately
sees that the fiber $H'_{x}$ over $x\in N'$ is isomorphic to
$H_{f(x)}$. More precisely, there is a cartesian diagram
$$\xymatrix{
H'\ar[r] \ar[d] & H \ar[d] \\
N'\ar[r]^f & N.
}$$

Recall also from the proof of Proposition \ref{pro:eq}
 that the central extension $S^1\to \tilde{H}
\to H$ is obtained by the  pull back of the
 central extension of $\tilde{\gm}\to \gm[U]$ 
by the groupoid morphism $\varphi: H\to \gm [U]$.
The central extension $S^1\to \tilde{H}'\to H'$ admits  a similar description.
It is simple to see that the diagram
$$\xymatrix{
H'\ar[r]\ar[d]^{\varphi'} & H\ar[d]^\varphi\\
f^*\gm[U]\ar[r]& \gm[U]
}$$
commutes.  It thus follows that there exists an $S^1$-equivariant map
 ${\tilde{H}}'\to \tilde{H}$ such that
$$\xymatrix{
{\tilde{H}'}\ar[r]\ar[d] &{\tilde{H}}\ar[d]\\
N'\ar[r]^f & N
}$$
commutes.
In other words, considering $\tilde{H}\to N$ as an $S^1$-central extension
of groupoids fibered over $N$, the pull-back by $f$ of
$(\tilde{H}\to N)$ is isomorphic to $(\tilde{H}'\to N')$.
The conclusion thus follows.

The last assertion is easy to check.
\end{pf}

\begin{lem}\label{prop:cc'}
Let $c$, $c'\in \check{Z}^2(\sigma\uU;\sS^1)$, where $\uU$ is a cover of
$(N\rtimes \gm)\lcom$. Then $A_c\otimes_{C_0(\gm_0)} A_{c'}\iso
 A_{c+c'}$ canonically.
\end{lem}

\begin{pf}
If $\eE$ (resp. $\eE'$) is the $S^1$-central extension whose $C^*$-algebra
is $A_c$ (resp. $A_{c'}$), then $A_{c+c'}$ is the $C^*$-algebra
of $\eE+\eE'$.
 The latter is canonically isomorphic to $A_c\otimes_{C_0(\gm_0)} A_{c'}$.
\end{pf}

\begin{prop}\label{prop:p1cp2c}
Let $M$ and  $N$ be  $\gm$-spaces, $c$ and  $c'$
$\sS^1$-valued \v{C}ech 2-cocycles on $M\rtimes\gm$  and
$N\rtimes\gm$, respectively. Let  $P=M\times_{\gm_0}N$ and
$p_1:P\to M$, $p_2:P\to N$ be natural projections. Then
$A_{p_1^*c+p_2^*c}\iso  A_c\otimes_{C_0(\gm_0)} A_{c'}$
canonically.  Moreover this isomorphism is compatible with
 the maps $\Phi_b$.
\end{prop}
\begin{pf}
We have
\begin{eqnarray*}
A_{p_1^*c+p_2^*c'}&\cong& A_{p_1^*c}\otimes_{C_0(P)}
A_{p_2^*c'}\quad\mbox{(by Lemma~\ref{prop:cc'} for $P\rtimes\gm$)}\\
&\cong& p_1^*A_c\otimes_{C_0(P)}p_2^*A_{c'}\quad
\mbox{(by Proposition~\ref{prop:fAc})}\\
&\cong& A_c\otimes_{C_0(\gm_0)} A_{c'}.
\end{eqnarray*}
\end{pf}

\subsection{Equivariant twisted $K$-theory and the map $\Phi_b$}

First we recall the definition of twisted $K$-theory  
 of a stack \cite{TXL04}.
Let $\gm$ be a locally compact groupoid with Haar system, $\uU$ an open
cover of $\gm\lcom$ and $c\in\check{Z}^2(\sigma\uU;\sS^1)$ a 2-cocycle.
Recall from Section~\ref{sec:4.2} that $c$ determines
a $C^*$-algebra $C^*_r(\gm,c)$.

\begin{defn}\label{def:twisted-kth}\cite{TXL04}
The twisted $K$-theory group $K^i_c (\gm )$
is defined as $K^i(C^*_r(\gm,c))$.
\end{defn}

\begin{rmk}
In \cite{TXL04}, we defined, for any 
 $\alpha\in H^2(\gm\lcom;\sS^1)$,
the twisted $K$-theory group $K^i_\alpha(\gm)$ as the $K$-theory group
of the $C^*$-algebra of any $S^1$-central extension representing $\alpha$.
This is indeed well-defined, as the $C^*$-algebras of two such central
extensions are Morita equivalent. However, the shortcoming of
such a   definition is that the Morita equivalence is not canonical (see
Proposition~\ref{pro:cocycle-iso}(3)).
On the other hand, 
Definition~\ref{def:twisted-kth} is canonical, as the 
$2$-cocycle $c$ determines
a unique $C^*$-algebra rather than  a  Morita equivalence class
of $C^*$-algebras.
This is extremely  important for our purpose in this paper
 since we need to study the functorial properties of twisted $K$-theory.
\end{rmk}

In particular, if $\gm\toto \gm_0$ is a Lie groupoid
acting on  a manifold $N$ with momentum map
 $J: N\to \gm_0$, and $c\in \check{Z}^2(\sigma\uU;\sS^1)$
is a 2-cocycle of the corresponding transformation groupoid
$N\rtimes \gm : N\rtimes \gm\toto N$. Then
the twisted equivariant $K$-theory group is defined as   follows \cite{TXL04}:
$$K^i_{c, \gm} (N)=K^i_c (N \rtimes \gm).$$

As an immediate consequence of Lemma \ref{lem:cocycle-iso-gpoid},
we have the following

\begin{prop}\label{pro:cocycle-iso2}
For any open covers $\uU$ of $\gm\lcom$,
$c^1$, $c^2\in \check{Z}^2(\sigma\uU;\sS^1)$, $b\in \check{C}^1(\sigma\uU;\sS^1)$
%\comment{maybe $\check{C}^2(\sigma \uU;\sS^1)$?
%There are a few places below that  $\sigma \uU$
%and $ \uU$ are mixed up. Perhaps we should just
%make a rmk that we use $\check{C}^2(\uU;\sS^1)$
%to denote  $\check{C}^2(\sigma \uU;\sS^1)$?}
satisfying $c^1-c^2=\del b$, there is a canonical
isomorphism
 $$\Phi^\uU_{c^2,b,c^1}: K^*_{c^1}(\gm)\to K^*_{c^2}(\gm), $$
which satisfies  the following properties:

\begin{enumerate}
\item if $c^1-c^2=\del b$ and $c^2-c^3=\del b'$,
then
$$\Phi^\uU_{c^3,b',c^2}\smalcirc \Phi^\uU_{c^2,b,c^1} = \Phi^\uU_{c^3,b+b',c^1}; $$
\item for any $a\in \check{C}^0(\uU, \sS^1)$,
 we have
$$\Phi^\uU_{c, \del a ,c}= {\mathrm{Id}}; $$
\item if $\uU'$ is a refinement of $\uU$ and ${c'}^1$, ${c'}^2$ and
$b'$ are the pull-backs of $c^1$, $c^2$ and $b$ under the canonical map
$\check{C}^*(\sigma\uU';\sS^1)\to \check{C}^*(\sigma\uU;\sS^1)$, then
$K^*_{c^j}(\gm)$ is canonically isomorphic to
$K^*_{{c'}^j}(\gm)$, and  the following diagram
$$\xymatrix{
K^*_{{c'}^1}(\gm) \ar[r]\ar[d]_{\Phi^{\uU'}_{{c'}^2,b',{c'}^1}} & 
K^*_{c^1}(\gm)\ar[d]^{\Phi^\uU_{c^2,b,c^1}} \\
K^*_{{c'}^2}(\gm)\ar[r] &   K^*_{c^2}(\gm)
}$$
commutes.
\end{enumerate}
\end{prop}
\begin{pf}
(1) and (3) are immediate consequences of  Lemma \ref{lem:cocycle-iso-gpoid}.
For (2), on the groupoid level, according to
 Lemma \ref{lem:cocycle-iso-gpoid},  
$$\Phi^\uU_{c,\del a,c}(\gamma)=\gamma
a(s(\gamma))a(t(\gamma))^{-1}\quad , \ \  \forall 
\gamma\in\tilde{\gm}^c.$$
This induces an automorphism of $C^*_r(\gm, c)$ given by
$\Phi(f)(\gamma)=a(t(\gamma))^{-1}f(\gamma)a(s(\gamma))$,
$\forall f\in C_c^\infty (\tilde{\gm}^c, L^c)$.
Let $(Uf)(\gamma)=a(t(\gamma))f(\gamma)$. Then $U$ is an
 unitary multiplier
of $C^*_r(\gm,c)$ and $\Phi(f)=U^*fU$. The conclusion thus follows.
\end{pf}

Recall that $\check{C}^n(\gm\lcom;\sS^1)=\lim_{\uU}\check{C}^n(\sigma\uU;\sS^1)$,
where $\uU$ runs over covers of the form
$(U^p_\gamma)_{\gamma\in\gm_p}$, $\gamma\in U^p_\gamma$,
and $\uU'$ is finer than $\uU$ if ${U'}^p_\gamma \subset U^p_\gamma$
for all $p$ and all $\gamma\in\gm_p$. We are now ready to 
introduce

\begin{defn}\label{def:twisted-kth2}
For any
 $c\in \check{Z}^2(\gm\lcom;\sS^1)$, we define
 the twisted $K$-theory
group $K^i_c(\gm)$ as $K^i_{c'}(\gm)$, where $c'\in Z^2(\sigma\uU';\sS^1)$
is a $2$-cocycle on a cover $\uU'$ of the above
 form which 
corresponds to $c$
 in $\check{C}^2(\gm\lcom;\sS^1)=\lim_{\uU}\check{C}^2(\sigma\uU;\sS^1)$.
\end{defn}

Note that this definition is valid  since
if $\uU''$ is another such
cover, then $K^i_{c'}(\gm)$ and $K^i_{c''}(\gm)$ are canonically
isomorphic according to Proposition~\ref{pro:cocycle-iso2}(3).
Now the following is an immediate consequence of
Proposition~\ref{pro:cocycle-iso2}.

\begin{prop}\label{pro:cocycle-iso}
Assume  that $c^1$ and $c^2$ are two $\sS^1$-valued \v{C}ech 2-cocycles
on a groupoid $\gm\toto \gm_0$, and $b$ is an $\sS^1$-valued 2-cochain
such that $c^1-c^2=\del b$. Then there is a canonical isomorphism
$$\Phi_{c^2,b,c^1}: K^*_{c^1}(\gm)\to K^*_{c^2}(\gm), $$
which satisfies the  following properties:
\begin{enumerate}
\item if $c^1-c^2=\del b$ and $c^2-c^3=\del b'$,
then 
$$\Phi_{c^3,b',c^2}\smalcirc \Phi_{c^2,b,c^1} = \Phi_{c^3,b+b',c^1}; $$
\item for any $a\in \check{C}^0(\gm\lcom, \sS^1)$,  we have
$$\Phi_{c,b,c}=\Phi_{c,b+\del a,c}.$$
\end{enumerate}
\end{prop}

As a consequence, we have

\begin{cor}
The map $[b]\mapsto \Phi_{c,b,c}$ defines
a morphism $H^1(\gm;\sS^1)\to\mbox{Aut}(K^*_c(\gm))$.
\end{cor}

For simplicity, in the sequel
 we will  write $\Phi_b$ instead of $\Phi_{c^2,b,c^1}$
 whenever there is no ambiguity.

\subsection{External Kasparov product}

Let us first recall a few basic facts about
equivariant $KK$-theory, as introduced by Kasparov \cite{Kas88}
and generalized by Le Gall to groupoids \cite{Leg}.

Assume that $\gm$ is a locally compact $\sigma$-compact groupoid
with Haar system (for instance, a Lie groupoid), $A$, $A_1$, $B$, $B_1$, and
$D$ are $\gm-C^*$-algebras,
 i.e $C^*$-algebras endowed with an action of
$\gm$. Then there is a 
 bifunctor 
$$(A,B)\mapsto KK_\gm^i(A,B), \ \ \ i=0,1$$
 covariant in $B$ and  contravariant in $A$.
There is also a suspension map 
$$\sigma_{\gm_0,D}: KK_\gm^i(A,B)\to KK_\gm^i(D\otimes_{C_0(\gm_0)}A,
D\otimes_{C_0(\gm_0)}B), $$
 and an associative product
\begin{eqnarray}
KK_\gm^i(A,B)\otimes KK_\gm^j(B,C)&\to &
 KK_\gm^{i+j}(A,C), \label{eq:abb1}\\
(\alpha,\beta)&\mapsto& \alpha\otimes_{B}\beta.
\end{eqnarray}
More generally, the map
$$KK_\gm^i(A,B\otimes_{C_0(\gm_0)} D)\otimes KK_\gm^j(A_1 \otimesc D,B_1)
\to KK_\gm^{i+j}(A\otimesc A_1, B\otimesc B_1)$$
defined by
 $$\alpha\otimes_D\beta:=\sigma_{\gm_0,A_1}(\alpha)
 \otimes_{A_1\otimesc B\otimesc D}\sigma_{\gm_0,B}(\beta)$$
 is associative. In particular, when $D=C_0(\gm_0)$,
we obtain an associative and (graded) commutative  product
\begin{equation}
\label{eq:k}
KK_\gm^i(A,B)\otimes KK_\gm^j(A_1,B_1)\to KK_\gm^{i+j}(A\otimesc A_1,
B\otimesc B_1),
\end{equation}
which is called the {\em external Kasparov product}.

We also recall \cite{bch94,tu00}
that if $A$ is a $\gm-C^*$-algebra, there is a
group morphism, the Baum--Connes assembly map
$$\mu_r: K_i^{top}(\gm;A)\to K_i(A\rtimes_r\gm).$$
Here the left-hand side is  defined by
$\lim_{Y\subset\underline{E}\gm} KK_\gm^i(C_0(Y),A)$,
where $Y$ runs over saturated subsets of a certain proper $\gm$-space
$\underline{E}\gm$ such that $Y/\gm$ is compact (if $\gm$ is proper,
then $\underline{E}\Gamma=\Gamma_0$). The map $\mu_r$ is
defined by $\mu_r(x)=\lambda_{Y\rtimes \gm}\otimes_{C_0(Y)\rtimes_r\gm}
j_{\gm,r}(x)$ for all $x\in KK_\gm^i(C_0(Y),A)$, where
$$j_{\gm,r}:KK_\gm^i(C_0(Y),A)\to KK^i(C_0(Y)\rtimes_r\gm,A\rtimes_r\gm)$$
is the Kasparov's descent morphism and $\lambda_{Y\rtimes\gm}
\in KK(\cc,C_0(Y)\rtimes_r\gm)=K_0(C^*_r(Y\rtimes\gm))$ is a canonical
element.

The Baum--Connes assembly map is an isomorphism
when $\gm$ is proper, or more generally, amenable \cite{tu99a}
(as well as in many other cases such as connected Lie groups \cite{cen}).
One does not however expect $\mu_r$ to be an isomorphism for every
Hausdorff Lie groupoid, in view of the counterexamples in \cite{hls}
(these counter examples are, however, neither Hausdorff nor  smooth).

%\comment{do we need this map commute with  Gysin map and $\Phi_b$?
%Perhaps we need to prove  a similar result like 
%Lemma 4.1 for $A_c$ (or the groupoid extension which 
%gives rise to $A_c$). See remark
%below}

\begin{prop}
\label{pro:4.13}
Assume  that $\gm\toto \gm_0$ is a proper Lie  groupoid (or, more generally,
a groupoid such that the Baum--Connes map is an isomorphism)
acting on a space $N$.
Let $c\in \check{Z}^2(\uU, \sS^1 )$ be any 2-cocycle,
where $\uU$ is a cover of $(N\rtimes\gm)\lcom$ and
$N\rtimes \gm $ is the transformation groupoid.
Then there is a canonical isomorphism
\begin{equation}
K_i^{top}(\gm; A_c) \longiso K_{c, \gm}^i (N),
\end{equation}
where  $A_c$  is the $C^*$-algebra defined in Corollary~\ref{cor:Ac}.
Moreover, under the same assumptions as in Proposition~\ref{pro:cocycle-iso},
denote  again by $\alpha_{c^2,b,c^1}$ the element  in
$KK_\gm(A_{c^1},A_{c^2})$ determined by the $\gm$-equivariant
Morita equivalence $\alpha_{c^2,b,c^1}:A_{c^1}\to A_{c^2}$. Then
we have a commutative diagram
$$\xymatrix{
K_i^{top}(\gm;A_{c^1})\ar[d]\ar[r]^\cong & K_{c^1,\gm}^i(N)
\ar[d]^{\Phi_{c^2,b,c^1}}\\
K_i^{top}(\gm;A_{c^2})\ar[r]^\cong & K_{c^2,\gm}^i(N),
}$$
where the left vertical arrow is the multiplication by the element
$\alpha_{c^2,b,c^1}$.
\end{prop}
\begin{pf}
We have
$$K_i^{top}(\gm,A_c)\cong K_i(A_c\rtimes_r\gm)
\cong  K_i(C^*_r(N\rtimes\gm,c))= K_{c,\gm}^i(N), $$
where the first isomorphism is the
Baum--Connes assembly map $\mu_r$,  and the second one follows from
Corollary~\ref{cor:Ac}.
Note that each isomorphism above is canonical.
That the diagram commutes follows from the equality
$\mu_r(x\otimes \alpha_{c^2,b,c^1})=j_{\gm,r}(\lambda_{Y\rtimes\gm}\otimes
x\otimes \alpha_{c^2,b,c^1})=j_{\gm,r}(\lambda_{Y\rtimes\gm}\otimes
x)\otimes j_{\gm,r}(\alpha_{c^2,b,c^1})$  together with Corollary~\ref{cor:Ac}.
\end{pf}

If $c\in \check{Z}^2(\gm\lcom;\sS^1)$ lives on a cover
$\uU$,
 and  $c' \in  \check{Z}^2(\sigma \uU;\sS^1) $ is the corresponding
2-cocycle, then $K^i_{c,\gm}(N)=K_i(A_{c'} \rtimes \gm)$
is canonically defined. 
Below, we will write $A_c$ instead of $A_{c'}$ whenever there is no
ambiguity  (keeping in mind that the $\gm$-$C^*$-algebra $A_c$ is
determined up to a unique Morita equivalence).

Suppose now that $\gm$ is proper. Let $N_2= N \times_{\gm_0} N$,
 and $p_i:N_2\to N$, $i=1, 2$ be the natural projections.
Take  $A=C=C_0(Y)$, and  $B=D=A_c$. Thus we have
$B\otimes_{C_0(\gm_0)} D\cong A_{p_1^*c+p_2^*c}$
according to  Proposition~\ref{prop:p1cp2c}. 
Thus  Proposition \ref{pro:4.13} together with
the external Kasparov product
Eq. \eqref{eq:k} implies the following

\begin{prop}
\label{pro:4.18}
Assume  that $\gm\toto \gm_0$ is a proper Lie  groupoid
acting on a space $N$. Let $c\in \check{Z}^2(\uU, \sS^1 )$ be any 2-cocycle,
where $\uU$ is a cover of $(N\rtimes\gm)\lcom$ and
$N\rtimes \gm $ is the transformation groupoid.
We have a map
\begin{equation}
K^i_{c,\gm}(N)\otimes K^j_{c,\gm}(N)\to
K^{i+j}_{p_1^*c+p_2^*c,\gm}(N_2) .
\end{equation}
\end{prop}

\subsection{Gysin maps}
%Construction of the ring structure}

%\begin{lem}
%If $\psi:\gm \to\gm'$ is a groupoid isomorphism,
%$c'\in \check{C}^2(\gm';\sS^1)$, $c=\psi^*c'\in
%\check{C}^2(\gm;\sS^1)$, then $\psi$ defines a canonical map
%$\psi^*: K^*_{c'}(\gm'\lcom)\to K^*_{c}(\gm\lcom)$ such that
%the relation $(\psi\circ\psi')^*=(\psi')^*\circ \psi^*$ holds.
%\end{lem}
%
%\begin{pf}
%Obvious.
%\end{pf}

Finally, we will need Gysin (wrong-way functoriality) maps.
Such a map was studied extensively in \cite{CS84,HS87}.
Recall  that for manifolds $M$ and $N$, a  smooth map $f:M\to N$
 is said to be $K$-oriented if the normal bundle
$N_f=T^*M\oplus f^*(TN)\to M$ is $K$-oriented, i.e. it admits a
$\mbox{Spin}^c$-structure. When $f$ is a submersion,
$N_f$ can  be replaced by $\ker(df)\subseteq TM$ \cite{CS84}.

%Recall that if $M$ is a manifold endowed with an action of a Lie groupoid
%$\gm$ and if $E$ is a $\gm$-equivariant vector bundle on $M$,
%then $E$ is said to be $K$-oriented if it admits a
%$\gm$-equivariant $\mbox{Spin}^c$-structure.

Now let  $\gm$ be a Lie  groupoid, and both $M$ and  $N$ be    $\gm$-manifolds.
It is still  not clear what will be  an appropriate  notion
of $\gm$-equivariantly $K$-orientability
for a general $\gm$-equivariant smooth map $f:M\to N$,
 since  the normal bundle of $f$ does not necessarily
  admit   a $\gm$-action in general.
However  the  $\gm$-equivariant
$K$-orientability does make sense  in any of the following
 special  cases:

\begin{itemize}
\item[(a)] $\gm$ possesses a pseudo-\'etale structure
in the sense of ~Tang \cite{Tang} (which is also called
a flat structure by Behrend \cite{Kai}). In particular,
this includes the case that
 $\gm$ is \'etale, or is  a transformation groupoid  corresponding
to a Lie group   action;
\item[(b)] $f$ is a submersion (in this  case, $\ker(df)$ is endowed
with an action of $\gm$);
\item[(c)] both  momentum maps $J:M\to \gm_0$
and $J':N\to \gm_0$ are submersions. Denote by $VM=
\ker J_*:TM\to T\gm_0$ and $VN=\ker J'_*:TN\to T\gm_0$ the vertical
tangent bundles. Then  $f$ is said to be $K$-oriented if  the vector
bundle $V^*M\oplus f^* VN\to M$ admits a $\gm$-equivariant
 $K$-orientation.
\end{itemize}

From now on, we will assume  that $f:M\to N$ is a
 $\gm$-equivariantly $K$-oriented
submersion as in (b). Indeed, although Gysin maps in twisted
 $K$-theory may be constructed in other cases as well (see
for instance \cite{CW}), we will only develop the case of submersions
 since other cases are not needed in this paper.

As we will see,
 $f$ induces an element $f_!\in KK_\gm^d(C_0(M),C_0(N))$, where
$d=\dim M-\dim N $,
such that the relation $(g\smalcirc f)_! = g_!\smalcirc f_!$ holds.
More generally, we have

\begin{lem}\label{prop:Gysin}
Let $\gm$ be a Lie groupoid,  and $M$, $N$ two proper $\gm$-manifolds. 
 Assume that  $f:M\to N$ is a $\gm$-equivariantly  
 $K$-oriented submersion.
Then for any  $\sS^1$-valued \v{C}ech 2-cocycle $c$
on the groupoid $N\rtimes \gm$,  there is
a Gysin element
$f^c_!\in KK_\gm^d(A_{f^*c},A_c)$,  where $d=\dim N-\dim M$,
satisfying the property:
\begin{equation}
\label{eq:com}
f^{g^*c}_!\otimes_{A_{g^*c}} g_!^c = (g\smalcirc f)^c_!,
\end{equation}
for any $\gm$-equivariant  $K$-oriented maps
$f:M\to N$ and $g:N\to P$. The $K$-orientation of $g\circ f$ is
induced from that of $f$ and $g$. 
Here both sides are considered
as elements in $KK_\gm^{d''}(A_{(g\circ f)^*c},A_c)$ and
$d''=\dim P-\dim M$.
\end{lem}
\begin{pf}
It is standard that any $K$-oriented submersion
 $f: M\to N$ yields  a Gysin element
$f_!\in KK^d(C_0(M),C_0(N))$ \cite{CS84,HS87}.
When  $\gm$ is a  Lie group, an equivariant version
was  proved by Kasparov-Skandalis \cite[\S 4.3]{KS91}:
 any  $\gm$-equivariantly $K$-oriented map $f: M\to N$
 determines an element $f_!\in KK^d_\gm(C_0(M),C_0(N))$.
A similar argument can be adapted to show that the same assertion holds
when $\gm$ is a Lie groupoid, and  $KK_\gm$ is Le Gall's groupoid
equivariant $KK$-theory \cite{Leg}.
In fact, since $f$ is also $N\rtimes \gm$-equivariant, it is easy to see that
one obtains an element
$f_!\in KK^d_{N\rtimes\gm}(C_0(M),C_0(N))$.

Consider 
$$\sigma_{N,A_c}(f_!) \in KK^d_{N\rtimes \gm}
(C_0(M)\otimes_{C_0(N)} A_c,C_0(N)\otimes_{C_0(N)}A_c).$$
By Proposition \ref{prop:fAc}, we have $C_0(M)\otimes_{C_0(N)} A_c
\cong A_{f^*c}$. Thus $\sigma_{N,A_c}(f_!) \in
KK^d_{N\rtimes\gm}(A_{f^*c},A_c)$.
 We define $f^c_!$ as the image
of this element under
 the forgetful functor $KK_{N\rtimes\gm}\to KK_\gm$.

We need to check that $f^c_!$ is well-defined.
Recall that $A_c$ is defined only up to Morita equivalence.
 In the definition of $A_c$, we have implicitly chosen an open cover $\uU$,
 on which the $2$-cocycle $c$ lives.
Assume  that $c'$ is another such a $2$-cocycle,
 which is defined
on another cover $\uU'$. Let $\alpha\in KK_{N\rtimes\gm}(A_c,A_{c'})$
be the element determined by the Morita equivalence between
$A_c$ and $A_{c'}$. Since all the constructions are natural,
$f^*\alpha\in KK_{M\rtimes\gm}(A_{f^*c},A_{f^*c'})
=KK_{M\rtimes\gm}(C_0(M)\otimes_{C_0(N)}A_c,C_0(M)\otimes_{C_0(N)}A_{c'})$
is the  element induced by
the Morita equivalence between $A_{f^*c}$ and $A_{f^*(c')}$.

We need to show that $f^c_!$ and $f^{c'}_!$ can be identified, i.e.
 $$\sigma_{N,A_c}(f_!)\otimes_{A_c}\alpha
=f^*\alpha\otimes_{f^*A_{c'}}\sigma_{N,A_{c'}}(f_!)\in
KK_{N\rtimes\gm}^*(f^*A_c,A_{c'}). $$
The left hand side is $f_!\otimes_{C_0(N)}\alpha$, while the
right hand side is equal to 
 $\sigma_{N,C_0(M)}(\alpha)\otimes_{C_0(M)\otimes_{C_0(N)}A_{c'}}
\sigma_{N,A_{c'}}(f_!)=\alpha\otimes_{C_0(N)} f_!$.
Thus the equality follows from the commutativity 
of the external Kasparov product.

Eq. \eqref{eq:com} now
 follows from the compatibility of the suspension
$\sigma$ with the Kasparov product
$$f_!^{g^*c}\otimes_{A_{g^*c}} g_!^c=\sigma_{P,A_c}(f_!)\otimes_{A_{g^*c}}
\sigma_{P,A_c}(g_!)=\sigma_{P,A_c}(f_!\otimes_{C_0(N)}g_!)
=\sigma_{P,A_c}((g\smalcirc f)_!)=(g\smalcirc f)^c_! .$$
\end{pf}

\begin{cor}
\label{cor:Gysin}
Under the same hypothesis as in Lemma \ref{prop:Gysin},
there is a Gysin map
$$f_!:K^i_{f^*c, \gm}(M)\to K^{i+d}_{c, \gm} (N),$$
with $d= \dim N-\dim M$,
which satisfies $g_!\smalcirc f_! = (g\smalcirc f)_!$.
\end{cor}
\begin{pf} We have
$f_! (\beta )=\beta\otimes_{f^*A_c} f_!^c$.
\end{pf}
The following proposition describes the
naturality property
of the Gysin map with respect to the cocycle $c$.

%With the same assumptions, suppose that $c'$ is a cocycle which is
%cohomologous to $c$, i.e. $c-c'=\del u$ for some $u\in \check{C}^1(
%(N_1\rtimes\gm)\lcom;\sS^1)$. We have the following naturality property
%of the Gysin map with respect to the cocycle $c$:

\begin{prop}\label{prop:gysin-phi}
Under the same hypothesis as in Corollary \ref{cor:Gysin},
assume that $c'$ is another $2$-cocycle which is
cohomologous to $c$, i.e. $c-c'=\del u$ for some $u\in \check{C}^1(
(N\rtimes\gm)\lcom;\sS^1)$. Then
the diagram
$$\xymatrix{
K_{f^*c,\gm}^i (M)\ar[r]^{f_!}\ar[d]^{\Phi_{f^*u}}
 & K_{c,\gm}^{i+d}(N)\ar[d]^{\Phi_u}\\
K_{f^*c',\gm}^i (M)\ar[r]^{f_!}
 & K_{c',\gm}^{i+d}(N)
}$$
 commutes.
\end{prop}

\begin{pf}
We can assume that the relation $c-c'=\del u$ holds in some
complex $C^*(\sigma\uU;\sS^1)$. Therefore
 the isomorphism $\Phi_u$ is implemented
by an isomorphism of  $S^1$-central extensions,  and
thus by an invertible element
$\alpha\in KK_{N\rtimes\gm}(A_c,A_{c'})$. It is clear,
 by  naturality of the constructions,
 that $\Phi_{f^*u}$ is implemented by
$f^*\alpha\in KK_{M\rtimes\gm}(A_{f^*c},A_{f^*c'})$.
Then for all $\beta\in K^*_{f^*c,\gm}(M)$, we have
$(\Phi_u\smalcirc f_! )(\beta)=\beta\otimes_{f^*A_c}f^c_!\otimes_{A_c}\alpha$,
while $(f_!\smalcirc \Phi_{f^*u})(\beta)=\beta\otimes_{f^*A_c}
f^*\alpha\otimes_{f^*A_{c'}}f^{c'}_!$. Thus
 we need the identity $\sigma_{N,A_c}(f_!)\otimes_{A_c}\alpha
=f^*\alpha\otimes_{f^*A_{c'}}\sigma_{N,A_{c'}}(f_!)\in
KK_{N\rtimes\gm}^*(f^*A_c,A_{c'})$, which
can be proved   exactly in  the same way  as
in Proposition~\ref{prop:Gysin}.
\end{pf}

%\begin{tiny}
%It is thus important to know whether the
%  manifolds involved admit ${\mathrm{Spin}^c}$
%structures.
%
%Until the end of the section, $G$ is a compact, connected, simply
%connected and simple Lie group and $N=G$ with $G$ acting on $N$
%by conjugation. Note that every cohomology class in
%$H^2((N\rtimes G)\lcom;\sS^1)$ can be represented by a
%$2$-multiplicative cocycle $c$ since $T_1:H^3(G\lcom;\sS^1)\to
%H^2((N\rtimes G)\lcom;\sS^1)$ is an isomorphism
%(Corollary~\ref{cor:T1iso}).
%\begin{lem}\label{lem:spinc}
%$N$ admits a $G$-equivariant Spin structure (and thus a $G$-equivariant
%Spin${}^c$-structure).
%\end{lem}
%\begin{pf}
%The tangent bundle of $N$ is $G$-equivariantly isomorphic to $N\times \GG$,
%where $\GG$ is the Lie algebra of $G$ endowed with the adjoint action.
%Endow $\GG$ with an ${\mathrm{Ad}}$-invariant scalar product.
%Then $G$ acts on $N$ by isometries which preserve the orientation
%(since $G$ is connected), thus $G$ acts on the direct
%orthonormal frame bundle $N\times SO(\GG)$ by $g\cdot (x,T)
%=(gxg^{-1},Ad(g)T)$. Moreover,
%since $G$ is simply connected, the group morphism $Ad:G\to SO(\GG)$
%lifts to the universal cover $\widetilde{SO(\GG)}$.
%\end{pf}
%\end{tiny}

\subsection{Ring structure on equivariant twisted $K$-theory group}

%A groupoid $\gm$ is said to be
% longitudinally $K$-oriented
%if $s:\gm_1\to\gm_0$ is $K$-oriented. \comment{are you sure this
%is enough? any compatibility condition with the
%groupoid structure?}
%According to \cite[Proposition B.6]{CS84},  this is also equivalent to
%that the subbundle $\ker(s^*)$ of $T\gm_1\to\gm_1$ admits a
%$\mbox{Spin}^c$-structure, which means that each manifold $\gm_x$
%is $K$-oriented, and the $K$-orientations vary smoothly with respect
%to $x\in \gm_0$. \comment{do we require $\mbox{Spin}^c$-structure
%invriant under the left translation}

Let $N\stackrel{\varphi}{\to} \gm$ be a crossed module,
 where $\gm\toto \gm_0$ is a proper  Lie groupoid such that
$s:N\to N_0$ is $\gm$-equivariantly $K$-oriented.
%This means that there is a $\mbox{Spin}^c$-structure on each manifold
%$N_x$ ($x\in N_0$) which varies smoothly with respect to $x$, and which
%is $\gm$-equivariant in the sense that for each $g\in\gm$, the
%$\mbox{Spin}^c$-structures on $N_{s(g)}$ and on $N_{t(g)}$ are
%the same, once $N_{s(g)}$ and $N_{t(g)}$ are identified via
%the diffeomorphism determined by $g$.
Assume that $(c, b, a)$ is a multiplicator  as in  Definition
\ref{def:multiplicator}.

\begin{defn}\label{def:product}
Define
$$K^{i+d}_{c,\gm}(N)\otimes K^{j+d}_{c,\gm}(N)\to K^{i+j+d}_{c,\gm}(N),$$
where $ d=\dim N -\dim N_0$,
as the composition of the external Kasparov product
$K^{i+d}_{c,\gm}(N)\otimes K^{j+d}_{c,\gm}(N)
\to K^{i+j}_{p_1^*c+p_2^*c,\gm}(N_2)$ as in Proposition~\ref{pro:4.18}
with  the maps
$$K^{i+j}_{p_1^*c+p_2^*c,\gm}(N_2)
\stackrel{\Phi_b}{\longrightarrow} K^{i+j}_{m^*c,\gm}(N_2)
\stackrel{m_!}{\longrightarrow}K^{i+j+d}_{c,\gm}(N).$$
Here $m_!$ is the Gysin map  corresponding to $m: N_2\to N$.
\end{defn}

Note that in the above definition if $c$ is defined on some cover $\uU$,
then $p_1^*c$ and $p_2^*c$ are defined on different covers.
Therefore to define $p_1^*c+p_2^*c$,
 one needs to pass to a refinement. This
is the reason why we choose to define
twisted $K$-theory with respect to a \v{C}ech cocycle
$c\in \check{C}^2((N\rtimes\gm)\lcom;\sS^1)$ as in
Definition~\ref{def:twisted-kth2} instead of a cocycle
in $\check{C}^2(\sigma\uU;\sS^1)$ for some fixed cover.

Also, note that  $s:N\to N_0$ is $\gm$-equivariantly
$K$-oriented implies that $pr_2: N\times_{N_0}N\to N$
is $\gm$-equivariantly $K$-oriented. Since $N_2$ is   $\gm$-equivariantly
diffeomorphic to $N\times_{N_0}N$ by 
$\psi:(x,y)\mapsto (x,xy)$, it follows that
$m=pr_2\circ\psi:N_2\to N$ is  $\gm$-equivariantly $K$-oriented.
Moreover, define maps $N_3\to N_2$ by $m_{12}(x,y,z)=(xy,z)$
and $m_{23}(x,y,z)=(x,yz)$. It is not hard to check from the assumptions
that the $K$-orientations of $m\smalcirc m_{12}$ and $m\smalcirc m_{23}$
coincide. Let $m_{123}=m \smalcirc m_{12}=m\smalcirc m_{23}$ be  endowed
with this $K$-orientation, then we obtain  the following

\begin{prop}
\label{pro:m}
$${m_{12}}_!\smalcirc m_!={m_{123}}_!={m_{23}}_!\smalcirc m_!.$$
\end{prop}

\begin{rmk}\label{rmk:non-connected}
In Definition~\ref{def:product} above, if $N$ is not connected, 
$K^{i+d}_{c,\gm}(N)$ needs to  be replaced by
$\oplus_l K^{i+d_l}_{(c,\gm)}(N^{[l]})$, where
 $N=\coprod_l N^{[l]}$ is  the  decomposition into
its  connected components, and  $d_l$ is the dimension of
$N^{[l]}$. An example where a disconnected crossed module naturally
appears is when $N$ is the space of closed loops
$S\gm=\{g\in\gm\vert\; s(g)=t(g)\}$,
and $\gm$ is an \'etale proper groupoid.
 The results formulated
below are still valid in this context.
\end{rmk}

Now we can state the main theorem of this paper.

\begin{them}
Let $N\stackrel{\varphi}{\to} \gm$ be a crossed module,
 where $\gm\toto \gm_0$ is a proper  Lie groupoid such that
$s:N\to N_0$ is $\gm$-equivariantly  $K$-oriented.
Assume that $(c, b, a)$ is a multiplicator,
where $c\in \check{C}^2((N\rtimes \gm)\lcom, \sS^1)$,
$b\in \check{C}^1((N_2\rtimes \gm)\lcom, \sS^1)$ and
$a\in \check{C}^0((N_3\rtimes \gm)\lcom, \sS^1)$.
Then the product 
$$K^{i+d}_{c,\gm}(N)\otimes K^{j+d}_{c,\gm}(N)\to K^{i+j+d}_{c,\gm}(N)$$
is associative.
\end{them} 

%\begin{tiny}
%Before we show associativity, let us prove
%Proposition~\ref{prop:Psi}.
%
%Let $c=T_1e$ and $c'=T_1e'$. Since $\del$ and $T_1$ anti-commute
%(see (\ref{eqn:delT})), we have $c-c'=-\del T_1 u$. Define
%$$\Psi_{e,u,e'}=\Phi_{c,-T_1u,c'}. $$
%All the assertions of Proposition~\ref{prop:Psi},
%except the fact that $\Psi_{e,u,e'}$ is a ring morphism,
%follow from Lemma~\ref{lem:cocycle-iso}.
%
%It remains to check that $\Psi_{e,u,e'}$ preserves the product.
%We thus have to check that the following diagram commutes:
%$$\xymatrix{
%K^*_{p_1^*c+p_2^*c,G}(N^2)
%       \ar[r]^{\Phi_b}\ar[d]^{\Phi_{-p_1^*T_1u-p_2^*T_1u}}
%  & K^*_{m^*c,G}(N^2)\ar[r]^{m!}\ar[d]^{\Phi_{-m^* T_1 u}}
%  & K^*_{c,G}(N)\ar[d]^{\Phi_{-T_1^* u}}\\
%K^*_{p_1^*c'+p_2^*c',G}(N^2)
%        \ar[r]^{\Phi_{b'}}
%  & K^*_{m^*c',G}(N^2)\ar[r]^{m!}
%  & K^*_{c',G}(N)
%}$$
%
%For the second square, this is easy. \comment{SPELL IT OUT}
% For the first square,
%it suffices, using  Lemma~\ref{lem:cocycle-iso},  to check that
%$b-m^*T_1 u$ and $-p_1^* T_1u-p_2^*T_1u+b'$ differ by a \v{C}ech
%coboundary. Now, using from (\ref{eqn:delT}) the relation
%$\del' T_1=T_2\del-\del T_2$, we get
%\begin{eqnarray*}
%\lefteqn{(b-m^* T_1u)-(-p_1^*T_1 u -p_2^*T_1u+b')
%= b-b'+\del'T_1 u}\\
%& =& -T_2(e-e')+\del'T_1u
%=-T_2\del u +\del'T_1u=-\del T_2 u.
%\end{eqnarray*}
%This completes the proof of Proposition~\ref{prop:Psi}.
%\end{tiny}
\begin{pf}
Let us study the product $(x,y,z)\mapsto x(yz)$ from
$K^{i+d}_{c,\gm}(N)\otimes K^{j+d}_{c,\gm}(N)\otimes K^{k+d}_{c,\gm}(N)$
to $K_{c,\gm}^{i+j+k+d}(N)$.
To simplify notations, we assume $i=j=k=0$.
The product $x(yz)$ is obtained by composing the external product
$K_{c,\gm}^d(N)\otimes K_{c,\gm}^d(N)\otimes K_{c,\gm}^d(N)
\to K^{d}_{p_1^*c+p_2^*c+p_3^*c,\gm}(N_3)$ with
$\Phi_{p_{23}^*b}$,
followed by the down-right composition of the diagram below:
$$\xymatrix{
K^d_{p_1^*c+m_{23}^*c,\gm}(N_3)\ar[r]^{\Phi_{m_{23}^*b}}\ar[d]_{(m_{23})_!}
  & K^d_{m_{123}^*c,\gm}(N_3)\ar[d]^{(m_{23})_!}\ar[dr]^{(m_{123})_!} & \\
K^0_{p_1^*c+p_2^*c,\gm}(N_2)\ar[r]^{\Phi_b}
  & K^0_{m^*c,\gm}(N_2)\ar[r]^{m_!}& K^d_{c,\gm}(N).
}$$
(Recall that $m_{23}(x,y,z)=(x,yz)$, $m_{12}(x,y,z)=(xy,z)$,
$m_{123}(x,y,z)=xyz$.)
\par\medskip

It follows from Proposition~\ref{prop:gysin-phi}, Corollary~\ref{cor:Gysin}
and Proposition~\ref{pro:m} that the  diagram above commutes.
 Hence
the product $x(yz)$ can also be obtained from
the composition of the external product
$K_{c,\gm}^d(N)\otimes K_{c,\gm}^d(N)\otimes K_{c,\gm}^d(N)
\to K^d_{p_1^*c+p_2^*c+p_3^*c,\gm}(N_3)$
with the maps
$$K_{p_1^*c+p_2^*c+p_3^*c,\gm}^d(N_3)
\stackrel{\Phi_{p_{23}^*b+m_{23}^*b}}{\longrightarrow}
K^d_{m_{123}^*c,\gm}(N_3)
\stackrel{(m_{123})_!}{\rightarrow}K^d_{c,\gm}(N).$$
\par\medskip

Similarly, the product $(xy)z$ is obtained by
composing the external product
$K_{c,\gm}^d(N)\otimes K_{c,\gm}^d(N)\otimes K_{c,\gm}^d(N)
\to K^d_{p_1^*c+p_2^*c+p_3^*c,\gm}(N_3)$
with  the maps
$$K_{p_1^*c+p_2^*c+p_3^*c,\gm}^d(N_3)
\stackrel{\Phi_{p_{12}^*b+m_{12}^*b}}{\longrightarrow}
K^d_{m_{123}^*c,\gm}(N_3)
\stackrel{(m_{123})_!}{\rightarrow}K^d_{c,\gm}(N).$$
Now since $p_{23}^*b+m_{23}^*b -
 p_{12}^*b+m_{12}^*b=\del' b =\del a$, it follows
from  Proposition~\ref{pro:cocycle-iso}  that
$\Phi_{p_{23}^*b+m_{23}^*b}=\Phi_{p_{12}^*b+m_{12}^*b}$.
 This completes the proof of the theorem.
\end{pf}

\subsection{Ring structure on  the $K$-theory group twisted by 2-gerbes}

%As an application of Theorem \ref{}, we consider twisted $K$-theory group of
%an inertia groupoid.  Let $G_1\toto G_0$ be a
%Lie groupoid and $\Lambda G: SG_1\rtimes G_1\toto SG_1\rtimes$
%be its inertia groupoid.   Due to the transgression map
%$T^1: H^3(G\lcom, \sS^1)\to H^2(\Lambda G\lcom, \sS^1)$,
%any element in the image   of $T_1$ is $2$-multiplicative.
%Thus one has a ring structure on the twisted $K$-theory group.
%Since $H^3(G\lcom, \sS^1)$ classifies $2$-gerbes, 
%we conclude that the twisted $K$-theory group on the inertia
%stack twisted by   a $2$-gerbe over the stack admits
%a ring structure.

Thanks to the transgression maps, for any crossed
module $N\to \gm$,  one can produce a canonical
multiplicator from $\check{Z}^3(  \gm\lcom;\sS^1)$.
Therefore one obtains  a ring structure on the corresponding
twisted $K$-theory group.

More precisely, given any $e\in \check{Z}^3( \gm\lcom;\sS^1)$,
$T_1 e\in \check{Z}^2( (N\rtimes\gm)\lcom;\sS^1)$ is 2-multiplicative.
Define
$$K_{e, \gm}^*(N): =K_{T_1 e, \gm}^*(N). $$
Let $c=T_1 e$, $b=-T_2 e$ and $a=-T_3e$. 
Then $(c, b, a)$ is a multiplicator according to Corollary
\ref{cor:3.11}.  Thus  $K_{e, \gm}^*(N)$ naturally admits
 a ring structure.

\begin{them}
\label{thm:main}
Let $N\stackrel{\varphi}{\to} \gm$ be a crossed module,
where $\gm\toto \gm_0$ is a proper  Lie groupoid such that
$s:N\to N_0$ is $\gm$-equivariantly $K$-oriented.
Let $d=\dim N-\dim N_0$ (see also Remark~\ref{rmk:non-connected}).
\begin{enumerate}
\item
For any $e\in \check{Z}^3( \gm\lcom;\sS^1)$, the
twisted  $K$-theory group ${K}^*_{e,\gm}(N)$ 
 is endowed with a ring structure
$$K^{i+d}_{e,\gm}(N)\otimes K^{j+d}_{e,\gm}(N)\to K^{i+j+d}_{e,\gm}(N),$$
where $d=\dim N -\dim N_0$.
\item Assume that $e$ and  $e'\in \check{Z}^3(\gm\lcom;\sS^1)$ satisfy
$e-e'=\del u$ for some $u\in \check{C}^2(\gm\lcom;\sS^1)$.
Then there is a ring isomorphism
$$\Psi_{e',u,e}:K^*_{e,\gm}(N)\to K^*_{e',\gm}(N)$$
such that 
\begin{itemize}
\item  if $e-e'=\del u$ and $e'-e''=\del u'$, then
$$\Psi_{e'',u',e'}\smalcirc\Psi_{e',u,e}=\Psi_{e'',u+u',e};$$
\item  for any $v\in  \check{C}^1(\gm\lcom;\sS^1)$,
$$\Psi_{e',u,e}=\Psi_{e',u+\del v,e}.$$
\end{itemize}
\item
% As a consequence, in particular when taking $e=e'$,
  There is a morphism
$$H^2(\gm\lcom;\sS^1)\to \mathop{Aut} K^*_{e,\gm}(N).$$
The ring structure
 on  $K^{*+d}_{e,\gm}(N)$, up to an isomorphism, 
 depends only  on the cohomology class
$[e]\in H^3(\gm\lcom;\sS^1)$.
The isomorphism is unique up to an automorphism of
$K^{*+d}_{e,\gm}(N)$ induced   from $H^2(\gm\lcom;\sS^1)$.
\end{enumerate}
\end{them}

\begin{pf}
Let $c=T_1e$ and $c'=T_1e'$. Since $\del$ and $T_1$ anti-commute
according to Eq.  (\ref{eq:T1par}), we have $c-c'=-\del T_1 u$. Define
$$\Psi_{e',u,e}=\Phi_{c',-T_1u,c}. $$
All the assertions of the theorem
except  that $\Psi_{e',u,e}$ is a ring morphism,
follow from Proposition~\ref{pro:cocycle-iso}.

It remains to prove that $\Psi_{e',u,e}$ preserves the product.
We thus need to check that the following diagram:
$$\xymatrix{
K^*_{c,\gm}(N)\otimes K^*_{c,\gm}(N)\ar[r]\ar[d]^{\Phi_{-T_1u}\otimes
     \Phi_{-T_1u}}
  &K^*_{p_1^*c+p_2^*c,\gm}(N_2)
       \ar[r]^{\Phi_b}\ar[d]^{\Phi_{-p_1^*T_1u-p_2^*T_1u}}
  & K^*_{m^*c,\gm}(N_2)\ar[r]^{m!}\ar[d]^{\Phi_{-m^* T_1 u}}
  & K^*_{c,\gm}(N)\ar[d]^{\Phi_{-T_1 u}}\\
K^*_{c',\gm}(N)\otimes K^*_{c',\gm}(N)\ar[r]
  & K^*_{p_1^*c'+p_2^*c',\gm}(N_2)
        \ar[r]^{\Phi_{b'}}
  & K^*_{m^*c',\gm}(N_2)\ar[r]^{m!}
  & K^*_{c',\gm}(N)
} $$
commutes.

The commutativity of the
 third square follows from Proposition~\ref{prop:gysin-phi},
and the commutativity of the first square follows from
Proposition~\ref{prop:p1cp2c}.

 For the second square,
it suffices, using Proposition~\ref{pro:cocycle-iso},  to check that
$b-m^*T_1 u$ and $-p_1^* T_1u-p_2^*T_1u+b'$ differ by a \v{C}ech
coboundary. Now  we have,
 using the relation
$\del' T_1=T_2\del-\del T_2$ (see Eq. \eqref{eq:10}) that, 
\begin{eqnarray*}
\lefteqn{(b-m^* T_1u)-(-p_1^*T_1 u -p_2^*T_1u+b')
= b-b'+\del'T_1 u}\\
& =& -T_2(e-e')+\del'T_1u
=-T_2\del u +\del'T_1u=-\del T_2 u.
\end{eqnarray*}
This concludes the proof.
\end{pf}

As an application,
we consider twisted $K$-theory group of
an inertia groupoid.  Let $ \gm\toto \gm_0$ be a
Lie groupoid and consider the
crossed module $S\gm \to \gm$.
As before, $\Lambda \gm: S\gm\rtimes \gm\toto S\gm$
denote  the inertia groupoid of $\gm$.
Any element in the image   of the transgression map
$T_1: H^3(\gm\lcom; \sS^1)\to H^2(\Lambda \gm\lcom; \sS^1)$
 is $2$-multiplicative according to Remark \ref{rmk:3.10}.
Thus one obtains a ring structure on the corresponding
 twisted $K$-theory group.
Since $H^3(\gm\lcom; \sS^1)$ classifies $2$-gerbes \cite{Breen},
we conclude that the twisted $K$-theory group on the inertia
stack twisted by   a $2$-gerbe over the stack admits
a ring structure.

\begin{them}
\label{thm:main2}
Let $\gm\toto \gm_0$ be a  proper  Lie groupoid such that
$S\gm$ is a manifold and $S\gm\to \gm_0$
is $\gm$-equivariantly $K$-oriented
(these assumptions hold for instance when $\gm$ is proper
and \'etale, or when $\gm$ is a compact connected and simply connected
Lie group). Let $d=\dim S\gm-\dim\gm_0$ (see also
Remark~\ref{rmk:non-connected}).

\begin{enumerate}
\item
For any $e\in \check{Z}^3( \gm\lcom;\sS^1)$, the 
twisted $K$-theory group $K^{*+d}_{e,\gm}(S\gm)$
 is endowed with a ring structure
$$K^{i+d}_{e,\gm}(S\gm)\otimes K^{j+d}_{e,\gm}(S\gm)
\to K^{i+j+d}_{e,\gm}(S\gm).$$
\item Assume that $e$ and  $e'\in \check{Z}^3(\gm\lcom;\sS^1)$ satisfy
$e-e'=\del u$ for some $u\in \check{C}^2(\gm\lcom;\sS^1)$.
Then there is a ring isomorphism
$$\Psi_{e',u,e}:K^*_{e,\gm}(S\gm)\to K^*_{e',\gm}(S\gm)$$
satisfying the properties:
\begin{itemize}
\item  if $e-e'=\del u$ and $e'-e''=\del u'$, then
$$\Psi_{e'',u',e'}\smalcirc\Psi_{e',u,e}=\Psi_{e'',u+u',e}$$
\item  for any $v\in  \check{C}^1(\gm\lcom;\sS^1)$,
$$\Psi_{e',u,e}=\Psi_{e',u+\del v,e}.$$
\end{itemize}
\item 
  There is a morphism
$$H^2(\gm\lcom;\sS^1)\to \mathop{Aut} K^*_{e,\gm}(S\gm).$$
The ring structure on  $K^{*+d}_{e,\gm}(S\gm)$, up to an isomorphism,
depends only  on the cohomology class
$[e]\in H^3(\gm\lcom;\sS^1)$.
The isomorphism is unique up to an automorphism of
$K^{*+d}_{e,\gm}(S\gm )$ induced   from $H^2(\gm\lcom;\sS^1)$.
\end{enumerate}
\end{them}

\begin{rmk}
In \cite{ARZ}, Adem-Ruan-Zhang introduced
 an associative stringy product for the twisted
orbifold $K$--theory of a compact, almost complex orbifold
$\mathcal X$. This product is defined on the twisted $K$--theory
 of the inertia orbifold
$\Lambda \mathcal X$, where the twisting class $\tau$ is assumed 
to be in the image of the  transgression $H^4({\mathcal X}_\com, \mathbb Z)
\to H^3(\Lambda {\mathcal X}_\com , \mathbb Z)$. It would be interesting
to  investigate the relation between the Adem-Ruan-Zhang's
stringy product and the product we introduced above.
\end{rmk}

Now assume that $G$ is a compact, connected  and simply
connected  simple Lie group. Then $SG\cong G$ and
the $G$-action on $G$ is the conjugation. 

\begin{lem}\label{lem:spinc}
$G$ admits a $G$-equivariant Spin structure,
and thus a $G$-equivariant
Spin${}^c$-structure.
\end{lem}
\begin{pf}
Identify $TG$ with $G\times  \Gg$ by left translations, where
 $\Gg$ is the Lie algebra of $G$. Under this identification,
the $G$-action on $TG$, which is the lift of  the conjugation action on $G$,
becomes  $g\cdot (x,T) =(gxg^{-1},Ad(g)T)$, $\forall g\in G$ and
$(x, T)\in G\times  \Gg$. 
Endow $\Gg$ with an ${\mathrm{Ad}}$-invariant inner product.
Since $G$ is connected and simply connected,
 the group morphism $Ad:G\to SO(\Gg)$ thus
lifts to  $G\to\widetilde{SO(\Gg)}$.
\end{pf}

According to Corollary~\ref{cor:T1iso},
the transgression map $T_1:H^3(G\lcom;\sS^1)\to
H^2((G\rtimes G)\lcom;\sS^1)$ is an isomorphism.
Thus every cohomology class in
$H^2((G\rtimes G)\lcom;\sS^1)$  can be represented by a
$2$-multiplicative cocycle. Moreover, it is known that
$H^2(G\lcom;\sS^1)=0$. Thus we are led to the following

\begin{them}\label{thm:simply-connected}
 Let $G$ be  a compact, connected  and simply
connected  simple Lie group, and $[c]\in H^2((G\rtimes G)\lcom;\sS^1)
\cong \zz$. Then the equivariant twisted $K$-theory
group $K_{[c], G}^*   (G)$ is  endowed with
 a canonical ring structure
$$K^{i+d}_{[c],G}(G)\otimes K^{j+d}_{[c],G}(G)\to K^{i+j+d}_{[c], G}(G),$$
where $d=\dim G$, in the sense that there is a canonical
isomorphism of  the rings  when  using any
two  2-cocycles in $\check{Z}^2((G\rtimes G)\lcom;\sS^1)$
which are images under  the transgression $T_1$.
\end{them}

\begin{rmk}
In general, if $G$ is  a compact and connected Lie group, 
$G$ admits a $G$-equivariant Spin${}^c$-structure if and only if
there exists an  infinitesimal character $\psi:\Gg\to \mbox{Lie}(U(1))$ such that
$\rho_G+\psi$ is a weight of $T$, where $\rho_G$ is
the half-sum  of  positive roots (for any choice of maximal torus $T$) 
 \footnote{We are grateful to Meinrenken for
pointing this out to us.}. In that case, $e\in \check{Z}^3((G\rtimes G)\lcom;\sS^1)$
determines a ring structure on $K_{[T_1 e],G}^*(G)$, but this ring structure a
priori depends on the choice of the cocycle $e$ and not just of the cohomology
class $[e]\in H^3((G\rtimes G)\lcom;\sS^1)$, since
$H^2(G\lcom;\sS^1)\ne 0$ in general.
\end{rmk}

%\begin{them}
% Let $G$ be  a compact connected   Lie group, and
% $[c]\in H^2((G\rtimes G)\lcom;\sS^1)$ be in the image
%of the transgression map $T_1:H^3(G\lcom;\sS^1)\to
%H^2((G\rtimes G)\lcom;\sS^1)$. 
%Then the equivariant twisted $K$-theory
%group $K_{[c], G}^*   (G)$ is  endowed with
% a canonical ring structure
%$$K^{i+d}_{[c],G}(G)\otimes K^{j+d}_{[c],G}(G)\to K^{i+j+d}_{[c], G}(G),$$
%where $d=\dim G$, in the sense that there is a canonical
%isomorphism of  the rings  when  using any
%two  2-cocycles in $\check{Z}^2((G\rtimes G)\lcom;\sS^1)$
%which are images under  the transgression $T_1$.
%\end{them}

%%%PING

\begin{rmk}
\begin{enumerate}
\item There is a similar product in $K$-homology. Let $G$ be a compact group,
$N\to G$ a crossed module such that $N$ is $G$-equivariantly $K$-oriented.
Let us define
$K_{i,c,G}(N)$ as $KK^i_G(A_c,\cc)$.  One can define an associative
product as the composition of the external product
$K_{i,c,G}(N)\otimes K_{j,c,G}(N)\to K_{i+j,p_1^*c+p_2^*c,G}(N)$,
the isomorphism $K_{i+j,p_1^*c+p_2^*c,G}(N)\to K_{i+j,m^*c,G}(N)$
given by $b$ such that $\del' c=\del b$, and the map $m^*:
K_{i+j,m^*c,G}(N)\to K_{i+j,c,G}(N)$ coming from
$m: A_c\to A_{m^*c}$.

It should not be hard to prove that the product in $K$-theory and
the one in $K$-homology are related by Poincar\'e duality.

\item  For the case $M=G$ and $c=0$, 
$K^*_G(G)$ is explicitly computed by
Brylinski-Zhang \cite{BZ}. It would be
interesting to investigate how to express the ring structure
 $K^i_G (G)\otimes K^j_G (G)\to K^{i+j+d}_G (G)$ explicitly
in this case.

\item
Freed-Hopkins-Teleman  have proved  a remarkable theorem \cite{FHT0, FHT1}
that the equivariant twisted
$K$-theory group  is isomorphic to  Verlinde algebra 
for a compact connected Lie group. See also \cite{CW1} for
related discussion from a different perspective.
Here we define equivariant twisted $K$-theory group from a different viewpoint
 by using $K$-theory of groupoid $C^*$-algebras. It would be interesting
to explore the connection between the
ring structure on $K_{[c], G}^* (G)$  using our   construction
and the ones in \cite{FHT0} and  \cite{CW1}.
\end{enumerate}
\end{rmk}

\end{document}